\documentclass{article}
\newtheorem{theorem}{Theorem}
\newtheorem{lemma}{Lemma}

\usepackage{graphicx}
\DeclareSymbolFont{AMSa}{U}{msa}{m}{n}
\DeclareMathDelimiter\ulcorner{\mathopen} {AMSa}{"70}{AMSa}{"70}
\DeclareMathDelimiter\urcorner{\mathclose}{AMSa}{"71}{AMSa}{"71}
\makeatletter
\def\uufill{$\m@th\mathopen\ulcorner\mkern-7mu%
  \cleaders\hbox{\rule[6pt]{1dd}{1dd}}\hfill
  \mkern-7mu\mathclose\urcorner$}
\def\overbrack#1{\vbox{\m@th\ialign{##\crcr
      \uufill\crcr\noalign{\kern-\p@\nointerlineskip}%
      $\hfil\displaystyle{#1}\hfil$\crcr}}}
\makeatother

\title{Quantum Invariant of 3-Manifolds and  Poincar\'{e} Conjecture}

\author{Sze Kui Ng
\\ {\small Department of Mathematics,
Hong Kong Baptist University, Hong Kong }\\{\small Email: szekuing@yahoo.com.hk}
}

\begin{document}
\date{}
\maketitle
\begin{abstract}
A new quantum gauge model is proposed.
From this quantum gauge model we derive a quantum invariant
of 3-manifolds. We show that this invariant
gives a classification of closed (orientable and connected)
3-manifolds.
 From this classification we prove the Poincar\'{e} Conjecture.


\end{abstract}

\section{Introduction}\label{sec00}

In 1989 Witten derived the Jones polynomial from quantum field theory based on the
Chern-Simon Lagrangian \cite{Jon}\cite{Witten}.
Inspired by Witten's work
in this paper we shall derive a way to construct quantum knots and knot invariant from a
quantum gauge model of electrodynamics  and its nonabelian
generalization.

From this quantum gauge model we derive a conformal
field theory which includes the Kac-Moody algebra
and the Knizhnik-Zamolodchikov equation \cite{Kni}. Here as a
difference from the usual conformal field theory
we can derive two quantum 
Knizhnik-Zamolodchikov (KZ) equations
which are dual to each other.
These two quantum KZ equations are equations for the product of $n$ Wilson lines $W(z,z')$ which are defined by the gauge model. 
These two quantum KZ equations can be regarded as a quantum Yang-Mills equation since it is analogous to the classical Yang-Mills equation derived from the classical Yang-Mills gauge model.

From  
the two quantum KZ
equations we derive the
skein relation of the HOMFLY
polynomial \cite{Koh}-\cite{Kau}.
In this derivation we represent the uppercrossing, zero crossing
and undercrossing of two pieces of curves by the  products of
two quantum Wilson lines $W(z_1,z_2)$ and $W(z_3,z_4)$.

Then by the braiding of quantum Wilson lines we construct generalized
Wilson loops. We show that these generalized
Wilson loops can represent knot and link diagrams and thus can be defined as quantum  knots and links. 

From quantum  knots and links we derive quantum invariant of knots and links.
We show that this  invariant gives a  classification of knots and
links. From this  invariant we then derive quantum invariant of
closed 3-manifolds. We show that this quantum invariant of closed
3-manifolds give a  classification of closed manifolds. From this
classification we then prove the Poincar\'{e} conjecture.

This paper is organized as follows. In section 2 we give a brief
description of a quantum gauge model of electrodynamics  and its
nonabelian generalization. In this paper we shall consider a
nonabelian generalization with a $SU(2)$ gauge symmetry. 
In section 3 we define the classical Wilson loop.

In section 4 
we derive the defintion of the generator of the Wilson line. From this definition we derive a conformal field theory which includes the affine
Kac-Moody algebra, the Virasoro energy operator (or the Virasoro energy-momentum tensor) and the Virasoro algebra.
In section 5 we derive the quantum KZ equation in dual form.
In section 6 we compute the solutions of the quantum KZ equation in dual form.
In section 7 we compute the quantum Wilson line.
In section 8 
we represent the braiding 
of two pieces of curves by defining the braiding of
 two quantum Wilson lines. 
In section 9 
we derive the skein relation for the HOMFLY polynomial.
In section 10  
we compute the quantum Wilson loop.
In section 11 we define generalized Wilson loops which will be as quantum knots.
In section 12  
we give some examples of generalized Wilson loops and show that they have the properties of the corresponding knot diagram and thus may be regarded as quantum knots.
In  section 13  
we show
that this generalized Wilson loop is a complete copy of the corresponding knot diagram and thus we may call a generalized Wilson loop as a quantum knot. From quantum knots we have a knot invariant of the form $Tr R^{-m} W(z,z)$ where $W(z,z)$ denotes
a Wilson loop  and $R$ is the
braiding
matrix and is the monodromy of the quantum KZ equation and $m$ is an integer. We show that this knot invariant classifies knots and that knots can be one-to-one assigned with the integer $m$.
In section 14 
we give more computations of quantum knots and their knot invariant.
Then in section 15 and 16  
with the integer $m$ we give a
classification table
of knots  where we show that 
prime knots (and only prime knots) are assigned with prime integer $m$. 

In section 17 we give examples of invariant of links.
 In section 18  
 we give a classification of links. 
 In section 19 
we construct a quantum invariant of 3-manifolds. We first construct quantum
invariant of closed three-manifolds obtained by Dehn surgery on
framed knots.  We then introduce the concept of minimal link to
construct quantum invariant of closed three-manifolds obtained by
Dehn surgery on framed links. Then by using the Lickorish-Wallace theorem which states that any closed  3-manifold $M$ can be obtained from a Dehn surgery on a
framed link $L$ we show that this invariant gives a classification of closed 3-manifolds. Then in
section  20  
we use this classification to prove the Poincar\'{e} conjecture.

\section{A Quantum Gauge Model}\label{sec2}

We shall first establish a quantum gauge model. This quantum gauge model will be as a physical motivation for introducing operators which will be called Wilson loop and Wilsom line as analogous to the Wilson loops in the existing quantum field theories. Then the definition of Wilson loop and Wilson line and the definition of a generator $J$ of the Wilson line will be as the basis of the mathematical foundation of this paper (In order to simplify the mathematics of this paper we treat this quantum gauge model as a physical motivation instead of as the mathematical foundation of this paper). 

We shall show that the generator $J$ gives an affine Kac-Moody algebra and a Virasoro energy operator $T$ with central charge $c$.
From  $J$ and $T$ we shall derive the quantum KZ equation in dual form which will be regarded as the quantum Yang-Mills equation.
From this quantum KZ equation we then construct generalized Wilson loops which will be as quantum  knots and links.

Let us  
construct a quantum gauge model, as follows. In
probability theory we have the Wiener measure $\nu$ which is a
measure on the space $C[t_0,t_1]$ of continuous functions
\cite{Jaf}. This measure is a well defined mathematical theory for
the Brownian motion and it may be symbolically written in the
following form:
\begin{equation}
d\nu =e^{-L_0}dx
\label{wiener}
\end{equation}
where $L_0 :=
\frac12\int_{t_0}^{t_1}\left(\frac{dx}{dt}\right)^2dt$ is the
energy integral of the Brownian particle and $dx =
\frac{1}{N}\prod_{t}dx(t)$ is symbolically a product of Lebesgue
measures $dx(t)$ and $N$ is a normalized constant.

Once the Wiener measure is defined we may then define other
measures on $C[t_0,t_1]$ as follows\cite{Jaf}. Let a potential
term $\frac12\int_{t_0}^{t_1}Vdt$ be added to $L_0$. Then we have
a measure $\nu_1$ on $C[t_0,t_1]$ defined by:
\begin{equation}
d\nu_1 =e^{-\frac12\int_{t_0}^{t_1}Vdt}d \nu
\label{wiener2}
\end{equation}
Under some condition on $V$ we have that $\nu_1$ is well defined
on $C[t_0,t_1]$. Let us call (\ref{wiener2}) as the Feymann-Kac
formula \cite{Jaf}.

Let us then follow this formula to construct a quantum  model of
electrodynamics, as follows. 
Then similar to the formula (\ref{wiener2}) we construct a quantum
model of electrodynamics from the following energy
integral:
\begin{equation}
 \frac12\int_{s_0}^{s_1}[
\frac12\left(\frac{\partial A_1}{\partial x^2}-\frac{\partial
A_2}{\partial x^1}\right)^* \left(\frac{\partial A_1}{\partial
x^2}-\frac{\partial A_2}{\partial x^1}\right)
 +\sum_{j=1}^2
\left(\frac{\partial Z^*}{\partial
x^j}+ieA_jZ^*\right)\left(\frac{\partial Z}{\partial
x^j}-ieA_jZ\right)]ds \label{1.1}
\end{equation}
where the complex variable $Z=Z(z(s))$ and the real variables
$A_1=A_1(z(s))$ and $A_2=A_2(z(s))$ are continuous functions in a form that they are in terms of an arbitrary (continuously differentiable) closed curve $z(s)=C(s)=(x^1(s),x^2(s)), s_0\leq
s\leq s_1, z(s_0)=z(s_1)$ in the complex plane where
$s$ is a parameter representing the proper time in relativity (We shall also write $z(s)$ in the complex variable form $C(s)=z(s)=x^1(s)+ix^2(s),s_0\leq s\leq s_1$). The complex variable $Z=Z(z(s))$ represents a field of matter( such as the electron)
($Z^*$ denotes its complex conjugate) and the real variables
$A_1=A_1(z(s))$ and $A_2=A_2(z(s))$ represent a connection (or the gauge field of the photon) and $e$
denotes the electric charge.

The integral (\ref{1.1}) has the following gauge symmetry:
\begin{equation}
\begin{array}{rl}
Z^{\prime}(z(s)) & := Z(z(s))e^{iea(z(s))} \\
A'_j(z(s)) & := A_j(z(s))+\frac{\partial a}{\partial x^j} \quad
j=1,2
\end{array}
\label{1.2}
\end{equation}
where $a=a(z)$ is a continuously differentiable real-valued
function of $z$.

We remark that this model is similar to the usual Yang-Mills gauge model. A  feature of (\ref{1.1}) is that it is not
formulated with the four-dimensional space-time but is formulated
with the one dimensional proper time. This one dimensional nature
let this model avoid the usual utraviolet divergence difficulty of
quantum fields. 

Similar to the usual Yang-Mills gauge theory we can generalize this gauge model with $U(1)$ gauge symmetry to
nonabelian gauge models. As an illustration let us consider
$SU(2)$ gauge symmetry. Similar to (\ref{1.1}) we consider the
following energy integral:
\begin{equation}
L := \frac12\int_{s_0}^{s_1}
[\frac12 Tr (D_1A_2-D_2A_1)^{*}(D_1A_2-D_2A_1) +
(D_1^*Z^*)(D_1Z)+(D_2^*Z^*)(D_2Z)]ds
\label{n1}
\end{equation}
where $Z= (z_1, z_2)^{T}$ is a two dimensional complex vector;
$A_j =\sum_{k=1}^{3}A_j^k t^k $ $(j=1,2)$ where $A_j^k$ denotes a
component of a gauge field $A^k$; $t^k$ denotes a generator of
$SU(2)$ (Here for simplicity we choose a convention that the
complex $i$ is absorbed by $t^k$); and
$D_j=\frac{\partial}{\partial x^j}-gA_j$, $(j=1,2)$ where $g$
denotes the charge of interaction (For simplicity let us set
$g=1$).

From (\ref{n1}) we can develop a nonabelian gauge model as similar
to that for the above abelian gauge model.
We have that (\ref{n1}) is invariant under the following
gauge transformation:
\begin{equation}
\begin{array}{rl}
Z^{\prime}(z(s)) & :=U(a(z(s)))Z(z(s)) \\
A_j^{\prime}(z(s)) & := U(a(z(s)))A_j(z(s))U^{-1}(a(z(s)))+
 U(a(z(s)))\frac{\partial U^{-1}}{\partial x^j}(a(z(s))),
j =1,2
\end{array}
\label{n2}
\end{equation}
where $U(a(z(s)))=e^{-a(z(s))}$ and $a(z(s))=\sum_k a^k
(z(s))t^k$. We shall mainly consider the case that $a$ is a
function of the form $a(z(s))
=\sum_k \mbox{Re}\,
\omega^k(z(s))t^k$ where $\omega^k$ are 
analytic functions of $z$ (We let $\omega(z(s)):=\sum_k\omega^k(z(s))t^k$ and we write $a(z)=\mbox{Re}\,\omega(z)$). 

The above gauge model is based on the Banach space $X$ of continuous functions $Z(z(s)), A_j(z(s)), j=1,2, s_0\leq s\leq s_1$
on the one dimensional interval $\lbrack s_0, s_1 \rbrack$. 

Since $L$ is positive and the model is one dimensional (and thus is simpler than the usual two dimensional Yang-Mills gauge model) we have that this gauge model is similar to the Wiener measure except that this gauge model has a gauge symmetry. This gauge symmetry gives a degenerate degree of freedom. In the physics literature the usual way to treat the degenerate degree of freedom of gauge symmetry is to introduce a gauge fixing condition to eliminate the degenerate degree of freedom where each gauge fixing will give equivalent physical results\cite{Fad}. There are various gauge fixing conditions such as the Lorentz gauge condition, the Feynman gauge condition, etc. We shall later in
section \ref{sec6} (on the Kac-Moody algebra) adopt a gauge fixing
condition for the above gauge model. This gauge fixing condition
will also be used to derive the quantum KZ equation in dual form which will be regarded as a quantum Yang-Mill equation since its role will be similar to the classical Yang-Mill equation derived from the classical Yang-Mill gauge model.

Since $L$ is positive we have that without gauge fixing condition the above gauge model is a positive linear functional on the Banach space $C(X)$ of continuous functions on $X$ and is multivalued in the sense that each gauge fixing gives a value.

{\bf Remark}. In this paper the main aim of introducing this quantum gauge model is to derive the quantum KZ equation in dual form which will be regarded as a quantum Yang-Mills equation (or as a quantum Euler-Lagrange equation). From this quantum KZ equation in dual form we then construct quantum knots and links. From quantum knots and links we then prove the Poincare Conjecture.

\section{Classical Wilson Loop } \label{sec4}

Similar to the Wilson loop in quantum field theory \cite{Witten} from our
quantum model we introduce an analogue of Wilson loop, as follows.

{\bf Definition}.
 A classical Wilson loop $W_R(C)$ is defined by :
\begin{equation}
W_R(C):= W(z_0, z_1):= Pe^{\int_C A_jdx^j} \label{n4}
\end{equation}
where $R$ denotes a representation of $SU(2)$; $C(\cdot)=z(\cdot)$
is a fixed curve where the quantum gauge models are based on it
as specified in the above section.
As
usual the notation $P$ in the definition of $W_R(C)$ denotes a
path-ordered product \cite{Witten}\cite{Kau}\cite{Baez}.

Let us give some remarks on the above definition of Wilson loop,
as follows.

1) We use the notation $W(z_0, z_1)$ to mean the Wilson loop
$W_R(C)$ which is based on the whole closed curve $z(\cdot)$. Here for
convenience we  use only the end points $z_0$ and $z_1$ of the
curve $z(\cdot)$ to denote this Wilson loop 
(We keep in mind that the definition of $W(z_0, z_1)$ depends on the whole curve $z(\cdot)$ connecting $z_0$ and $z_1$).

Then we  extend the definition of $W_R(C)$ to the case that
$z(\cdot)$ is not a closed curve with $z_0\neq z_1$. When
$z(\cdot)$ is not a closed curve we shall call $W(z_0, z_1)$ as a
Wilson line.
 
2) In constructing the Wilson loop we need to choose a
representation $R$ of the $SU(2)$ group. We shall see that because a Wilson line 
$W(z_0, z_1)$ is with two variables $z_0$ and $z_1$ a natural representation of a Wilson line or a Wilson loop is the tensor product of the usual two dimensional representation of  $SU(2)$ for constructing the Wilson loop. $\diamond$

We first have the following theorem on $W(z_0,z_1)$:
\begin{theorem}
For a given continuous path $A_i, i=1,2$  on $[s_0, s_1]$ 
the Wilson line $W(z_0,z_1)$ exists on this path and has the
following transition property:
\begin{equation}
W(z_0,z_1)=W(z_0,z)W(z,z_1)
 \label{df2}
\end{equation}
where $W(z_0,z_1)$ denotes the Wilson line of a
curve $z(\cdot)$ which is with $z_0$ as the starting
point and $z_1$ as the ending point and $z$ is a
point on $z(\cdot)$ between $z_0$ and $z_1$.
\end{theorem}

{\bf Proof}. We have that $W(z_0,z_1)$ is a limit
(whenever exists)
of ordered product of $e^{A_i\triangle x^i}$ and thus can be
written in the following form:
\begin{equation}
\begin{array}{rl}
W(z_0,z_1)= & I +
\int_{s^{\prime}}^{s^{\prime\prime}}
A_i(z(s))\frac{dx^i(s)}{ds}ds \\
 & + \int_{s^{\prime}}^{s^{\prime\prime}}
[\int_{s^{\prime}}^{s_1} A_i(z(s_1))\frac{dx^i(s_1)}{ds}ds_1]
A_i(z(s_2))\frac{dx^i(s_2)}{ds}ds_2 +\cdot\cdot\cdot
\end{array}
\label{df3}
\end{equation}
where $z(s^{\prime})=z_0$ and $z(s^{\prime\prime})=z_1$. Then
since $A_i$ are continuous on $[s^{\prime}, s^{\prime\prime}]$ and
$x^i(z(\cdot))$ are continuously differentiable on $[s^{\prime},
s^{\prime\prime}]$ we have that the series in (\ref{df3}) is
absolutely convergent. Thus the Wilson line $W(z_0,z_1)$ exists.
Then since $W(z_0,z_1)$ is the limit of ordered
product  we can write $W(z_0,z_1)$ in the form $W(z_0,z)W(z,z_1)$
by dividing $z(\cdot)$ into two parts at $z$. This proves the
theorem. $\diamond$

{\bf Remark (Classical and quantum Wilson loop)}. This theorm means that the Wilson line $W(z_0,z_1)$ exists in the classical pathwise sense where $A_i$ are as classical paths on $[s_0, s_1]$. This pathwise version of the Wilson line $W(z_0,z_1)$; from the Feymann path integral point of view; is as a partial description of the quantum version of the Wilson line $W(z_0,z_1)$ which is as an operator when $A_i$ are as operators.
We shall in the next section derive and define a quantum generator $J$ of $W(z_0,z_1)$ from the quantum gauge model. Then by using this generator $J$ we shall compute the quantum version of the Wilson line $W(z_0,z_1)$. 

We shall denote both the classical version and quantum version of Wilson line by the same notation  $W(z_0,z_1)$ when there is no confusion.
$\diamond$

{\bf Remark}. We remark again that in order to simplify the mathematics of this paper we treat the above quantum gauge model as a physical motivation instead of as the mathematical foundation of this paper. The mathematical foundation of this paper will base on the definition of the Wilson line $W(z_0,z_1)$ and the generator $J$ instead of the above quantum gauge model which is as a physical motivation for introducing the Wilson line $W(z_0,z_1)$ and the generator $J$. $\diamond$

By following the
usual approach 
from a gauge
transformation  we have the following symmetry on Wilson lines (This symmetry is sometimes called the chiral symmetry) \cite{Kau}:

\begin{theorem}
Under an analytic gauge
transformation with an analytic function $\omega$ 
we have the following symmetry:
\begin{equation}
W(z_0, z_1) \mapsto W^{\prime}(z_0, z_1)=U(\omega(z_1))
W(z_0, z_1)U^{-1}(\omega(z_0))
\label{n5}
\end{equation}
where $W^{\prime}(z_0, z_1)$ is a Wilson line with gauge field $A_{\mu}^{\prime} =  \frac{\partial U(z)}{\partial x^{\mu}}U^{-1}(z) + U(z)A_{\mu}U^{-1}(z)$.
\end{theorem}

{\bf Proof}. Let us prove this symmetry as follows. Let $U(z):=
U(\omega(z(s)))$ and $U(z+dz)\approx U(z)+\frac{\partial
U(z)}{\partial x^{\mu}}dx^{\mu}$ where $dz=(dx^1,dx^2)$. Following
Kauffman \cite{Kau} we have
\begin{equation}
\begin{array}{rl}
& U(z+ dz)(1+ dx^{\mu}A_{\mu})U^{-1}(z)\\
=& U(z+ dz)U^{-1}(z)
+ dx^{\mu}U(z+dz)A_{\mu}U^{-1}(s) \\
\approx & 1+ \frac{\partial U(z)}{\partial
x^{\mu}}U^{-1}(z)dx^{\mu}
  + dx^{\mu}U(z+ dz)A_{\mu}U^{-1}(s) \\
\approx & 1+ \frac{\partial U(z)}{\partial
x^{\mu}}U^{-1}(z)dx^{\mu}
+ dx^{\mu}U(z)A_{\mu}U^{-1}(z) \\
=: & 1+  \frac{\partial U(z)}{\partial x^{\mu}}U^{-1}(z)dx^{\mu}
+ dx^{\mu}U(z)A_{\mu}U^{-1}(z)\\
=:& 1 + dx^{\mu}A_{\mu}^{\prime}
\end{array}
\label{n5b}
\end{equation}

From (\ref{n5b}) we have that (\ref{n5}) holds since (\ref{n5}) is
the limit of ordered product in which the left-side factor $U(z_i+
dz_i)$ in (\ref{n5b}) with $z_i=z$ is canceled by the right-side
factor $U^{-1}(z_{i+1})$ of (\ref{n5b}) where $ z_{i+1}=z_i+ dz_i$
with $z_{i+1}=z$. This proves the theorem. $\diamond$

As analogous to  
the WZW model in
conformal field theory \cite{Fra}\cite{Fuc}
from the above symmetry  we have the following formulas for the
variations $\delta_{\omega}W$ and $\delta_{\omega^{\prime}}W$ with
respect to this symmetry:
\begin{equation}
\delta_{\omega}W(z,z')=W(z,z')\omega(z)
\label{k1}
\end{equation}
and
\begin{equation}
\delta_{\omega^{\prime}}W(z,z')=-\omega^{\prime}(z')W(z,z')
\label{k2}
\end{equation}
where $z$ and $z'$ are independent variables and
$\omega^{\prime}(z')=\omega(z)$ when $z'=z$. In (\ref{k1}) the
variation is with respect to the $z$ variable while in (\ref{k2})
the variation is with respect to the $z'$ variable. This
two-side-variations when $z\neq z'$ can be derived as follows. For
the left variation we may let $\omega$ be analytic in a
neighborhood of $z$ and continuous differentiably extended to a
neighborhood of $z'$ such that $\omega(z')=0$ in this neighborhood
of $z'$. Then from (\ref{n5}) we have that (\ref{k1}) holds.
Similarly we may let $\omega^{\prime}$ be analytic in a
neighborhood of $z'$ and continuous differentiably extended to a
neighborhood of $z$ such that $\omega^{\prime}(z)=0$ in this
neighborhood of $z$. Then we have that (\ref{k2}) holds.

\section{A Gauge Fixing Condition and Affine Kac-Moody Algebra} \label{sec6}

This section has two related purposes. One purpose is to find a
gauge fixing condition for eliminating the degenerate degree of
freedom from the gauge invariance of the above quantum gauge model
in section 2. Then another purpose is to find an equation for
defining a generator $J$ of the Wilson line $W(z,z')$. This
defining equation of $J$ can then be used as a gauge fixing
condition. Thus with this defining equation of $J$ the
construction of the quantum gauge model in section 2 is then
completed (We remark that we shall let the definition of the Wilson line and the definition
of the generator $J$ as the mathematical foundation of this paper while the quantum gauge model is as a physical motivation for deriving these two definitions).

 We shall derive a quantum loop algebra (or the
affine Kac-Moody algebra) structure from the Wilson line $W(z,z')$
for the generator $J$ of $W(z,z')$. To this end let us first
consider the classical case. Since $W(z,z')$ is constructed from $
SU(2)$ we have that the mapping $z \to W(z,z')$ (We consider
$W(z,z')$ as a function of $z$ with $z'$ being fixed) has a loop
group structure \cite{Lus}\cite{Seg}. For a loop group we have the
following generators:
\begin{equation}
J_n^a = t^a z^n \qquad n=0, \pm 1, \pm 2, ...
\label{km1}
\end{equation}
These generators satisfy the following algebra:
\begin{equation}
[J_m^a, J_n^b] =
if_{abc}J_{m+n}^c
\label{km2}
\end{equation}
This is  the so called loop algebra \cite{Lus}\cite{Seg}. Let us
then introduce the following generating function $J$:
\begin{equation}
J(w) = \sum_a J^a(w)=\sum_a j^a(w) t^a
\label{km3}
\end{equation}
where we define
\begin{equation}
J^a(w)= j^a(w) t^a :=
\sum_{n=-\infty}^{\infty}J_n^a(z) (w-z)^{-n-1}
\label{km3a}
\end{equation}

From $J$ we have
\begin{equation}
J_n^a=  \frac{1}{2\pi i}\oint_z dw (w-z)^{n}J^a(w)
\label{km4}
\end{equation}
where $\oint_z$ denotes a closed contour integral  with center $z$. This formula
can be interpreted as that
$J$ is the generator of the loop group and that
$J_n^a$ is the directional generator in the direction
$\omega^a(w)= (w-z)^n$. We may generalize $(\ref{km4})$
to the following  directional generator:
\begin{equation}
  \frac{1}{2\pi i}\oint_z dw \omega(w)J(w)
\label{km5}
\end{equation}
where the analytic function
$\omega(w)=\sum_a \omega^a(w)t^a$ is regarded
as a direction and we define
\begin{equation}
 \omega(w)J(w):= \sum_a \omega^a(w)J^a
\label{km5a}
\end{equation}

Then since $W(z,z')\in SU(2)$, from the variational formula
(\ref{km5}) for the loop algebra of the loop group of $SU(2)$ we
have that the variation of $W(z,z')$ in the direction $\omega(w)$
is given by
\begin{equation}
W(z,z')
  \frac{1}{2\pi i}\oint_z dw \omega(w)J(w)
\label{km6}
\end{equation}

Now let us consider the quantum case which is based on the quantum
gauge model in section 2. For this quantum case we shall define a
quantum generator $J$ which is analogous to the $J$ in
(\ref{km3}). We shall choose the equations (\ref{n8b}) and
(\ref{n6}) as  the equations for defining the quantum generator
$J$.  
Let us first give a
formal derivation of the equation (\ref{n8b}), as follows.
 Let us consider the
following formal functional integration:
\begin{equation}
\langle W(z,z')A(z) \rangle := \int dA_1dA_2dZ^{*}dZ  e^{-L}
W(z,z')A(z) \label{n8a}
\end{equation}
where $A(z)$ denotes a field from the quantum gauge model (We
first let $z'$ be fixed as a parameter).

Let us  do a calculus of variation on this integral to derive a variational
equation by applying a gauge transformation on (\ref{n8a}) as follows
(We remark that such variational equations are usually called the
Ward identity in the  physics literature).

Let $(A_1,A_2,Z)$ be regarded as a coordinate system of the integral
(\ref{n8a}).
Under a gauge transformation (regarded as
a change of coordinate) with gauge function $a(z(s))$ this coordinate
is changed to another coordinate denoted by
$(A_1^{\prime}, A_2^{\prime}, Z^{\prime})$.
As similar to the usual change of variable for integration we have that
the integral  (\ref{n8a}) is unchanged
under a change of variable and we have the following
equality:
\begin{equation}
\begin{array}{rl}
& \int dA_1^{\prime}
 dA_2^{\prime}dZ^{\prime *}dZ^{\prime}
 e^{-L^{\prime}} W^{\prime}(z,z')A^{\prime}(z) \\
= & \int dA_1dA_2dZ^{*}dZ  e^{-L} W(z,z')A(z)
\end{array}
\label{int}
\end{equation}
where $W^{\prime}(z,z')$ denotes the Wilson line based on
$A_1^{\prime}$ and $A_2^{\prime}$ and similarly $A^{\prime}(z)$
denotes  the field obtained from $A(z)$ with $(A_1, A_2,Z)$
replaced by $(A_1^{\prime}, A_2^{\prime},Z^{\prime})$.

Then it can be shown that the differential is unchanged under a
gauge transformation \cite{Fad}:
\begin{equation}
dA_1^{\prime}
dA_2^{\prime}dZ^{\prime *}dZ^{\prime}
= dA_1dA_2dZ^{*}dZ
\label{int2}
\end{equation}
Also by the gauge invariance property the factor $e^{-L}$ is
unchanged under a gauge transformation. Thus from (\ref{int}) we
have
\begin{equation}
0 = \langle W^{\prime}(z,z')A^{\prime}(z)\rangle -
  \langle W(z,z')A(z)\rangle
\label{w1}
\end{equation}
where the correlation notation
$\langle \rangle$ denotes the integral with
respect to the differential
\begin{equation}
e^{-L}dA_1dA_2dZ^{*}dZ
\label{w1a}
\end{equation}

We can now carry out the calculus of variation. From the gauge
transformation we have the formula
$W^{\prime}(z,z')=U(a(z))W(z,z')U^{-1}(a(z'))$ 
($a(z)=\mbox{Re}\,\omega(z)$). This
gauge transformation gives a variation of $W(z,z')$ with the
gauge function $a(z)$  
as the variational direction $a$  
in the variational formulas (\ref{km5}) and  (\ref{km6}). Thus analogous
to the variational formula (\ref{km6}) we have that the variation
of $W(z,z')$ under this gauge transformation is given by
\begin{equation}
W(z,z')
  \frac{1}{2\pi i}\oint_z dw a(w)J(w) 
\label{int3}
\end{equation}
where the generator $J$ for this variation is to
be specified. This $J$ will be a quantum generator
which generalizes the classical generator $J$ in
(\ref{km6}).

Thus under a gauge transformation with gauge function $a(z)$ from (\ref{w1}) we have the
following variational equation:
\begin{equation}
0= \langle W(z,z')[\delta_{a}A(z)+\frac{1}{2\pi i}\oint_z
dw a(w)J(w)A(z)]\rangle 
\label{w2}
\end{equation}
where $\delta_{a}A(z) $ 
denotes the variation of the field
$A(z)$ in the direction $a(z)$. 
From this equation an ansatz of
$J$ is that $J$ satisfies the following equation:
\begin{equation}
W(z,z')[\delta_{a}A(z)+\frac{1}{2\pi i}\oint_z
dw a(w)J(w)A(z)] =0 \label{n8bb}
\end{equation}
From this equation we have the following variational equation:
\begin{equation}
\delta_{a}A(z)=\frac{-1}{2\pi i}\oint_z dw a(w)J(w)A(z)
\label{n8bre}
\end{equation}
This completes the formal calculus of variation. Now (with the
above derivation as a guide) we choose the following equation (\ref{n8b}) as one of the
equation for defining the generator $J$:
\begin{equation}
\delta_{\omega}A(z)=\frac{-1}{2\pi i}\oint_z dw\omega(w)J(w)A(z)
\label{n8b}
\end{equation}
where we generalize the direction $a(z)=\mbox{Re}\,\omega(z)$ to the analytic direction $\omega(z)$
(This generalization has the effect of extending the real measure to include the complex Feymann path integral).

Let us now choose one more equation for determine the generator
$J$ in (\ref{n8b}). This choice will be as  
a gauge fixing
condition. As analogous to the WZW model in conformal field
theory \cite{Fra}\cite{Fuc} \cite{Kni}  let us consider a $J$
given by
\begin{equation}
J(z) := -k W^{-1}(z, z')\partial_z W(z, z') \label{n6}
\end{equation}
where we define $\partial_z=\partial_{x^1} +i\partial_{x^2} $ and we set $z'=z$ after the differentiation with respect to $z$;
$ k>0 $ is a constant which is fixed when the $J$ is determined to
be of the form (\ref{n6}) and the minus sign is chosen by
convention. In the WZW model \cite{Fra}\cite{Kni}
 the $J$ of the form (\ref{n6})
is the  generator  of the chiral symmetry of the WZW model. We can
write the $J$ in (\ref{n6}) in the following form:
\begin{equation}
 J(w) = \sum_a J^a(w) =
\sum_a j^a(w) t^a  
\label{km8}
\end{equation}
We see that the generators $t^a$ of $SU(2)$ appear in this form of
$J$ and  this form is analogous to the classical $J$ in
(\ref{km3}). This shows that
 this $J$ is a possible candidate for the generator
$J$ in (\ref{n8b}).

Since $W(z,z')$ is constructed by gauge field we need to have a
gauge fixing for the computations related to $W(z,z')$. Then since
the $J$ in (\ref{n8b}) and (\ref{n6}) is constructed from
$W(z,z')$ we have that in defining this $J$ as the generator $J$
of $W(z,z')$ we have chosen a condition for the gauge fixing. In
this paper we shall always choose this defining equations
(\ref{n8b}) and (\ref{n6}) for $J$ as the gauge fixing condition.

In summary we introduce the following definition.

{\bf Definition} The generator $J$ of the quantum Wilson line $W(z,z')$ whose classical version is defined by (\ref{n4}), is an operator defined by the two conditions (\ref{n8b}) and (\ref{n6}).
$\diamond$

{\bf Remark}. We remark that the condition (\ref{n6}) first defines $J$ classically. Then the condition (\ref{n8b}) raises this classical $J$ to the quantum generator $J$. $\diamond$

Now we want to show that this generator $J$ in (\ref{n8b}) and
(\ref{n6}) can be uniquely solved (This means that the gauge
fixing condition has already fixed the gauge that the degenerate
degree of freedom of gauge invariance has been eliminated so that
we can carry out computation). Before solving $J$ we give the following remark.

{\bf Remark}. We remark again that in the above of this paper we have introduced a quantum gauge model as a physical motivation for introducing the Wilson loop and Wilson line defined by (\ref{n4}) and the generator $J$ defined by the two conditions (\ref{n8b}) and (\ref{n6}). In the following of this paper all the mathematics will be based on these two definitions. 
Thus we let these two definitions be as the mathematical foundation of this paper and treat the quantum gauge model as a physical motivation for deriving these two definitions. $\diamond$

Let us now solve $J$.
From (\ref{n5}) and (\ref{n6}) we
have that the variation $\delta_{\omega}J$ of the generator $J$ in
(\ref{n6}) is given by \cite{Fra}(p.622) \cite{Kni}:
\begin{equation}
\delta_{\omega}J= \lbrack J, \omega\rbrack -k\partial_z \omega
\label{n8c}
\end{equation}

From (\ref{n8b}) and (\ref{n8c}) we have that $J$ satisfies the
following relation of current algebra
\cite{Fra}\cite{Fuc}\cite{Kni}:
\begin{equation}
J^a(w)J^b(z)=\frac{k\delta_{ab}}{(w-z)^2}
+\sum_{c}if_{abc}\frac{J^c(z)}{(w-z)} \label{n8d}
\end{equation}
where as a convention the regular term of the product
$J^a(w)J^b(z)$ is omitted. Then by following
\cite{Fra}\cite{Fuc}\cite{Kni} from (\ref{n8d}) and (\ref{km8})
we can show that the $J_n^a$ in (\ref{km3})  for the corresponding Laurent series of the quantum generator $J$ 
satisfy the following  Kac-Moody algebra:
\begin{equation}
[J_m^a, J_n^b] = if_{abc}J_{m+n}^c + km\delta_{ab}\delta_{m+n, 0}
\label{n8}
\end{equation}
where $k$ is  usually called the central extension or the level of
the Kac-Moody algebra.

{\bf Remark}. Let us also consider the other side of the chiral 
symmetry.
Similar to the $J$ in (\ref{n6}) we define a generator
$J^{\prime}$ by:
\begin{equation}
J^{\prime}(z')= k\partial_{z'}W(z, z')W^{-1}(z, z') \label{d1}
\end{equation}
where after differentiation with respect to $z'$ we set $z=z'$.
Let us then consider
 the following formal correlation:
\begin{equation}
\langle A(z')W(z,z') \rangle := \int
dA_1dA_2dZ^{*}dZ
  A(z')W(z,z')e^{-L}
\label{n8aa}
\end{equation}
where $z$ is fixed. By an approach similar to the above derivation
of (\ref{n8b}) we have the following  variational equation:
\begin{equation}
\delta_{\omega^{\prime}}A(z') =\frac{-1}{2\pi i}\oint_{z^{\prime}}
dwA(z')J^{\prime}(w) \omega^{\prime}(w) \label{n8b1}
\end{equation}
where as a gauge fixing we choose the $J^{\prime}$ in (\ref{n8b1})
be the $J^{\prime}$ in (\ref{d1}). Then similar to (\ref{n8c}) we
also have
\begin{equation}
\delta_{\omega^{\prime}}J^{\prime}= \lbrack  J^{\prime},
\omega^{\prime}\rbrack -k\partial_{z'} \omega^{\prime}
\label{n8c1}
\end{equation}
Then from (\ref{n8b1}) and (\ref{n8c1}) we can derive the current
algebra and the Kac-Moody algebra for $J^{\prime}$ which are of
the same form of (\ref{n8d}) and (\ref{n8}). From this we  have
$J^{\prime}=J$. $\diamond$


\section{Quantum Knizhnik-Zamolodchikov Equation In Dual Form} \label{sec7}

With the above current algebra $J$ and the formula (\ref{n8b}) we can now
follow the usual approach
in conformal field theory to derive a quantum
Knizhnik-Zamolodchikov (KZ) equation for the product of
primary fields in a conformal field theory \cite{Fra}\cite{Fuc}\cite{Kni}.
We shall derive the KZ equation for the product of $n$ Wilson lines $W(z, z')$.
Here an important point is that from the two sides of
$W(z, z')$  we can derive two quantum KZ equations which are
dual to each other. These two quantum KZ equations are different from the usual KZ equation in that they are equations for the quantum operators $W(z, z')$ while the usual KZ equation is for the correlations of quantum operators.

With this difference the following derivation of  KZ equation for deriving these two quantum KZ equations is well known  in conformal field theory  \cite{Fra}\cite{Fuc}. The reader may skip this derivation of  KZ equation and just look at the form of the Virasoro energy operator $T(z)$ (which is usually called the Virasoro energy-momentum tensor) and the Virasoro algebra and the form of these two quantum KZ equations.

Let us first consider (\ref{k1}).
From (\ref{n8b}) and (\ref{k1}) we have
\begin{equation}
J^a(z)W(w, w') = \frac{-t^aW(w,w')}{z-w}
\label{k3}
\end{equation}
where as a convention the regular term of the product $J^a(z)W(w, w')$
is omitted.

Following \cite{Fra} and \cite{Fuc}
let us define an energy operator $T(z)$ by
\begin{equation}
T(z) := \frac{1}{2(k+g)}\sum_a :J^a(z)J^a(z):
\label{k4}
\end{equation}
where $g$ is the dual Coxter number of $SU(2)$ \cite{Fra}. In (\ref{k4})
the symbol $:J^a(z)J^a(z):$ denotes the normal ordering  of the operator $J^a(z)J^a(z)$ which can be
defined as follows \cite{Fra}\cite{Fuc}. Let a product of operators
$A(z)B(w)$ be written in the following Laurent series
form:
\begin{equation}
A(z)B(w)= \sum_{n=-n_0}^{\infty}
 a_n(w)(z-w)^n
\label{normal}
\end{equation}
The singular part
of (\ref{normal}) is called the contraction
of $A(z)B(w)$ and will be denoted by $\overbrack{A(z)B}(w)$.
Then the term $a_0(w)$ is called the normal ordering
of $A(z)B(w)$ and we denote $a_0(w)$ by $:A(w)B(w):$.
These terms are originally from quantum field theory.
We remark that in \cite{Fra} the notation $(AB)$ is
used to generalize the original definition of
$:AB:$ for products of free fields. Here for simplicity
we shall always use the notation $:AB:$ to mean the
normal ordering of $AB$. From this definition of normal ordering we have the following
form of normal ordering \cite{Fra}:
\begin{equation}
:A(w)B(w): = \frac{1}{2\pi i}\oint \frac{dz}{z-w}A(z)B(w)
 \label{conform1a}
\end{equation}
This form can be checked by taking the contour integral on the
Laurent series expansion of $A(z)B(w)$. Alternatively we may let
(\ref{conform1a}) be the definition of normal ordering.
We  then define  (\ref{k4}) by (\ref{conform1a}) with $A=B=J^a$.

The above definition of the energy operator $T(z)$
is called the
Sugawara construction \cite{Fra}.
We first have the following well known theorem on $T(z)$ in conformal field theory \cite{Fra}:
\begin{theorem}
The operator product $T(z)T(w)$ is given by the following
formula:
\begin{equation}
T(z)T(w)=\frac{c}{2(z-w)^4}+
         \frac{2T(w)}{(z-w)^2}+\frac{\partial T(w)}{(z-w)}
\label{k5}
\end{equation}
for some constant $c=\frac{4k}{k+g}$ ($g=2$ for the group $SU(2)$) and as a convention we omit the regular term of this product.
\end{theorem}
{\bf Proof}. In \cite{Fra} there is a detail proof of this theorem. Here we want to remark that the formula (\ref{n8d}) for the product $J^a(z)J^b(x)$ is used for the proof of this theorem. $\diamond$

From this theorem we then have the following Virasoro algebra
of the mode expansion of $T(z)$ \cite{Fra}\cite{Fuc}:
\begin{theorem}
Let us write $T(z)$ in the following Laurent series form:
\begin{equation}
T(z)= \sum_{n=-\infty}^{\infty}(z-w)^{-n-2}L_n(w)
\label{conform9}
\end{equation}
This means that the modes $L_n(w)$ are defined by
\begin{equation}
L_n(w):= \frac{1}{2\pi i}\oint_w dz (z-w)^{n+1}T(z)
\label{conform10}
\end{equation}
Then we have that $L_n$ form a Virasoro algebra:
 \begin{equation}
[L_n, L_m] = (n-m)L_{n+m} + \frac{c}{12}n(n^2-1)\delta_{n+m,0}
\label{conform11}
\end{equation}
\end{theorem}

From  the formula (\ref{n8d}) for the product $J^a(z)J^b(w)$ we have the following operator product expansion \cite{Fra}:
\begin{equation}
\overbrack{T(z) J^a}(w)
= \frac{J^a(W)}{(z-w)^2}+ \frac{\partial J^a(W)}{(z-w)}
\label{conform7}
\end{equation}

Then we have the following operator product of $T(z)$ with an operator $A(w)$:
 \begin{equation}
T(z)A(w) = \sum_{n=-\infty}^{\infty}(z-w)^{-n-2}L_nA(w)
\label{conform13}
\end{equation}

From (\ref{conform7})  and (\ref{conform13}) we have that $L_{-1}J^a(w)= \partial J^a(w)$ and $L_{-1} =  \frac{\partial}{\partial z}$.
Thus we have
\begin{equation}
L_{-1}W(w, w')= \frac{\partial W(w, w')}{\partial w}
\label{conform14}
\end{equation}

On the other hand as shown in \cite{Fra} by using the Laurent series expansion of
$J^a(z)$ in the section on Kac-Moody algebra we can compute
the normal ordering $:J^a(z)J^a(z):$ from which we have the
Laurent series expansion of $T(z)$ with $L_{-1}$ given by \cite{Fra}:
\begin{equation}
L_{-1}=\frac{1}{2(k+g)}\sum_a \lbrack
\sum_{m\leq -1}J_m^aJ_{-1-m}^a +
\sum_{m\geq 0}J_{-1-m}^a J_m^a \rbrack
\label{conform15}
\end{equation}
where since $J_m^a$ and $J_{-1-m}^a$ commute each other the ordering of them
is irrelevant.

From (\ref{conform15}) we then have
\begin{equation}
\begin{array}{rl}
 & L_{-1}W(w,w') \\
=& \frac{1}{2(k+g)}\sum_a \lbrack
\sum_{m\leq -1}J_m^a(w)J_{-1-m}^a(w) +
\sum_{m\geq 0}J_{-1-m}^a(w) J_m^a(w) \rbrack W(w,w') \\
 = & \frac{1}{(k+g)}J_{-1}^a(w) J_{0}^a(w) W(w,w')
\end{array}
\label{conform16}
\end{equation}
since $J_m^a W(w,w')=0$ for $m>0$.

It follows from (\ref{conform14}) and (\ref{conform16}) that we have
the following equality:
\begin{equation}
\partial_w W(w, w')  = \frac{1}{(k+g)}J_{-1}^a(w) J_{0}^a(w) W(w,w')
\label{k7}
\end{equation}

Then form (\ref{k3}) we have
\begin{equation}
J_{0}^a(w)W(w,w')=-t^aW(w,w')
\label{k8}
\end{equation}
From (\ref{k7}) and (\ref{k8}) we then have
\begin{equation}
\partial_z W(z, z')=\frac{-1}{k+g}J_{-1}^a(z)t^aW(z,z')
\label{k9}
\end{equation}

Now let us
consider a product of $n$ Wilson lines:
$ W(z_1, z_1^{\prime})\cdot\cdot\cdot
 W(z_n, z_n^{\prime})$.
Let this product be represented as a tensor product when
$z_i$ and $z_j^{\prime}$, $i,j=1,...,n$ are all independent variables.
 Then from (\ref{k9}) we have
\begin{equation}
\begin{array}{rl}
& \partial_{z_{i}} W(z_1, z_1^{\prime})\cdot\cdot\cdot
W(z_i, z_i^{\prime})\cdot\cdot\cdot
 W(z_n, z_n^{\prime})\\
=& \frac{-1}{k+g}W(z_1, z_1^{\prime})\cdot\cdot\cdot
J_{-1}^a(z_i)t^aW(z_i,z_i^{\prime})\cdot\cdot\cdot
W(z_n, z_n^{\prime})\\
=& \frac{-1}{k+g}J_{-1}^a(z_i)t^a W(z_1, z_1^{\prime})\cdot\cdot\cdot
W(z_i,z_i^{\prime})\cdot\cdot\cdot
W(z_n, z_n^{\prime})
\end{array}
\label{k9a}
\end{equation}
where the second equality is from the definition
of tensor product for which we define
\begin{equation}
t^a W(z_1, z_1^{\prime})\cdot\cdot\cdot
W(z_i,z_i^{\prime})\cdot\cdot\cdot
W(z_n, z_n^{\prime})
:=W(z_1, z_1^{\prime})\cdot\cdot\cdot
[t^a W(z_i,z_i^{\prime})]\cdot\cdot\cdot
W(z_n, z_n^{\prime})
\label{k9aa}
\end{equation}

With this formula (\ref{k9a}) we can now follow \cite{Fra} and \cite{Fuc} to
derive the KZ equation.
For a easy reference let us present this derivation in \cite{Fra} and \cite{Fuc} as follows.
From the Laurent series of $J^a$ we have
\begin{equation}
J_{-1}^a(z_i) = \frac{1}{2\pi i}
\oint_{z_i} \frac{dz}{z-z_i}J^a(z)
\label{norm2}
\end{equation}
where the line integral is on a contour encircling  $z_i$. We also
let this contour encircles all other $z_j$ so that the effects
from Wilson lines $W(z_j, z_j^{\prime})$ for $j=1,...,n$ will all be counted.
Then we have
\begin{equation}
\begin{array}{rl}
 & J_{-1}^a(z_i)W(z_1, z_1^{\prime})\cdot\cdot\cdot
 W(z_n, z_n^{\prime}) \\
& \\
= & \frac{1}{2\pi i}\oint_{z_i} \frac{dz}{z-z_i}
 J^a(z)W(z_1, z_1^{\prime})\cdot\cdot\cdot
W(z_n, z_n^{\prime})
 \\
& \\
= & \frac{1}{2\pi i}\oint_{z_i} \frac{dz}{z-z_i}
 \sum_{j=1}^n W(z_1, z_1^{\prime})\cdot\cdot\cdot
 [\frac{-t^a}{z-z_j}W(z_j,z_j^{\prime})]\cdot\cdot\cdot
W(z_n, z_n^{\prime})
\end{array}
\label{norm3aa}
\end{equation}
where the second equality is from the $JW$ product formula (\ref{k3}).
Then by a deformation of the contour integral in (\ref{norm3aa})
into a sum of $n$ contour integrals such that each contour integral
encircles one and only one $z_j$ we have:
\begin{equation}
\begin{array}{rl}
& \sum_{j=1}^n
\frac{1}{2\pi i}\oint_{z_j} \frac{dz}{z-z_i}
\sum_{k=1}^n W(z_1, z_1^{\prime})\cdot\cdot\cdot
 [\frac{-t^a}{z-z_k}W(z_k,z_k^{\prime})]\cdot\cdot\cdot
W(z_n, z_n^{\prime}) \\
& \\
=& \sum_{j=1,j\neq i}^n
\frac{1}{z_j-z_i}W(z_1, z_1^{\prime})\cdot\cdot\cdot
 [-t^a W(z_j,z_j^{\prime})]\cdot\cdot\cdot
W(z_n, z_n^{\prime})\\
& \\
=& \sum_{j=1,j\neq i}^n
\frac{t_j^a}{z_i-z_j}W(z_1, z_1^{\prime})\cdot\cdot\cdot
 W(z_n, z_n^{\prime})
 \end{array}
\label{norm3}
\end{equation}
where for the second equality we have used the definition
of tensor product.
From (\ref{norm3})
 and by applying (\ref{k9a})
to $z_i$ for $i=1,...,n$ we have the following 
Knizhnik-Zamolodchikov equation \cite{Fra} \cite{Fuc}\cite{Kni}:
\begin{equation}
\partial_{z_i}
 W(z_1, z_1^{\prime})\cdot\cdot\cdot
W(z_n, z_n^{\prime})
=\frac{-1}{k+g}
\sum_{j\neq i}^{n}
\frac{\sum_a t_i^a \otimes t_j^a}{z_i-z_j}
 W(z_1, z_1^{\prime})\cdot\cdot\cdot
W(z_n, z_n^{\prime})
\label{n9}
\end{equation}
for $i=1, ..., n$. We remark that in (\ref{n9}) we have
defined $t_i^a:= t^a$ and
\begin{equation}
\begin{array}{rl}
 & t_i^a \otimes t_j^a W(z_1, z_1^{\prime})\cdot\cdot\cdot
W(z_n, z_n^{\prime}) \\
:=& W(z_1, z_1^{\prime})\cdot\cdot\cdot
 [t^aW(z_i, z_i^{\prime})]\cdot\cdot\cdot
[t^aW(z_j, z_j^{\prime})]\cdot\cdot\cdot
W(z_n, z_n^{\prime})
\end{array}
\label{n9a}
\end{equation}

 It is interesting and important that we also have
another KZ equation with respect to the
$z_i^{\prime}$ variables. The derivation of this  KZ equation is
dual to the above derivation in that the operator products and
their corresponding variables are with reverse order to that in
the above derivation.

From (\ref{k2}) and (\ref{n8b1}) we have a $WJ^{\prime}$ operator product given
by
\begin{equation}
W(w, w')J^{\prime a}(z') = \frac{-W(w, w')t^a}{w'-z'}
\label{d2}
\end{equation}
where we have omitted the regular term of the product.
Then similar to the above derivation of the KZ equation from (\ref{d2})
we can then derive the following Knizhnik-Zamolodchikov
equation which is dual to (\ref{n9}):
\begin{equation}
\partial_{z_i^{\prime}}
 W(z_1,z_1^{\prime})\cdot\cdot\cdot W(z_n,z_n^{\prime})
= \frac{-1}{k+g}\sum_{j\neq i}^{n}
 W(z_1, z_1^{\prime})\cdot\cdot\cdot
W(z_n, z_n^{\prime})
\frac{\sum_a t_i^a\otimes t_j^a}{z_j^{\prime}-z_i^{\prime}}
\label{d8}
\end{equation}
for $i=1, ..., n$ where we have defined:
\begin{equation}
\begin{array}{rl}
 &  W(z_1, z_1^{\prime})\cdot\cdot\cdot
W(z_n, z_n^{\prime})t_i^a \otimes t_j^a \\
:=& W(z_1, z_1^{\prime})\cdot\cdot\cdot
 [W(z_i, z_i^{\prime})t^a]\cdot\cdot\cdot
[W(z_j, z_j^{\prime})t^a]\cdot\cdot\cdot
W(z_n, z_n^{\prime})
\end{array}
\label{d8a}
\end{equation}

{\bf Rematk}. From the generator $J$ and the Kac-Moody algebra
we have derived a quantum KZ equation in dual form.
This quantum KZ equation in dual form may be consider as a quantum Yang-Mills equation since it is analogous to the classical Yang-Mills equation which is derived from the classical Yang-Mills gauge model. This quantum KZ equation in dual form will be as the starting point for the construction of quantum knots and links. $\diamond$

\section{Solving Quantum KZ Equation In Dual Form}\label{sec8a}

Let us consider the following product of two
quantum Wilson lines:
\begin{equation}
G(z_1,z_2, z_3, z_4):=
 W(z_1, z_2)W(z_3, z_4)
\label{m1}
\end{equation}
where the two quantum Wilson lines $W(z_1, z_2)$ and
$W(z_3, z_4)$ represent two pieces
of curves starting at $z_1$ and $z_3$ and ending at
$z_2$ and $z_4$ respectively.

We have that this product $G$ satisfies the KZ equation for the
variables $z_1$, $z_3$ and satisfies the dual KZ equation
for the variables $z_2$ and $z_4$.
Then
by solving the two-variables-KZ equation in (\ref{n9}) we have that a form of $G$ is
given by \cite{Koh}\cite{Dri}\cite{Chari}:
\begin{equation}
e^{-t\log [\pm (z_1-z_3)]}C_1
\label{m2}
\end{equation}
where $t:=\frac{1}{k+g}\sum_a t^a \otimes t^a$
and $C_1$ denotes a constant matrix which is independent
of the variable $z_1-z_3$.

We see that $G$ is a multivalued analytic function
where the determination of the $\pm$ sign depended on the choice of the
branch.

Similarly by solving the dual two-variable-KZ equation
 in (\ref{d8}) we have that
$G$ is of the form
\begin{equation}
C_2e^{t\log [\pm (z_4-z_2)]}
\label{m3}
\end{equation}
where $C_2$ denotes a constant matrix which is independent
of the variable $z_4-z_2$.

From (\ref{m2}), (\ref{m3}) and we let
$C_1=Ae^{t\log[\pm (z_4-z_2)]}$,
$C_2= e^{-t\log[\pm (z_1-z_3)]}A$ where $A$ is a constant matrix  we have that
$G$ is given by
\begin{equation}
G(z_1, z_2, z_3, z_4)=
e^{-t\log [\pm (z_1-z_3)]}Ae^{t\log [\pm (z_4-z_2)]}
\label{m4}
\end{equation}
where at the singular case that $z_1=z_3$ we simply define $\log [\pm (z_1-z_3)]=0$. Similarly for $z_2=z_4$.

Let us find a form of the initial operator $A$. We notice that there are two operators $\Phi_{\pm}(z_1-z_2):=e^{-t\log [\pm (z_1-z_3)]}$ and $\Psi_{\pm}(z_i^{\prime}-z_j^{\prime})$ acting on the two sides of $A$ respectively where  the two independent variables  $z_1, z_3$ of $\Phi_{\pm}$ are mixedly from the two quantum Wilson lines $W(z_1, z_2)$ and
$W(z_3, z_4)$ respectively and the the two independent variables  $z_2, z_4$ of $\Psi_{\pm}$ are mixedly from the two  quantum Wilson lines $W(z_1, z_2)$ and $W(z_3, z_4)$ respectively. From this we determine the form of $A$ as follows.

Let $D$ denote a representation of $SU(2)$. Let $D(g)$ represent an element $g$ of  $SU(2)$
and let $D(g)\otimes D(g)$ denote the tensor product representation of $SU(2)$. Then
in the KZ equation we define
\begin{equation}
[t^a\otimes t^a] [D(g_1)\otimes D(g_1)]\otimes
[D(g_2)\otimes D(g_2)]
:=[t^aD(g_1)\otimes D(g_1)]\otimes
[t^aD(g_2)\otimes D(g_2)]
\label{tensorproduct}
\end{equation}
and
\begin{equation}
[D(g_1)\otimes D(g_1)]\otimes
[D(g_2)\otimes D(g_2)][t^a\otimes t^a]
:=[D(g_1)\otimes D(g_1)t^a]\otimes
[D(g_2)\otimes D(g_2)t^a]
\label{tensorproduct2}
\end{equation}

Then we let $U({\bf a})$  
denote the universal
enveloping algebra 
where ${\bf a}$ denotes an algebra which is formed by the Lie
algebra $su(2)$ and the identity matrix. 

Now let the initial operator $A$ be of the form $A_1\otimes A_2\otimes A_3\otimes A_4$ with $A_i,i=1,...,4$
taking values in $U({\bf a})$. 
In
this case we have that in (\ref{m4}) the operator
$\Phi_{\pm}(z_1-z_2):=e^{-t\log [\pm (z_1-z_3)]}$ acts on $A$ from
the left via the following formula:
\begin{equation}
t^a\otimes t^a A=
[t^a A_1]\otimes A_2\otimes [t^a A_3]\otimes A_4
\label{ini2}
\end{equation}

Similarly the operator
$\Psi_{\pm}(z_1-z_2):=e^{t\log [\pm (z_1-z_3)]}$
in (\ref{m4}) acts on $A$ from the right via the following formula:
\begin{equation}
A t^a\otimes t^a =
A_1\otimes [A_2 t^a]\otimes A_3\otimes[A_4 t^a]
\label{ini3}
\end{equation}

We may generalize the above tensor product of
two quantum Wilson lines as follows.
Let us consider a tensor product of $n$ quantum Wilson lines:
$W(z_1, z_1^{\prime})\cdot\cdot\cdot W(z_n, z_n^{\prime})$
where the variables $z_i$, $z_i^{\prime}$
are all independent. By solving the two KZ equations
we have that this tensor product is given by:
\begin{equation}
W(z_1, z_1^{\prime})\cdot\cdot\cdot W(z_n, z_n^{\prime})
=\prod_{ij} \Phi_{\pm}(z_i-z_j)
A\prod_{ij}
\Psi_{\pm}(z_i^{\prime}-z_j^{\prime})
\label{tensor}
\end{equation}
where $\prod_{ij}$ denotes a product of
$\Phi_{\pm}(z_i-z_j)$ or
$\Psi_{\pm}(z_i^{\prime}-z_j^{\prime})$
for $i,j=1,...,n$ where $i\neq j$.
In (\ref{tensor}) the initial operator
$A$ is represented as a tensor product of operators $A_{iji^{\prime}j^{\prime}}, i,j,i^{\prime}, j^{\prime}=1,...,n$ where each $A_{iji^{\prime}j^{\prime}}$ is of the form of the initial operator $A$ in the above tensor product of two-Wilson-lines case and is acted by $\Phi_{\pm}(z_i-z_j)$ or
$\Psi_{\pm}(z_i^{\prime}-z_j^{\prime})$ on its two sides respectively.

\section{Computation of Quantum Wilson Lines}\label{sec 8aa}

Let us consider the following product of two quantum
Wilson lines:
\begin{equation}
G(z_1,z_2, z_3, z_4):=
W(z_1, z_2)W(z_3, z_4)
\label{h1}
\end{equation}
where the two quantum Wilson lines $W(z_1, z_2)$ and
$W(z_3, z_4)$ represent two pieces
of curves starting at $z_1$ and $z_3$ and ending at
$z_2$ and $z_4$ respectively.
As shown in the above section we have that $G$
is given by the following formula:
\begin{equation}
G(z_1, z_2, z_3, z_4)=
e^{-t\log [\pm (z_1-z_3)]}Ae^{t\log [\pm (z_4-z_2)]}
\label{m4a}
\end{equation}
where the product is  
a 4-tensor.

Let us set $z_2=z_3$. Then  
the 4-tensor $W(z_1, z_2)W(z_3, z_4)$ is reduced to the 2-tensor $W(z_1, z_2)W(z_2, z_4)$. By using (\ref{m4a}) the 2-tensor $W(z_1, z_2)W(z_2, z_4)$ is given by:
\begin{equation}
W(z_1, z_2)W(z_2, z_4)
=e^{-t\log [\pm (z_1-z_2)]}A_{14}e^{t\log [\pm (z_4-z_2)]}
\label{closed1}
\end{equation}
where $A_{14}=A_1\otimes A_4$ is a 2-tensor reduced from the 4-tensor 
$A=A_1\otimes A_2\otimes A_3\otimes A_4$ in (\ref{m4a}). In this reduction the $t$ operator of $\Phi=e^{-t\log [\pm (z_1-z_2)]}$ acting on the left side of $A_1$ and $A_3$ in $A$ is reduced to acting on the left side of $A_1$ and $A_4$ in $A_{14}$. Similarly the $t$ operator of $\Psi=e^{-t\log [\pm (z_4-z_2)]}$ acting on the right side of $A_2$ and $A_4$ in $A$ is reduced to acting on the right side of $A_1$ and $A_4$ in $A_{14}$.

Then since $t$ is a 2-tensor operator we have that $t$ is as a matrix acting on the two sides of the 2-tensor $A_{14}$ which is also as a matrix with the same dimension as $t$.
Thus $\Phi$ and $\Psi$ are as matrices of the same dimension as the matrix
$A_{14}$  acting on $A_{14}$ by the usual matrix operation.
Then since $t$ is a Casimir operator for the 2-tensor group representation of $SU(2)$ we have that
$\Phi $  and $\Psi $ commute
with $A_{14}$ since  $\Phi $  and $\Psi$ are exponentials
of $t$ (We remark that $\Phi $  and $\Psi $ are in general not commute with the 4-tensor initial operator $A$).
Thus we have
\begin{equation}
e^{-t\log [\pm (z_1-z_2)]}A_{14}e^{t\log[\pm (z_4-z_2)]}
=e^{-t\log [\pm (z_1-z_2)]}e^{t\log[\pm (z_4-z_2)]}A_{14}
\label{closed1a}
\end{equation}

We let $W(z_1, z_2)W(z_2, z_4)$ be as a representation of the quantum Wilson line $W(z_1,z_4)$ and we write $W(z_1,z_4)=W(z_1, z_2)W(z_2, z_4)$.
Then we have the following representation of $W(z_1,z_4)$:
\begin{equation}
W(z_1,z_4)=W(z_1,w_1)W(w_1,z_4)=e^{-t\log [\pm (z_1-w_1)]}e^{t\log[\pm (z_4-w_1)]}A_{14}
\label{closed1a1}
\end{equation}
This representation of the quantum Wilson line $W(z_1,z_4)$ means that the line (or path) with end points $z_1$ and $z_4$ is specified that it passes the intermediate point $w_1=z_2$. This representation shows the quantum nature that the path is  not specified at other intermediate points except the intermediate point $w_1=z_2$. This unspecification of the path is of the same quantum nature of the Feymann path description of quantum mechanics.

Then let us consider another representation of the quantum Wilson line $W(z_1,z_4)$. We consider $W(z_1,w_1)W(w_1,w_2)W(w_2,z_4)$ which is obtained from the tensor $W(z_1,w_1)W(u_1,w_2)W(u_2,z_4)$ by two reductions where $z_j$, $w_j$, $u_j$, $j=1,2$ are independent variables. For this representation we have: 
\begin{equation}
W(z_1,w_1)W(w_1,w_2)W(w_2,z_4)
=e^{-t\log [\pm (z_1-w_1)]}e^{-t\log [\pm (z_1-w_2)]}
e^{t\log[\pm (z_4-w_1)]}e^{t\log[\pm (z_4-w_2)]}A_{14}
\label{closed1a2}
\end{equation}
This representation of the quantum Wilson line $W(z_1,z_4)$ means that the line (or path) with end points $z_1$ and $z_4$ is specified that it passes the intermediate points $w_1$ and $w_2$. This representation shows the quantum nature that the path is  not specified at other intermediate points except the intermediate points $w_1$ and $w_2$. This unspecification of the path is of the same quantum nature of the Feymann path description of quantum mechanics.

Similarly we may represent the quantum Wilson line $W(z_1,z_4)$ by path with end points $z_1$ and $z_4$ and is specified only to pass at finitely many intermediate points. Then we let the quantum Wilson line $W(z_1,z_4)$ as an equivalent class of all these representations. Thus we may write $W(z_1,z_4)=W(z_1,w_1)W(w_1,z_4)=W(z_1,w_1)W(w_1,w_2)W(w_2,z_4)=\cdot\cdot\cdot$.

{\bf Remark}. Since $A_{14}$ is a 2-tensor  
we have that a natural group representation for the Wilson line $W(z_1,z_4)$ is the 2-tensor group representation of the group $SU(2)$.

\section{Representing Braiding of Curves by Quantum Wilson Lines}\label{sec 9aa}

Consider again the product $G(z_1, z_2, z_3, z_4)=W(z_1,z_2)W(z_3,z_4)$.
We have that $G$ is a multivalued analytic function
where the determination of the $\pm$ sign depended on the choice of the
branch.

Let the two pieces of curves be crossing at $w$. Then we have $W(z_1,z_2)=W(z_1,w)W(w,z_2)$ and
 $W(z_3,z_4)=W(z_3,w)W(w,z_4)$. Thus we have
\begin{equation}
W(z_1,z_2)W(z_3,z_4)=
W(z_1,w)W(w,z_2)W(z_3,w)W(w,z_4)
\label{h2}
\end{equation}

If we interchange $z_1$ and $z_3$, then from
(\ref{h2}) we have the following ordering:
\begin{equation}
 W(z_3,w)W(w, z_2)W(z_1,w)W(w,z_4)
\label{h3}
\end{equation}

Now let us choose a  branch. Suppose that
these two curves are cut from a knot and that
following the orientation of a knot the
curve represented by  $W(z_1,z_2)$ is before the
curve represented by  $W(z_3,z_4)$. Then we fix a branch such that the  product in (\ref{m4a}) is
with two positive signs :
\begin{equation}
W(z_1,z_2)W(z_3,z_4)=
e^{-t\log(z_1-z_3)}Ae^{t\log(z_4-z_2)}
\label{h4}
\end{equation}

Then if we interchange $z_1$ and $z_3$ we have
\begin{equation}
W(z_3,w)W(w, z_2)W(z_1,w)W(w,z_4) =
e^{-t\log[-(z_1-z_3)]}Ae^{t\log(z_4-z_2)}
\label{h5}
\end{equation}
From (\ref{h4}) and (\ref{h5}) as a choice of branch we have
\begin{equation}
W(z_3,w)W(w, z_2)W(z_1,w)W(w,z_4) =
R W(z_1,w)W(w,z_2)W(z_3,w)W(w,z_4)
\label{m7a}
\end{equation}
where $R=e^{-i\pi t}$ is the monodromy of the KZ equation.
In (\ref{m7a}) $z_1$ and $z_3$ denote two points on a closed curve
such that along the direction of the curve the point
$z_1$ is before the point $z_3$ and in this case we choose
a branch such that the angle of $z_3-z_1$ minus the angle
of $z_1-z_3$ is equal to $\pi$.

{\bf Remark}. We may use other representations of $W(z_1,z_2)W(z_3,z_4)$. For example we may use the following representation:
\begin{equation}
\begin{array}{rl}
 &W(z_1,w)W(w, z_2)W(z_3,w)W(w,z_4)\\
= &e^{-t\log(z_1-z_3)}e^{-2t\log(z_1-w)}e^{-t2\log(z_3-w)}Ae^{t\log(z_4-z_2)}e^{2t\log(z_4-w)}e^{2t\log(z_2-w)}
\end{array}
\label{h4a}
\end{equation}
Then the interchange of $z_1$ and $z_3$ changes only $z_1-z_3$ to $z_3-z_1$. Thus the formula (\ref{m7a}) holds. Similarly
all other representations of $W(z_1,z_2)W(z_3,z_4)$ will give the same result. $\diamond$

Now from (\ref{m7a}) we can take a convention that the ordering (\ref{h3}) represents that
the curve represented by  $W(z_1,z_2)$ is upcrossing
the curve represented by  $W(z_3,z_4)$ while
(\ref{h2}) represents zero crossing of these two
curves.

Similarly from the dual KZ equation as a choice of branch which
is consistent with the above formula we have
\begin{equation}
W(z_1,w)W(w,z_4)W(z_3,w)W(w,z_2)=
W(z_1,w)W(w,z_2)W(z_3,w)W(w,z_4)R^{-1}
\label{m8a}
\end{equation}
where $z_2$ is before $z_4$. We take a convention that the ordering (\ref{m8a}) represents that
the curve represented by $W(z_1,z_2)$ is undercrossing the curve represented by $W(z_3,z_4)$.
Here along the orientation of a closed curve the piece of curve
represented by $W(z_1,z_2)$ is before the piece of curve represented by
$W(z_3,z_4)$. In this case since the angle of $z_3-z_1$ minus the angle
of $z_1-z_3$ is equal to $\pi$ we have that the
angle of $z_4-z_2$ minus the angle of $z_2-z_4$ is
also equal to $\pi$ and this gives the $R^{-1}$ in this formula
(\ref{m8a}).

From (\ref{m7a}) and (\ref{m8a}) we have
\begin{equation}
 W(z_3,z_4)W(z_1,z_2)=
RW(z_1,z_2)W(z_3,z_4)R^{-1} \label{m9}
\end{equation}
where $z_1$ and $z_2$ denote the end points of a curve which is before a curve with end points $z_3$ and $z_4$.
From (\ref{m9}) we see that the algebraic structure of these
quantum Wilson lines $W(z,z')$
is analogous to the quasi-triangular quantum
group \cite{Fuc}\cite{Chari}.

\section{Skein Relation for the HOMFLY Polynomial }\label{sec9a}

In this section let us apply the above result which is from the KZ equation in dual form
to derive the skein relation for the HOMFLY polynomial. 
From this relation we then have the skein relation for the Jones polynomial which
is a special case of the HOMFLY polynomial \cite{Lic}\cite{Mur}\cite{Kau}.  

It is well known that from the one-side KZ equation
we can derive a braid group representation which is
related to the the derivation of the
skein relation for the Jones polynomial \cite{Koh}\cite{Kan}\cite{Dri}. We shall see
here that by applying the two KZ equations of the KZ equation in dual form
we also have a  way to derive the
skein relation for the HOMFLY polynomial.

Let us first consider the following theorem
of Kohno and Drinfield \cite{Koh}\cite{Dri}\cite{Chari}:
\begin{theorem}[Kohno-Drinfield]
Let $R$ be the monodromy of the KZ equation for the group $SU(2)$
and let $\bar R$ denotes the $R$-matrix of
the quantum group $U_q (su(2))$ where $su(2)$ denotes the Lie algebra
of $SU(2)$ and $q=e^{\frac{i2\pi}{k+g}}$ where $g=2$. Then there exists a twisting
$F\in U_q(su(2))\otimes U_q(su(2))$ such that
\begin{equation}
\bar R=F^{-1}RF^{-1}
\label{R}
\end{equation}
From this relation we have that the braid group representations obtained from the quantum group $U_q(su(2)$
and obtained from the one-side KZ equation are equivalent.
\end{theorem}
We shall use only the relation (\ref{R}) of this theorem  to derive
the skein relation of the HOMFLY polynomial.

From the property of the quantum group $U_q(su(2))$ we
also have the following formula \cite{Fuc}\cite{Chari}\cite{Kak}:
\begin{equation}
\bar R^{ 2} -(q^{\frac12}-q^{-\frac12}) \bar R
 - I
=0
\label{h6}
\end{equation}
Thus we have
\begin{equation}
\bar R -(q^{\frac12}-q^{-\frac12})I - \bar R^{ -1}
=0
\label{h6a}
\end{equation}
By using this formula we have
\begin{equation}
[\bar R -(q^{\frac12}-q^{-\frac12})I - \bar R^{-1 }]
 FW(z_1,w)W(w,z_2)W(z_3,w)W(w,z_4)F^{-1}
=0
\label{h6b}
\end{equation}
Thus by using the relation (\ref{R}) we have
\begin{equation}
\begin{array}{rl}
& Tr F^{-1}R W(z_1,w)W(w,z_2)W(z_3,w)W(w,z_4)F^{-1}
  \\
- &
(q^{\frac12}-q^{-\frac12})
Tr FW(z_1,z_2)W(z_3,z_4)F^{-1} \\
-
& Tr FW(z_1,w)W(w,z_2)W(z_3,w)W(w,z_4)R^{-1}F
=0
\end{array}
\label{h7a}
\end{equation}
Then by using  the formulas
(\ref{m7a}) and (\ref{m8a})
for upcrossing and undercrossing from (\ref{h6a})
we have
\begin{equation}
\begin{array}{rl}
& Tr
F^{-1}W(z_3,w)W(w,z_2)W(z_1,w)W(w,z_4)F^{-1} \\
- &
(q^{\frac12}-q^{-\frac12})
Tr  W(z_1,z_2)W(z_3,z_4) \\
- &
 Tr\langle FW(z_1,w)W(w,z_4)W(z_3,w)W(w,z_2)F\rangle
=0
\end{array}
\label{h7}
\end{equation}
Let us make a further twist that replace $F^2$ by
$F^2x$  where $x$ denotes a nonzero variable. Then from
(\ref{h7}) we have the following skein relation
for the HOMFLY polynomial:
\begin{equation}
xL_{+}+yL_{0}-x^{-1}L_{-}=0 \label{hh7}
\end{equation}
where we define $y=q^{-\frac12}-q^{\frac12}$ and
that $L_{+}$, $L_{0}$
and $L_{-}$ are defined by
\begin{equation}
\begin{array}{rl}
L_{+}=&
Tr F^{-2}x^{-1}
W(z_3,w)W(w,z_2)W(z_1,w)W(w,z_4)
\\
L_{0}=&
Tr W(z_1,z_2)W(z_3,z_4)
\\
L_{-}=&
Tr xF^{2}W(z_1,w)W(w,z_4)W(z_3,w)W(w,z_2)
\end{array}
\label{h8}
\end{equation}
which are as the HOMFLY polynomials for upcrossing, zero crossing
and undercrossing respectively.

\section{Computation of Quantum Wilson Loop}\label{sec10a}

Let us consider again the quantum Wilson line $W(z_1,z_4)=W(z_1, z_2)W(z_2, z_4)$.
Let us set $z_1=z_4$. In this case the quantum Wilson line forms a closed loop.
Now in (\ref{closed1a}) with $z_1=z_4$ we have that $e^{-t\log  \pm (z_1-z_2)}$
and $e^{t\log \pm (z_1-z_2)}$ which come from the two-side KZ
equations cancel each other and from the multivalued property of
the $\log$ function we have
\begin{equation}
W(z_1, z_1) =R^{n}A_{14} \quad\quad n=0, \pm 1, \pm 2, ...
\label{closed2}
\end{equation}
where $R=e^{-i\pi t}$ is the monodromy of the KZ equation \cite{Chari}. 

{\bf Remark}. It is clear that if we use other representation of the quantum Wilson loop $W(z_1,z_1)$ (such as the representation $W(z_1,z_1)=W(z_1,w_1)W(w_1,w_2)W(w_2,z_1)$) then we will get the same result as (\ref{closed2}). 

{\bf Remark}. For simplicity we shall drop the subscript of $A_{14}$ in (\ref{closed2}) and simply write $A_{14}=A$.

 \section{Defining Quantum Knots and Knot Invariant}\label{sec10}
 
Now
we have that the quantum Wilson loop $W(z_1, z_1)$ corresponds to a closed
curve in the complex plane with starting and ending
point $z_1$.
Let this quantum Wilson loop $W(z_1, z_1)$ represents the unknot. We shall call $W(z_1, z_1)$ as the quantum unknot. Then from
(\ref{closed2}) we have the following invariant
for the unknot:
\begin{equation}
Tr W(z_1, z_1)= Tr R^{n}A \quad\quad n=0, \pm 1, \pm 2, ...
\label{m6}
\end{equation}
where $A=A_{14}$ is a 2-tensor constant matrix operator. 

In the following let us extend the definition (\ref{m6})
to a knot invariant for nontrivial knots.

Let $W(z_i,z_j)$ represent a piece of curve
with starting point $z_i$ and ending point $z_j$.
Then we let
\begin{equation}
W(z_1,z_2)W(z_3,z_4)
\label{m11}
\end{equation}
represent two pieces of uncrossing curve.
Then by interchanging $z_1$ and $z_3$ we have
\begin{equation}
W(z_3,w)W(w,z_2)W(z_1,w)W(w,z_4)
\label{m12}
\end{equation}
represent the curve specified by $W(z_1,z_2)$ upcrossing the
curve specified by $W(z_3,z_4)$.

Now for a given knot diagram we may cut it into a sum of
parts which are formed by two pieces of curves crossing  each other.
Each of these parts is represented
by  (\ref{m12})( For a knot diagram of the unknot
with zero crossings we simply do not need to cut the
knot diagram).
Then we define the trace of a knot with a
given knot diagram by the following form:
\begin{equation}
 Tr \cdot\cdot\cdot
 W(z_3,w)W(w,z_2)W(z_1,w)W(w,z_4)
\cdot\cdot\cdot
 \label{m14}
\end{equation}
where we use (\ref{m12})  to represent the state of the
two pieces of curves specified by
 $W(z_1,z_2)$ and
$W(z_3,z_4)$. The
 $\cdot\cdot\cdot$ means the product
of a sequence of parts represented by
(\ref{m12}) according to the state of
each part. The ordering of the sequence in (\ref{m14})
 follows the ordering of the parts given by the orientation of the
knot diagram. We shall call the sequence of crossings in
the trace (\ref{m14}) as the generalized Wilson
loop of the knot diagram. For the knot diagram of the unknot with zero crossings we simply
let it be $W(z,z)$ and call it the quantum Wilson loop.

We shall
 show that the generalized Wilson loop of a knot diagram has all the properties of the knot diagram  and that
(\ref{m14}) is  a knot invariant. From this we shall call a generalized Wilson loop as a quantum knot. 

\section{Examples of Quantum Knots} 

Before the proof that a generalized Wilson loop of a knot diagram has all the properties of the knot diagram 
in the following let us first consider
some examples to illustrate the way to define (\ref{m14}) and
the way of applying the  braiding formulas (\ref{m7a}),
 (\ref{m8a}) and (\ref{m9}) to
equivalently transform (\ref{m14}) to a simple
expression of the form  $Tr R^{-m}W(z,z)$ where $m$
is an integer. 

Let us first consider the knot in Fig.1.
For this knot we have that (\ref{m14}) is given by
\begin{equation}
Tr W(z_2,w)W(w,z_2)W(z_1,w)W(w,z_1)
\label{m15a}
\end{equation}
where the product of quantum Wilson lines  is from the definition (\ref{m12})
represented a crossing at $w$.
In applying (\ref{m12}) we let $z_1$ be the
starting and the ending point.

\begin{figure}[hbt]
\centering
\includegraphics[scale=0.6]{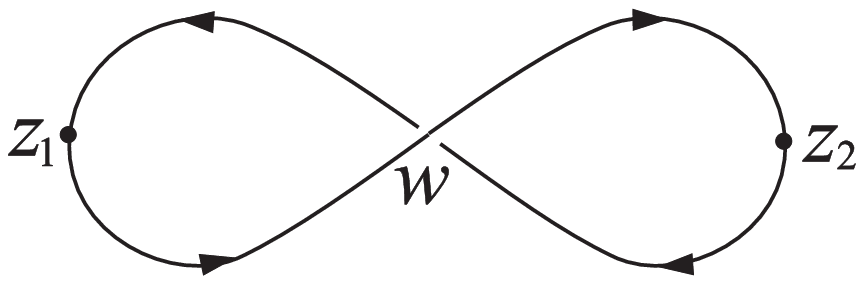}

                Fig.1
\end{figure}
Then  we have that (\ref{m15a}) is equal to
\begin{equation}
\begin{array}{rl}
&Tr W(w,z_2)W(z_1,w)W(w,z_1)W(z_2,w) \\
=&Tr RW(z_1,w)W(w,z_2)R^{-1}
RW(z_2,w)W(w,z_1)R^{-1} \\
=&Tr W(z_1,z_2)W(z_2,z_1) \\
=&Tr W(z_1,z_1)
\end{array}
\label{m16}
\end{equation}
where we have used (\ref{m9}).
We see that (\ref{m16}) is just the knot invariant (\ref{m6}) of
the unknot.
Thus the knot in Fig.1 is with the same knot invariant of the unknot and this
agrees with the fact that this knot is topologically equivalent
to the unknot.

Let us then consider a trefoil knot in Fig.2a. 
By (\ref{m12}) and similar to the above examples
we have that the definition (\ref{m14})
for this knot is given by:
\begin{equation}
\begin{array}{rl}
&Tr W(z_4,w_1)W(w_1,z_2)W(z_1,w_1)W(w_1,z_5)\cdot
W(z_2,w_2)W(w_2,z_6) \\
&W(z_5,w_2)W(w_2,z_3)\cdot
W(z_6,w_3)W(w_3,z_4)W(z_3,w_3)W(w_3,z_1) \\
=&Tr W(z_4,w_1)RW(z_1,w_1)W(w_1,z_2)
R^{-1}W(w_1,z_5)\cdot
W(z_2,w_2)RW(z_5,w_2) \\
&W(w_2,z_6)R^{-1}W(w_2,z_3)\cdot
W(z_6,w_3)RW(z_3,w_3)W(w_3,z_4)R^{-1}W(w_3,z_1) \\
=&Tr
W(z_4,w_1)RW(z_1,z_2)R^{-1}W(w_1,z_5)\cdot
W(z_2,w_2)RW(z_5,z_6)R^{-1}W(w_2,z_3)\cdot \\
&W(z_6,w_3)RW(z_3,z_4)R^{-1}W(w_3,z_1) \\
=&Tr
W(z_4,w_1)RW(z_1,z_2)
W(z_2,w_2)W(w_1,z_5)W(z_5,z_6)R^{-1}W(w_2,z_3)\cdot \\
&W(z_6,w_3)RW(z_3,z_4)R^{-1}W(w_3,z_1) \\
=&Tr
W(z_4,w_1)RW(z_1,w_2)W(w_1,z_6)R^{-1}W(w_2,z_3) \\&
W(z_6,w_3)RW(z_3,z_4)R^{-1}W(w_3,z_1) \\
=&Tr
W(z_4,w_1)W(w_1,z_6)W(z_1,w_2)W(w_2,z_3) \\&
W(z_6,w_3)RW(z_3,z_4)R^{-1}W(w_3,z_1) \\
=&Tr
W(z_4,z_6)W(z_1,z_3)
W(z_6,w_3)RW(z_3,z_4)R^{-1}W(w_3,z_1) \\
=&Tr R^{-1}W(w_3,z_1)
W(z_4,z_6)W(z_1,z_3)
W(z_6,w_3)RW(z_3,z_4) \\
=&Tr
W(z_4,z_6)W(w_3,z_1)R^{-1}W(z_1,z_3)
W(z_6,w_3)RW(z_3,z_4) \\
=&Tr
RW(z_3,z_6)W(w_3,z_1)R^{-1}W(z_1,z_3)
W(z_6,w_3) \\
=&Tr
W(w_3,z_1)W(z_3,z_6)W(z_1,z_3)
W(z_6,w_3) \\
=&Tr
W(z_6,z_1)W(z_3,z_6)W(z_1,z_3)
\end{array}
\label{m21}
\end{equation}
where we have repeatly used (\ref{m9}). Then
 we have that (\ref{m21}) is equal to:
\begin{equation}
\begin{array}{rl}
&TrW(z_6,w_3)W(w_3,z_1)W(z_3,w_3)W(w_3,z_6)W(z_1,z_3)
\\
=&Tr
RW(z_3,w_3)W(w_3,z_1)W(z_6,w_3)W(w_3,z_6)W(z_1,z_3)\\
=&Tr
RW(z_3,w_3)RW(z_6,w_3)W(w_3,z_1)
R^{-1}W(w_3,z_6)W(z_1,z_3)\\
=&Tr
W(z_3,w_3)RW(z_6,z_1)
R^{-1}W(w_3,z_6)W(z_1,z_3)R\\
=&Tr
W(z_3,w_3)RW(z_6,z_3)W(w_3,z_6)\\
=&Tr W(w_3,z_6)W(z_3,w_3)RW(z_6,z_3)\\
=&Tr RW(z_3,w_3)W(w_3,z_6)W(z_6,z_3)\\
=&Tr RW(z_3,z_3)
\end{array}
\label{m22}
\end{equation}
where we have used (\ref{m7a}) and (\ref{m9}).
This is as a knot invariant for the trefoil knot in Fig.2a. 
\begin{figure}[hbt]
\centering
\includegraphics[scale=0.5]{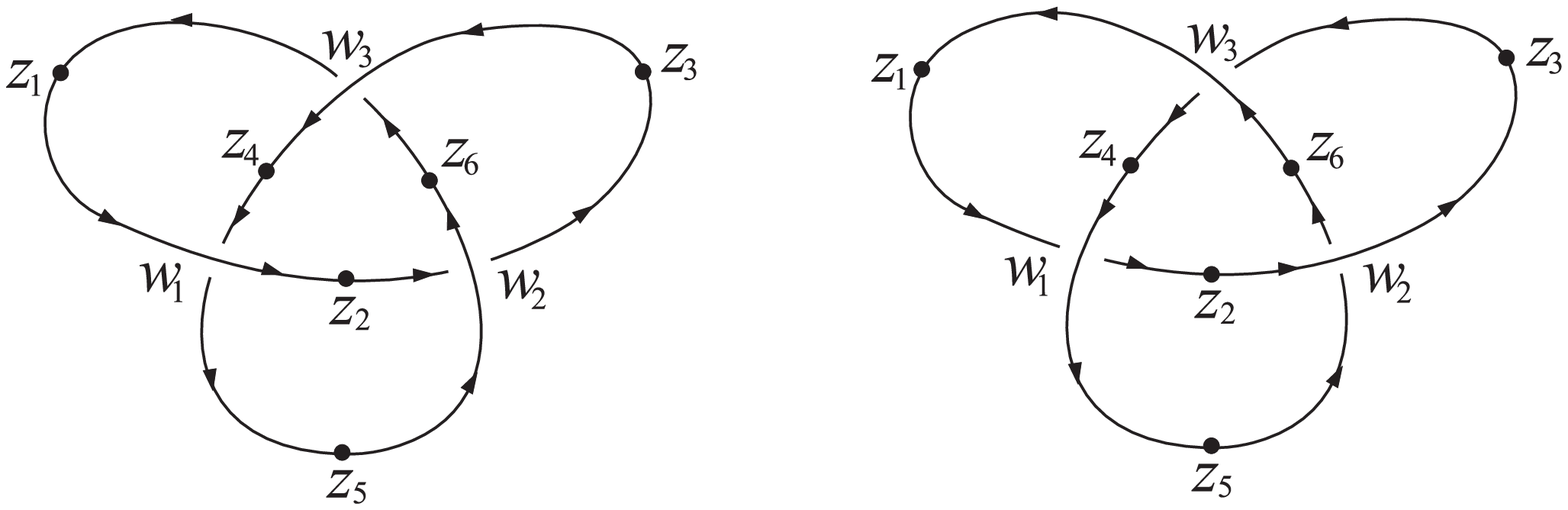}

             Fig.2a  \hspace*{5cm}  Fig.2b

\end{figure}

Then let us consider the trefoil knot in Fig. 2b 
which is the mirror image of the trefoil knot in Fig.2a. 
The definition (\ref{m14}) for this knot is given by:
\begin{equation}
\begin{array}{rl}
&Tr  W(z_1,w_1)W(w_1,z_5)W(z_4,w_1)W(w_1,z_2)\cdot
\\ & \qquad
W(z_5,w_2)W(w_2,z_3)W(z_2,w_2)W(w_2,z_6)\cdot \\
& \qquad W(z_3,w_3)W(w_3,z_1)W(z_6,w_3)W(w_3,z_4) \\
=&Tr W(z_5,z_1)W(z_2,z_5)W(z_1,z_2)
\end{array}
\label{m25}
\end{equation}
where similar to (\ref{m21}) we have repeatly used
(\ref{m9}).
Then we have that (\ref{m25}) is
equal to:
\begin{equation}
\begin{array}{rl}
&Tr W(z_5,z_1)W(z_2,w_1)W(w_1,z_5)W(z_1,w_1)W(w_1,z_2)\\
=&Tr
W(z_5,z_1)W(z_2,w_1)W(w_1,z_2)
W(z_1,w_1)W(w_1,z_5)R^{-1}\\
=&Tr
W(z_5,z_1)W(z_2,w_1)RW(z_1,w_1)W(w_1,z_2)
R^{-1}W(w_1,z_5)R^{-1}\\
=&Tr
R^{-1}W(z_5,z_1)W(z_2,w_1)RW(z_1,z_2)
R^{-1}W(w_1,z_5)\\
=&Tr
W(z_2,w_1)W(z_5,z_2)
R^{-1}W(w_1,z_5)\\
=&Tr
W(z_5,z_2)
R^{-1}W(w_1,z_5)W(z_2,w_1)\\
=&Tr
W(z_5,z_2)
W(z_2,w_1)W(w_1,z_5)R^{-1}\\
=&Tr
W(z_5,z_5)R^{-1}
\end{array}
\label{m26}
\end{equation}
where we have used (\ref{m8a}) and (\ref{m9}).
This is as a knot invariant for the trefoil knot in Fig.2b. 
We notice that
the knot invariants for the two
trefoil knots are different. This shows that these two
trefoil knots are not topologically equivalent.

More examples of the
above quantum knots and knot invariants will be given in a following section.

\section{Generalized Wilson Loops as Quantum Knots}\label{sec11}

Let us now show that the generalized Wilson loop of a knot diagram has all the properties of the knot diagram  and that
(\ref{m14}) is  a knot invariant.  
To this end let us first consider the structure of a knot. Let $K$ be
a knot. Then a knot diagram of $K$ consists of a sequence of
crossings of two pieces of curves cut from the knot $K$ where the
ordering of the crossings can be determined by the orientation of
the knot $K$.  As an example we may consider the two trefoil knots
in the above section. Each trefoil knot is represented by three
crossings of two pieces of curves. These three crossings are
ordered by the orientation of the trefoil knot starting at $z_1$.
Let us denote these three crossings by $1$, $2$ and $3$. Then the
sequence of these three crossings is given by $123$. On the other
hand if the ordering of the three crossings starts from other
$z_i$ on the knot diagram then we have sequences $231$ and $312$.
All these sequences give the same knot diagram and they can be
transformed to each other by circling as follows:
\begin{equation}
123\to 123(1) =231 \to 231(2)=312 \to 312(3)=123 \to \cdot\cdot\cdot
\label{s}
\end{equation}
where $(x)$ means that the number $x$ is to be moved to
the $(x)$ position as indicated.
Let us call (\ref{s}) as the circling property of
the trefoil knot.

As one more example let us consider the figure-eight
knot in Fig.3. 
The simplest knot diagram of this knot has
four crossings.
\begin{figure}[hbt]
\centering
\includegraphics[scale=0.5]{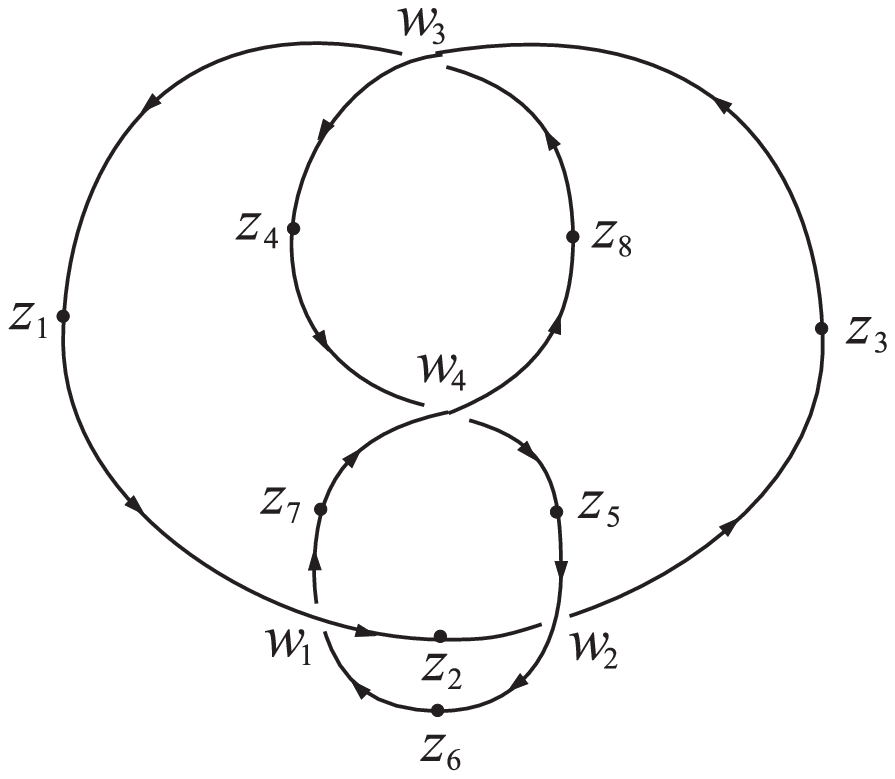}

Fig.3

\end{figure}

Starting at $z_1$ let us denote these crossings by
1, 2, 3 and 4. Then we have the following
circling property of the figure-eight
knot:
\begin{equation}
\begin{array}{rl}
& 1234\to 1234(2) =1342 \to 1342(1)=3421 \to 3421(4)=3214  \to
3214(3)=2143 \\
& \to 2143(1)=2431\to 2431(2)=4312
\to 4312(3)=4123\to 4123(4)=1234 \to \cdot\cdot\cdot
\end{array}
\label{s2}
\end{equation}
We notice that in this cirling of the figure-eight knot there are
subcirclings.

In summary we have that a knot diagram of a knot $K$
can be characterized  as a finite sequence of crossings
of curves which are cut from the knot diagram where the ordering
of the crossings is derived from the orientation of the
knot diagram and has a circling property for which
(\ref{s}) and (\ref{s2}) are examples.

Now let us represent a knot diagram of a knot $K$ by a sequence of
products of Wilson lines representing crossings as in the above
section. Let us call these products of Wilson lines by the term
W-product. Then we  call  this sequence of $W$-products as the
generalized Wilson loop of the knot diagram of a knot $K$.

Let us consider the following two $W$-products:
\begin{equation}
W(z_3,w)W(w,z_2)W(z_1,w)W(w,z_4) \quad \mbox{and}
\quad  W(z_1,z_2)W(z_3,z_4)
\label{w}
\end{equation}
 In the above section we have
shown that these two W-products faithfully represent two oriented
pieces of curves crossing or not crossing  each other where
$W(z_1,z_2)$ and $W(z_3,z_4)$ represent these two pieces of
curves.

Now there is a natural ordering of the  $W$-products
of crossings derived from the orientation of a knot as follows.
Let $W(z_1,z_2)$ and $W(z_3,z_4)$
represent two pieces of curves where the piece of curve
represented by $W(z_1,z_2)$ is before
the piece of curve
represented by $W(z_3,z_4)$  according to the
orientation of a knot. Then the ordering of these
two pieces of curves can be represented by the product
$W(z_1,z_2)W(z_3,z_4)$.
Now let 1 and 2 denote two $W$-products of crossings where we let
1  before 2 according to the orientation of a knot.
Then from the ordering of pieces of curves we have that the
product $12$ represents the ordering of the two
crossings 1 and 2.

Now let a knot diagram of a knot $K$ be given. Let the crossings
of this knot diagram be denoted by $1, 2, \cdot\cdot\cdot, n$ and
let this knot diagram be characterized by the sequence of
crossings $123\cdot\cdot\cdot n$ which is formed according to the
orientation of this knot diagram. On the other hand let us for
simplicity also denote the corresponding $W$-products of crossings
by $1, 2, \cdot\cdot\cdot, n$. Then the whole product of
$W$-products of crossings $123\cdot\cdot\cdot n$ represents the
sequence $123\cdot\cdot\cdot n$ of crossings which is identified
with the the knot diagram. This whole product $123\cdot\cdot\cdot
n$ of $W$-products of crossings is the generalized Wilson loop of
the knot diagram and we denote it by $W(K)$. In the following let us
show that this generalized Wilson loop $W(K)$ has the circling
property of the sequences of crossings of the knot diagram. It
then follows that this generalized Wilson loop represents all the
properties of the sequence $123\cdot\cdot\cdot n$ of crossings of
the knot diagram. Then since this sequence $123\cdot\cdot\cdot n$
of crossings of the knot diagram is identified with the knot
diagram we have that this generalized Wilson loop $W(K)$ can be
identified with the knot diagram and we have the following
theorem.

\begin{theorem}
Each knot $K$ can be faithfully represented by its
generalized Wilson loop $W(K)$
in the sense that if two knot diagrams have the
same generalized Wilson loop then these two knot diagrams must be topologically equivalent.
\end{theorem}
{\bf Proof.} Let us show that the generalized Wilson loop $W(K)$
of a knot diagram of $K$ has the circling property. Let us
consider a product $ W(z_1,z_2)W(z_3,z_4)$ where we first let $z_1,
z_2,z_3$ and $z_4$ be all independent. By solving the two KZ
equations as shown in the above sections we have
\begin{equation}
 W(z_1,z_2)W(z_3,z_4)
= e^{-t \log [\pm (z_3-z_1)]}
  Ae^{t \log [\pm (z_2-z_4)]}
\label{sa}
\end{equation}
where the initial operator $A$ is a 4-tensor as shown in the above sections.
The sign $\pm$ in
(\ref{sa}) reflects that solutions of the KZ
equations are complex multi-valued functions. 
(We remark that the 4-tensor initial operator $A$  in general may not commute with $\Phi_{\pm}(z_1-z_2)=e^{-t \log [\pm (z_1-z_2)]}$ and $\Psi_{\pm}(z_1-z_2)=e^{t \log [\pm (z_1-z_2)]}$).

Then the interchange of $ W(z_1,z_2)$ and $W(z_3,z_4)$
corresponds to that
$z_1$ and $z_3$ interchange their positions and
$z_2$ and $z_4$ interchange their positions respectively.
This interchange
gives a pair of sign changes:
\begin{equation}
(z_3-z_1) \to (z_1-z_3)\qquad
\mbox{and} \qquad
(z_2-z_4) \to (z_4-z_2)
\label{sign}
\end{equation}
From this we have that $W(z_3,z_4)W(z_1,z_2)$
is given by
\begin{equation}
 W(z_3,z_4)W(z_1,z_2)
= e^{-t \log [\pm (z_1-z_3)]}
  Ae^{t \log [\pm (z_4-z_2)]}
\label{sb}
\end{equation}

Now let us set $z_2=z_3$ and $z_1=z_4$ such that the two products
in (\ref{sa}) and (\ref{sb}) form a closed loop. In this case we
have that the initial operator $A$ is reduced from a 4-tensor to a 2-tensor 
and that
$\Phi_{\pm}$ and $\Psi_{\pm}$ act on $A$ by the usual matrix
operation where $A$, $\Phi_{\pm}$ and $\Psi_{\pm}$ are matrices of
the same dimension. 
In this case we have that $A$ commutes with $\Phi_{\pm}$ and $\Psi_{\pm}$ since $\Phi_{\pm}$ and $\Psi_{\pm}$ are Casimir operators on $SU(2)$.

Let us take a definite choice of branch such that
the sign change $z_3-z_1 \to z_1-z_3$ gives
a $i\pi$ difference from the multivalued function
$\log$. Then we have that
$\Phi_{\pm}(z_3-z_1)=R \Phi_{\pm}(z_1-z_3)$. Then since $W(z_1, z_2)$
and $W(z_3, z_4)$ represent two lines with $z_1$, $z_2$ and
$z_3$, $z_4$ as starting and ending points respectively we have that the sign change
$z_2-z_4 \to z_4-z_2$ also gives the same $i\pi$ difference from the
multivalued function $\log$. Thus we have that
$\Psi_{\pm}(z_4-z_2)=R^{-1} \Psi_{\pm}(z_2-z_4)$.
It follows from this pair of sign changes and  that $A$ commutes with
$\Phi_{\pm}$ and $\Psi_{\pm}$ we have that
$W(z_1,z_2)W(z_3,z_4)=W(z_3,z_4)W(z_1,z_2)$ when $z_2=z_3$ and
$z_1=z_4$. This proves the simplest circling property of generalized
Wilson loops.

We remark that in the above proof the pair of sign changes gives
two factors $R$ and $R^{-1}$ which cancel each other and gives the
circling property. We shall later apply the same reason of pair
sign changes to get the general circling property. We also remark
that the proof of this circling property is based on the same
reason as the derivation of the braiding formulas (\ref{m7a}),
(\ref{m8a}) and (\ref{m9}) as shown in the above sections.

Let us consider a product of $n$ quantum Wilson lines
$W(z_i,z_i^{\prime})$, $i=1,...,n$, with the property that the end
points $z_i$, $z_i^{\prime}$  of these quantum Wilson lines are connected
to form a closed loop. From the analysis in the above sections we
have that this product is reduced from a tensor product  
to a 2-tensor. It then follows from (\ref{tensor}) that this product is
of the following form:
\begin{equation}
\prod_{ij}\Phi_{\pm}(z_i-z_j)A
\prod_{ij}
\Psi_{\pm}(z_i^{\prime}-z_j^{\prime})
\label{class}
\end{equation}
where the initial operator $A$ is reduced to a 2-tensor 
and that the $\pm$ signs of
$\Phi_{\pm}(z_i-z_j)$ and $\Psi_{\pm}(z_i-z_j)$ are to be
determined. Then since  $\Phi_{\pm}(z_i-z_j)$ and
$\Psi_{\pm}(z_i-z_j)$ commute with $A$ we can write (\ref{class})
in the form
\begin{equation}
\prod_{ij}\Phi_{\pm}(z_i-z_j)
\prod_{ij}
\Psi_{\pm}(z_i^{\prime}-z_j^{\prime})A
\label{basic}
\end{equation}
where $i\neq j$.
From this formula let us derive the general circling
property as follows.

Let us consider two generalized Wilson lines
denoted by 1 and 2 respectively. Here by the term
generalized Wilson line we mean a product of quantum Wilson lines
with two open ends. As a simple example let us consider
the product $W(z,z_1)W(z_2,z)$. By definition this
is a generalized Wilson line with two open ends
$z_1$ and $z_2$ ($z$ is not an open end).
Suppose that  the two open ends of 1 and 2 are connected. Then
we want to show that $12=21$. This identity is a
generalization of the above interchange of
$W(z_1,z_2)$ and $W(z_3,z_4)$ with $z_2=z_3$ and
$z_1=z_4$.

Because $12$ and $21$ form closed loops we have that
$12$ and $21$ are
products of quantum Wilson lines $W(u_i, u_k)$
(where $u_i$ and $u_k$ denote some $z_p$ or $w_q$ where we use $w_q$
to denote crossing points) such that
for each pair of variables $u_i$ and $u_j$ appearing
at the left side of $W(u_i, u_k)$ and $W(u_j, u_l)$
there is exactly one pair of variables $u_i$ and $u_j$ appearing at
the right side of
$W(u_f, u_i)$ and $W(u_g, u_j)$.
Thus in the formula (\ref{basic}) (with the variables $z$, $z'$ in (\ref{basic})
denoted by variables $u$)
we have that the factors $\Phi_{\pm}(u_i-u_j)$ and $\Psi_{\pm}(u_i-u_j)$
appear in pairs.

As in the above case we have that the interchange
of the open ends of $12$ and $21$  interchanges
 $12$ to $21$. This  interchange gives
changes of the factors $\Phi_{\pm}(u_i-u_j)$
and $\Psi_{\pm}(u_i-u_j)$ as follows.

Let $z_1$ and $z_2$ be the open ends of $1$ and
$z_3$ and $z_4$ be the open ends of $2$
such that $z_1=z_4$ and $z_2=z_3$.
Consider a factor $\Phi_{\pm}(z_1-z_3)$. The interchange of $z_1$ and $z_3$
interchanges this factor to
$\Phi_{\pm}(z_3-z_1)$. Then there is another factor
$\Psi_{\pm}(z_2-z_4)$. The interchange of
 $z_2$ and $z_4$ interchanges this factor to
$\Psi_{\pm}(z_4-z_2)$. Thus this is a pair of sign changes. By the
same reason and the consistent choice of branch as in the above
case we have that the formula (\ref{basic}) is unchanged under
this pair of sign changes.

Then let us consider a factor $\Phi_{\pm}(u_i-u_j)$ of the form
$\Phi_{\pm}(z_1-u_j)$ where $u_i=z_1$ and $u_j$ is not an open
end. Corresponding to this factor we have the  factor
$\Phi_{\pm}(z_3-u_j)$. Then under the interchange of $z_1$ and
$z_3$ we have that $\Phi_{\pm}(z_1-u_j)$ and $\Phi_{\pm}(z_3-u_j)$
change to $\Phi_{\pm}(z_3-u_j)$ and $\Phi_{\pm}(z_1-u_j)$
respectively which gives no change to the formula (\ref{basic}). A
similar result holds for the interchange of $z_2$ and $z_4$ for
factor $\Psi_{\pm}(z_2-u_j)$ and $\Psi_{\pm}(z_4-u_j)$.

It follows that under the interchange of the open ends
of $1$ and $2$ we have the pairs of sign changes
from which  the formula (\ref{basic}) is unchanged.
 This shows that
$ 12  =  21 $.

Then we consider two generalized Wilson products
of crossings which are products of crossings with
four open ends respectively.
Let us again denote them by $1$ and $2$.
Each such generalized Wilson crossing can be regarded
as the crossing of two generalized Wilson lines.
Then the interchange of two open ends of the two
generalized lines of $1$ with the two open ends of the two
generalized lines of $2$ respectively interchanges
$12$ to $21$. Then let us
suppose that the open ends of these two Wilson products
are connected in such a way that the products $12$ and $21$
form  closed loops. In this case we want to show that
$12=21$ which is a circling property of a knot
diagram.
The proof of this equality is again
similar to the above cases. In this case we also have that the interchange of
the open ends of the
two generalized Wilson crossings gives  pairs of sign changes
of the factors $\Phi_{\pm}(u_i-u_j)$ and $\Psi_{\pm}(u_i-u_j)$
in $ 12 $ and $ 12 $.
Then by using (\ref{basic}) we have $ 12  = 21 $.

Let us then consider two generalized Wilson products
of crossings denoted by $1$ and $2$ with open ends  connected in such a way
that two open ends of $1$ are
connected to two open ends of $2$ to form a closed loop.
We want to prove that $12=21$. This will give the subcircling property.

Since a closed loop is formed we have that
each open end of $1$ or of $2$ is connected to a closed loop. In this case
as the above cases we have that
 the products $12$ is with the initial operator $A$ being
 a 2-tensor 
since the open ends of $1$ or
$2$ do not cause $A$ to be a tensor with tensor degree more than $2$
by their connection to the closed loop. Indeed, let $z$ be an open
end of $1$ or $2$. Then it is an end point of a quantum Wilson line
$W(z,z')$ which is a part of $1$ and $2$ such that $z'$ is on the
closed loop formed by $1$ and $2$. Then we have that this quantum Wilson
line $W(z,z')$ is  connected with the closed loop at $z'$. Since
the loop is closed we have that this quantum Wilson line $W(z,z')$ and the
closed loop are connected into a connected line with orientation.
It follows that the open end $z$ gives no additional tensor degree
to the initial operator $A$ for the product $12$ or $21$.
 and that the initial operator $A$  is
still as the initial operator for the closed loop that it is a 2-tensor 
(We remark that in the above section on computation of quantum Wilson loop we see that 
an open quantum Wilson line $W(z_1,z_4)$ and a closed quantum Wilson loop $W(z_1,z_1)$ are with the same 2-tensor initial operator $A$. This shows that the open end $z_1$ of a quantum Wilson line $W(z_1,z_4)$ gives no additional tensor degree
to the initial operator $A$ of the closed quantum Wilson loop $W(z_1,z_1)$. 
This is the same reason that the open end $z$ of the quantum Wilson line
$W(z,z')$ gives no additional tensor degree
to the initial operator $A$ for the product $12$ or $21$).

Now since $A$ is a 2-tensor 
we have that $A$, $\Phi_{\pm}$ and $\Psi_{\pm}$ are as matrices of the same dimension. In this case we have that $A$ commutes
with $\Phi_{\pm}$ and $\Psi_{\pm}$. Then by interchange the open
ends of $1$ with open ends of $2$ we interchange $12$ to $21$.
This interchange again gives  pairs of sign changes. Then since
the initial operator $A$ commutes with $\Phi_{\pm}$ and
$\Psi_{\pm}$ we have that  $12=21$, as was to be proved. Then we
let $12$ and $21$ be connected to another generalized Wilson
product of crossing denoted by $3$ to form a closed loop. Then
from $12=21$ we have $312=321$ and $123=213$. This gives the
subcircling property of generalized Wilson loops. This subcircling
property has been illustrated in the knot diagram of the
fight-eight knot. Then from a case in the above we also have the
circling property $321=213$ between $3$ and $21$.

Continuing in this way we have the circling or
subcircling properties for generalized Wilson loops
whenever the open ends of a product of generalized Wilson lines
or crossings are connected in such a way that
among the open ends a closed loop is formed.
This shows that the generalized Wilson loop of a knot
diagram has the circling property of the knot
diagram.
With this circling property it then follows that the generalized Wilson loop of
a knot diagram completely describes the structure
of the knot diagram.

Now since the generalized Wilson loop of a knot diagram is a
complete copy of this knot diagram we have that two knot diagrams
which can be equivalently moved to each other if and only if the
corresponding generalized Wilson loops can be equivalently moved
to each other. Thus we have that if two knot diagrams have the
same generalized Wilson loop then these two knot diagrams must be
equivalent. This proves the theorem. $\diamond$

{\bf Examples of generalized Wilson loops}

As an example of generalized Wilson loops let us
consider the trefoil knots. Starting at $z_1$ let
the W-product of crossings be denoted by 1, 2 and 3.
Then we have the following circling property of
the generalized Wilson loops of the trefoil knots:
\begin{equation}
123 = 123(1) =231 = 231(2)=312 = 312(3)=123 = \cdot\cdot\cdot
\label{cr1}
\end{equation}
As one more example let us consider the figure-eight
knot.  Starting at $z_1$ let
the W-product of crossings be denoted by 1, 2, 3
and 4.
Then we have the following circling property of
the generalized Wilson loop of the figure-eight
knot:
\begin{equation}
\begin{array}{rl}
& 1234 = 1234(2) =1342 = 1342(1)=3421 = 3421(4)=3214 =
3214(3)=2143 \\
& = 2143(1)=2431 = 2431(2)=4312
= 4312(3)=4123 = 4123(4)=1234 = \cdot\cdot\cdot
\end{array}
\label{cr2}
\end{equation}
$\diamond$

{\bf Definition} We may call a generalized Wilson loop of a knot diagram as a quantum knot since by the above theorem this generalized Wilson loop is a complete copy of the knot diagram. $\diamond$

From the above theorem we have the following theorem.

\begin{theorem}
Let $W(K)$ denote the generalized Wilson loop of a knot $K$. Then
we can write $W(K)$ in the form $R^{-m} W(C)=R^{-m} W(z_1,z_1)$ for some integer $m$
where $C$ denotes a trivial knot and $W(C)=W(z_1,z_1)$ denotes a Wilson loop on $C$ with starting point $z_1$ and ending point $z_1$. From this form we have that the
trace $TrR^{-m}$
is a knot invariant which
classifies knots.  Thus knots can be
classified by the integer $m$.
\end{theorem}

{\bf Proof.} 
Since a generalized Wilson loop $ W(K)$ is in a closed and connected form we have that a generalized Wilson loop $ W(K)$ can be of the form
(\ref{basic}). Thus from the multivalued property of the $\log$
function  and the two-side cancelation in (\ref{basic}) we have
that $ W(K)$ can be of the following (multivalued) form
\begin{equation}
  W(K)= R^{-k}A
\label{BB}
\end{equation}
for some integer $k, k=0, \pm 1, \pm 2, \pm 3, ...$. Furthermore for nontrivial knot $K$ there are some factors $R^{-k_i}$ of $R^{-k}$ coming from the braidings of Wilson lines ( for which the generalized Wilson loop $W(K)$ is formed) by braiding operations such as (\ref{m7a})and (\ref{m8a}). Thus we can write the integer $k$ in the form $k=m+n$ for some integer $m$ and for some integer $n, n=0, \pm 1, \pm 2, ...$ where $n$ is obtained by the two-side cancelations in such a way that the cancelations are obtained when the 
 Wilson lines of the knot diagram for $K$ are connected together to form a  
 Wilson loop $W(C)$ where $C$ is a closed curve which is as an unknot and is of the same form as the knot diagram for $K$ when this knot diagram of $K$ is considered only as a closed curve in the plane (such that the upcrossings and undercrossings are changed to let $K$ be the unknot $C$). From this we have $W(C)=R^{-n}A$ for
$n=0, \pm 1, \pm 2, ...$. Thus $W(K)$ can be written in
the following form  for some $m$:
\begin{equation}
W(K) =  R^{-m} W(C)
\label{i}
\end{equation}

This number $m$ is unique since if there is another number $m_1$ such that 
$W(K)=R^{-m_1} W(C)$ then we have the equality:
\begin{equation}
R^{-m} W(C)=W(K)=R^{-m_1} W(C) 
\label{i1}
\end{equation}
This shows that $R^{-m}=R^{-m_1}$ and thus $m_1=m$.

From (\ref{i}) we also have
\begin{equation}
Tr  W(K) = Tr R^{-m} W(C)
\label{i2}
\end{equation}
for some integer $m$ and that $TrR^{-m} W(C)$ is a knot invariant.

Then let us show that the invariant $TrR^{-m} W(C)$  classifies
knots. Let $K_1$ and $K_2$ be two knots with the same invariant
$TrR^{-m} W(C)$. Then $K_1$ and $K_2$ are both with the same
invariant $R^{-m} W(C)$ where the trace is omitted.
Then by the above formula (\ref{i}) we have
\begin{equation}
 W(K_1) =  R^{-m}W(C) =  W(K_2)
\label{AA1}
\end{equation}
Thus $W(K_1)$ and $W(K_2)$ can be transformed to each other. Thus
$K_1$ and $K_2$ are equivalent. Thus the invariant $TrR^{-m}W(C)$  classifies knots.
It follows that the invariant $TrR^{-m}$ classifies knots and thus knots can be classified by the integer $m$, as was to be proved. $\diamond$

 \section{More Computations of Knot Invariants}

In this section let us give more computations of the knot invariant $Tr R^{-m}$.
We shall show by computation (with the chosen braiding formulae) that the fight-eight knot ${\bf 4_1}$ is assigned with the number $m=3$ and two composite knots composed by two trefoil knots (with the names reef knot and granny knot and denoted by ${\bf 3_1\star 3_1}$ and ${\bf 3_1\times 3_1}$ respectively) are assigned with the numbers $-m=4$ and $-m=9$ respectively. 
The computation is quite tedious. In the next section we shall have a more efficient way to determine the integer $m$. Readers may skip this section
for the first reading.

 Let us first consider the figure-eight knot. From the figure of this knot in a above section we
have that the knot invariant of this knot is given by:
\begin{equation}
\begin{array}{rl}
 & Tr  W(z_6, w_1)W(w_1, z_2)W(z_1, w_1)W(w_1,z_7)\cdot \\
 &W(z_2, w_2)W(w_2, z_6)W(z_5, w_2)W(w_2,z_3)\cdot \\
 &  W(z_8, w_3)W(w_3, z_4)W(z_3, w_3)W(w_3,z_1)\cdot \\
 & W(z_4, w_4)W(w_4, z_8)W(z_7, w_4)W(w_4,z_5)
\end{array}
\label{more}
\end{equation}

In the above computation we have chosen $z_1$ as the
staring point (By the circling property  we may choose any point as the starting point). By repeatedly applying the braiding formulas
(\ref{m7a}),(\ref{m8a}) and (\ref{m9})
 we have that this invariant
is equal to:
\begin{equation}
 Tr R^{-3}W(w_2, z_3)W(z_8, w_2)W(z_3, z_8)
\label{more2}
\end{equation}
Then we have that (\ref{more2}) is equal to
\begin{equation}
 Tr W(z_3, z_8)R^{-3}W(w_2, z)W(z, z_3)W(z_8,z_1)W(z_1, w_2)
\label{more3}
\end{equation}
where $W(w_2, z_3)=W(w_2, z)W(z, z_3)$ with $z$ being a point on the
line represented by  $W(w_2, z_3)$ and
that $W(z_8, w_2)=W(z_8,z_1)W(z_1, w_2)$.
Since $z_1$ is  as the starting and ending point and is an intermediate point we have the
following braiding formula:
\begin{equation}
\begin{array}{rl}
&W(w_2, z_3)W(z_8, w_2) \\
 =&W(w_2,z)W(z, z_3)W(z_8,z_1)W(z_1, w_2)\\
=& R^{-1}W(z_8,z_1)W(z, z_3)W(w_2, z)W(z_1, w_2)\\
=& R^{-1}W(z_8,z_1)W(z_1, w_2)W(w_2, z)W(z, z_3)R^{-1} \\
= & R^{-1}W(z_8, w_2)W(w_2, z_3)R^{-1}
\end{array}
\label{more4}
\end{equation}
 Thus we have that (\ref{more2}) is equal to
\begin{equation}
\begin{array}{rl}
 & Tr W(z_3, z_8)R^{-3}R^{-1}W(z_8, w_2)W(w_2, z_3)R^{-1}\\
=& Tr W(z_3, z_8)R^{-4}W(z_8, z_3)R^{-1} \\
=:&  Tr W(z_3, z_8)R^{-4} \bar W(z_8, z_3)
\end{array}
\label{more5}
\end{equation}
Then  in (\ref{more5}) we have that
\begin{equation}
\begin{array}{rl}
& \bar W(z_8, z_3)\\
 = & W(z_8, z_3)R^{-1} \\
=& W(z_8,z_1)W(z_1,w_1)W(w_1,z_2)W(z_2, z_3)R^{-1}\\
=& W(z_8,z_1)W(z_2, z_3)W(w_1,z_2)W(z_1,w_1)R R^{-1}\\
=& W(z_8,z_1)W(z_2, z_3)W(w_1,z_2)W(z_1,w_1)
\end{array}
\label{more6}
\end{equation}
This shows that $\bar W(z_8, z_3)$ is a generalized Wilson line.
Then since generalized Wilson lines are with the same braiding formulas as Wilson lines
we have that by a braiding formula similar to (\ref{more4}) (for $z_1$ as the starting and ending point and as an intermediate point) the formula (\ref{more5}) is equal to:
\begin{equation}
\begin{array}{rl}
& Tr R^{-4}\bar W (z_8, z_3)W(z_3, z_8)\\
 = & Tr R^{-4} R W(z_3, z_8)\bar W (z_8, z_3)R \\
=& Tr R^{-3}  W(z_3, z_8) W (z_8, z_3)\\
=& Tr R^{-3}  W(z_3, z_3)
\end{array}
\label{more7}
\end{equation}
where the first equality is by a braiding formula which is
similar to the braiding formula (\ref{more4}).
This is the  knot invariant for the figure-eight knot
and we have that $m=3$ for this knot.

Let us then consider the composite knot ${\bf 3_1\star 3_1}$ in Fig.4. 
The trace of the generalized loop of this knot is given by (In Fig.4 one of the two $w_3$ should be $w_1^{'}$):
\begin{equation}
\begin{array}{rl}
 & Tr  W(z_4,w_1)W(w_1,z_2)W(z_1,w_1)W(w_1,z_5)\cdot \\
 & W(z_2,w_2)W(w_2,z_6)W(z_5,w_2)W(w_2,z_3)\cdot \\

& W(z_3,w_1^{'})W(w_1^{'},z_5^{'})W(z_4^{'},w_1^{'})W(w_1^{'},z_2^{'})\cdot
\\
& W(z_5^{'},w_2^{'})W(w_2^{'},z_3^{'})W(z_2^{'},w_2^{'})W(w_2^{'},z_6^{'})
\cdot\\

& W(z_3^{'},w_3^{'})W(w_3^{'},z_1^{'})W(z_6^{'},w_3^{'})W(w_3^{'},z_4^{'})
\cdot \\
 & W(z_6,w_3)W(w_3,z_4)W(z_1^{'},w_3)W(w_3,z_1)
\end{array}
\label{more8}
\end{equation}

\begin{figure}[hbt]
\centering
\includegraphics[scale=0.5]{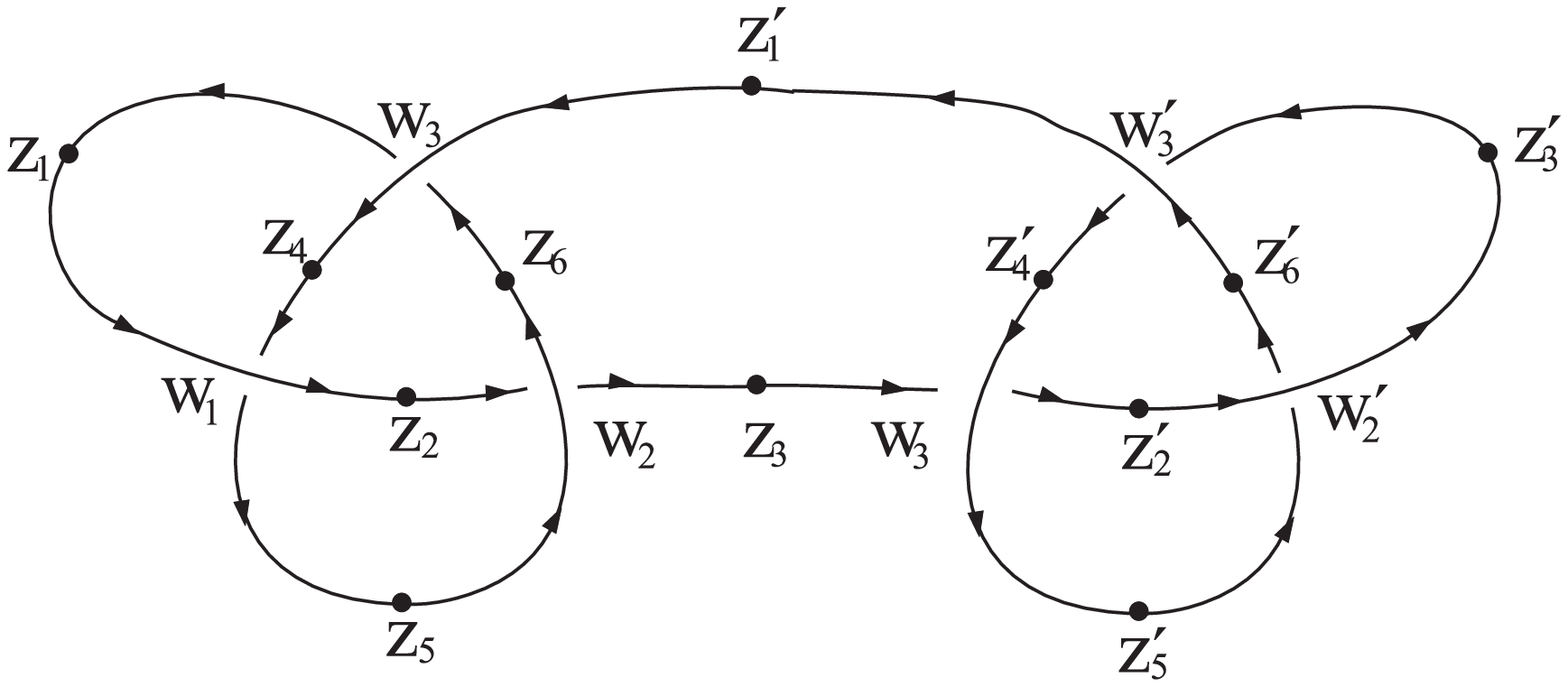}

Fig.4

\end{figure}

By repeatedly applying braiding formula (\ref{m9}) we have that
this invariant is equal to
\begin{equation}
\begin{array}{rl}
  & Tr  W(z_1,w_2^{'})W(z_5,z_1)W(w_2^{'},z_5)\\
=& Tr  W(z_1,w_2)W(w_2,w_2^{'})W(z_5,w_2)W(w_2,z_1)W(w_2^{'},z_5)\\
=& Tr W(z_1,w_2)W(w_2,z_1)W(z_5,w_2)W(w_2,w_2^{'})R^4W(w_2^{'},z_5)
\end{array}
\label{more8a}
\end{equation}
where the braiding of $W(w_2,z_1)$ and $W(w_2,w_2^{'})$ gives $R^4$.
This braiding formula comes from the fact that the Wilson line
$W(w_2,w_2^{'})$ represents a curve with end points $w_2$ and
$w_2^{'}$ such that one and a half loop
is formed which cannot be removed because the end point $w_2^{'}$
is attached to this curve itself to form the closed loop. This
closed loop gives a $3\pi$ phase angle  which is a topological
effect. Thus while the usual
braiding of two pieces of curves gives $R$ which is of a
$\pi$ phase angle  we have that the braiding of
$W(w_2,z_1)$ and $W(w_2,w_2^{'})$ gives $R$ and an additional
$3\pi$ phase angle  and thus gives $R^4$.

Then we have that (\ref{more8}) is equal to
\begin{equation}
\begin{array}{rl}
& Tr W(w_2^{'},z_5)W(z_1,w_2)W(w_2,z_1)W(z_5,w_2)W(w_2,w_2^{'})R^4\\
= & Tr W(w_2^{'},z_5)W(z_1,w_2)RW(z_5,w_2)W(w_2,z_1)
R^{-1}W(w_2,w_2^{'})R^4\\
= & Tr W(w_2^{'},z_5)W(z_1,w_2)RW(z_5,z_1)R^{-1}W(w_2,w_2^{'})R^4 \\
= & Tr RW(z_1,w_2)W(w_2^{'},z_5)W(z_5,z_1)R^{-1}W(w_2,w_2^{'})R^4 \\
= & Tr RW(z_1,w_2)W(w_2^{'},z_1)R^{-1}W(w_2,w_2^{'})R^4 \\
= & Tr W(w_2^{'},z_1)W(z_1,w_2)W(w_2,w_2^{'})R^4 \\
= & Tr W(w_2^{'},w_2)W(w_2,w_2^{'})R^4 \\
= & Tr W(w_2^{'},w_2^{'})R^4
\end{array}
\label{more9}
\end{equation}
This is the invariant of ${\bf 3_1\star 3_1}$.  
Thus we have that $-m=4$ for ${\bf 3_1\star 3_1}$.

Let us then consider the composite knot ${\bf 3_1\times 3_1}$ in Fig.5. 
We have that the trace of the generalized Wilson loop of
${\bf 3_1\times 3_1}$ is given by (In Fig.5 one of the two $w_3$ should be $w_1^{'}$):
\begin{equation}
\begin{array}{rl}
 & Tr  W(z_4,w_1)W(w_1,z_2)W(z_1,w_1)W(w_1,z_5)\cdot \\
 & W(z_2,w_2)W(w_2,z_6)W(z_5,w_2)W(w_2,z_3)\cdot \\
& W(z_4^{'},w_1^{'})W(w_1^{'},z_2^{'})W(z_3,w_1^{'})W(w_1^{'},z_5^{'})
\cdot \\
 & W(z_2^{'},w_2^{'})W(w_2^{'},z_6^{'})W(z_5^{'},w_2^{'})W(w_2^{'},z_3^{'})
\cdot \\
& W(z_6^{'},w_3^{'})W(w_3^{'},z_4^{'})W(z_3^{'},w_3^{'})W(w_3^{'},z_1^{'})
\cdot \\
 & W(z_6,w_3)W(w_3,z_4)W(z_1^{'},w_3)W(w_3,z_1)
\end{array}
\label{more10}
\end{equation}

\begin{figure}[hbt]
\centering
\includegraphics[scale=0.5]{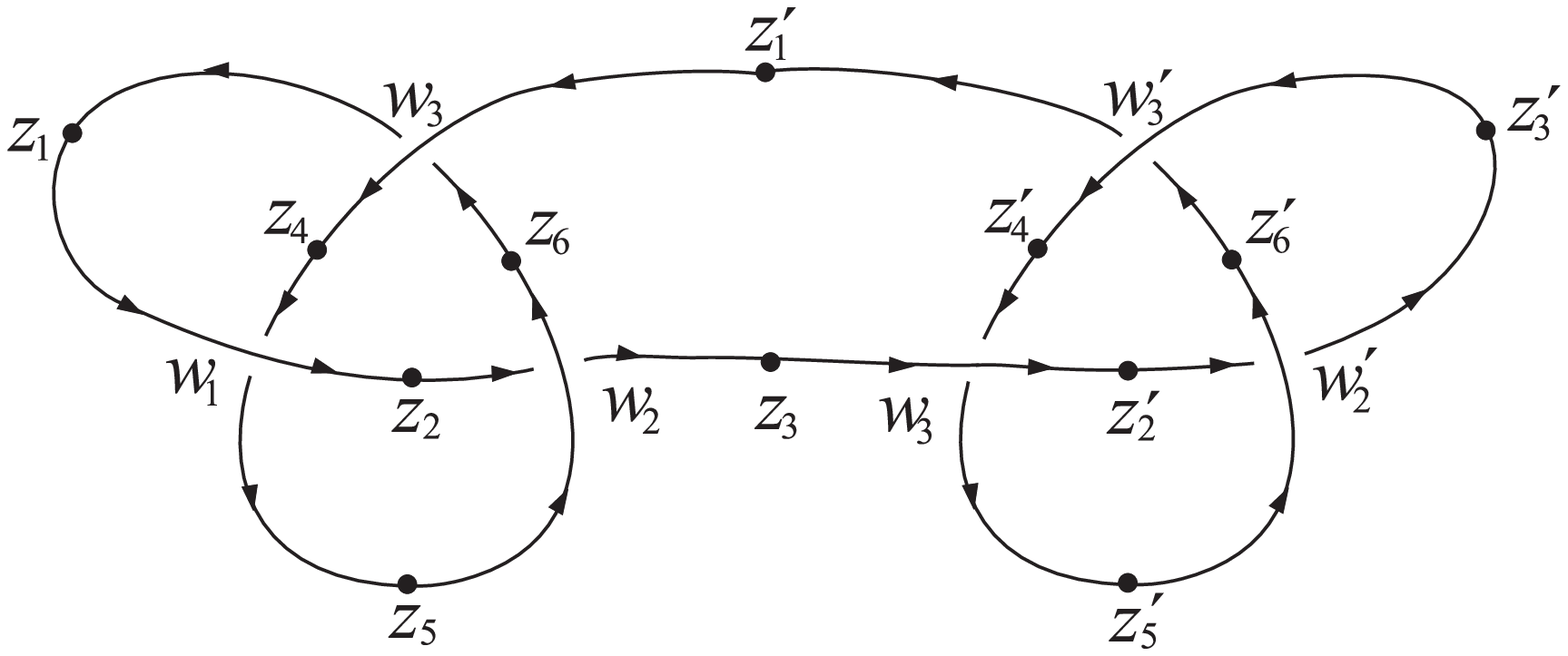}

Fig.5

\end{figure}

By repeatedly applying braiding formulas (\ref{m7a}), (\ref{m8a})
and (\ref{m9}) we have that this invariant is equal to
\begin{equation}
Tr RW(z_1^{'},w_1)R^2W(z_4^{'},z_6^{'})W(w_1,z_4^{'})W(z_6^{'},z_1^{'})
\label{more11}
\end{equation}
where the quantum Wilson line $W(w_1,z_4^{'})$ represents the piece of curve
which starts at $w_1$ and goes through $z_5$, $z_6$, $z_1$ and ends
at $z_4^{'}$. This curve includes a one and a half loop which cannot be
removed since $w_1$ is attached to this curve to form the loop.
This is of the same case as that in the knot ${\bf 3_1\star 3_1}$.
This is a topological property which gives a $3\pi$ phase angle.

We have that (\ref{more11}) is equal to
\begin{equation}
Tr RW(z_1^{'},w_1)R^2
W(z_4^{'},w_1^{'})W(w_1^{'},z_6^{'})
W(w_1,w_1^{'})W(w_1^{'},z_4^{'})W(z_6^{'},z_1^{'})
\label{more12}
\end{equation}
where the piece of curve represented by quantum Wilson line $W(w_1,w_1^{'})$ also
contains the closed loop. Now let this knot ${\bf 3_1\times 3_1}$ be
starting  and ending at $z_6^{'}$. 
Then by the braiding formula on $W(w_1,w_1^{'})$ and $W(z_4^{'},w_1^{'})$
as in the case of the knot ${\bf 3_1\star 3_1}$ we have that (\ref{more12})
is equal to
\begin{equation}
\begin{array}{rl}
& Tr RW(z_1^{'},w_1)R^2 \\
& R^4 W(w_1,w_1^{'})W(w_1^{'},z_6^{'})W(z_4^{'},w_1^{'})
W(w_1^{'},z_4^{'})W(z_6^{'},z_1^{'})\\
= & Tr RW(z_1^{'},w_1)R^6 \\
 & W(w_1,w_1^{'})RW(z_4^{'},w_1^{'})W(w_1^{'},z_6^{'})R^{-1}
 W(w_1^{'},z_4^{'})W(z_6^{'},z_1^{'})\\
= & Tr RW(z_1^{'},w_1)R^6 \\
 & W(w_1,w_1^{'})RW(z_4^{'},z_6^{'})
R^{-1}W(w_1^{'},z_4^{'})W(z_6^{'},z_1^{'})\\
= & Tr RW(z_1^{'},w_1)R^6 \\
 & W(w_1,w_1^{'})RW(z_4^{'},z_6^{'})
W(z_6^{'},z_1^{'})W(w_1^{'},z_4^{'})R^{-1}\\
= & Tr RW(z_1^{'},w_1)R^6 \\
& W(w_1,w_1^{'})RW(z_4^{'},z_1^{'})W(w_1^{'},z_4^{'})R^{-1}\\
= &Tr W(z_1^{'},w_1)R^6
 W(w_1,w_1^{'})RW(z_4^{'},z_1^{'})W(w_1^{'},z_4^{'})
\end{array}
\label{more13}
\end{equation}
where we have repeatedly applied the braiding formula (\ref{m9}).

Now let $z_4^{'}$ be the starting and ending point. Then we have that
(\ref{more13}) is equal to
\begin{equation}
\begin{array}{rl}
& Tr W(z_1^{'},w_1)R^6
 W(w_1,w_1^{'})W(w_1^{'},z_4^{'})W(z_4^{'},z_1^{'})R \\
= & Tr W(z_1^{'},w_1)R^6 W(w_1,z_1^{'})R \\
=& Tr W(z_1^{'},w_1)R^6
 W(w_1,w_1^{'})W(w_1^{'},z_4^{'})W(z_4^{'},w_1^{'})
W(w_1^{'},z_1^{'})R \\
=& Tr W(z_1^{'},w_1)R^6
 W(w_1,w_1^{'})W(w_1^{'},z_1^{'})W(z_4^{'},w_1^{'})W(w_1^{'},z_4^{'}) \\
=: & Tr W(z_1^{'},w_1)R^6 \bar W(w_1,z_1^{'}) \\
= & Tr R^6 \bar W (w_1,z_1^{'}) W(z_1^{'},w_1)
\end{array}
\label{more14}
\end{equation}
where $ \bar W(w_1,z_1^{'})$ denotes the following generalized Wilson line:
\begin{equation}
W(w_1,w_1^{'})W(w_1^{'},z_1^{'})W(z_4^{'},w_1^{'})W(w_1^{'},z_4^{'})
\label{more15}
\end{equation}
Then by the same braiding formula for generalized Wilson lines as that for Wilson lines (with $z_4^{'}$ as the starting and ending point and as an intermediate point)
we have that (\ref{more14}) is equal to:
\begin{equation}
\begin{array}{rl}
& Tr R^6 \bar W (w_1,z_1^{'}) W(z_1^{'},w_1) \\
= & Tr R^6 R W(z_1^{'},w_1) \bar W (w_1,z_1^{'})R \\
= & Tr R^6 R W(z_1^{'},w_1) W (w_1,z_1^{'})RR \\
= & Tr R^9 W(z_1^{'},z_1^{'})
\end{array}
\label{more15a}
\end{equation}
This is the knot invariant for the knot ${\bf 3_1\times 3_1}$. 
Thus we have that $-m=9$ for the knot
${\bf 3_1\times 3_1}$. 
Then we have that the image of ${\bf 3_1\times 3_1}$ is with the knot invariant
$Tr R^{-9} W(z_1^{'},z_1^{'})$. 

\section{A Classification Table of Knots I}\label{sec1a}

In the above sections the computations of the knot invariant $Tr R^{-m}$ is tedious. In this section let us use another method to determine the integer $m$ without carrying out the tedious computations. We shall use only the connected sum operation on knots to find out the integer $m$. For simplicity we use the positive integer $|m|$ to form a classification table of knots where $m$ is assigned to a knot while $-m$ is assigned to 
its mirror image if the knot is not equivalent to its mirror image.
Our main references on the connected sum operation on knots are \cite{Ada}-\cite{Rol}.

Let $\star$ denote the connected sum of two knots such that the
resulting total number of alternating crossings is equal to the
sum of alternating crossings of each of the two knots minus $2$.
As an example we have the reef knot (or the square knot) ${\bf
3_1\star 3_1}$ which is a composite knot composed with the knot
$\bf 3_1$ and its mirror image as in Fig.4.  
This square knot has $6$ crossings
and $4$ alternating crossings. Then let $\times$ denote the
connected sum for two knots such that the resulting total number
of alternating crossings is equal to the sum of alternating
crossings of each of the two knots. As an example we have the granny
knot ${\bf 3_1\times 3_1}$ which is a composite knot composed with
two identical knots $\bf 3_1$ as in Fig.5  
(For simplicity we use one notation $\bf 3_1$ to denote both the trefoil knot and its mirror image though these two knots are nonequivalent). This knot has $6$ alternating crossings which is equal to the total number of crossings. We have that the two operations $\star$ and $\times$ satisfy the commutative law and the associative law
 \cite{Ada}-\cite{Rol}. Further for each knot there is a unique factorization of this knot into a $\star$ and $\times$ operations of prime knots which is similar to the unique factorization of a number into a product of prime numbers \cite{Ada}-\cite{Rol}. We shall show that there is a deeper connection between these two factorizations.

We shall show that we can establish a classification table of knots where each knot is assigned with a number such that 
prime knots are bijectively assigned with prime numbers such that the prime number $2$ corresponds to the trefoil knot
(The trefoil knot will be assigned with the number $1$ and is related to the prime number $2$). We have shown by computation  that the knot $\bf
3_1$ is with $m=1$, the knot $\bf 4_1$ is with $m=3$. Thus there are
no knots assigned with the number $2$ since other knots are with
crossings more than these two knots. We have
shown by computation that the knot ${\bf 3_1\star 3_1}$ is assigned with the
number $4$. Thus we have $1\star 1=4$ (Since knots are assigned with integers we may regard the $\star$ and $\times$ as operations on the set of numbers). This shows that the number
$1$ plays the role of the number $2$. Thus while the knot $\bf
3_1$ is with $m=1$ we may regard this $m=1$ is as the even prime
number $2$. We shall have more to say about this phenomenon
of $1$ and $2$. This phenomenon reflects that the operation $\star$ has partial properties  of addition and multiplication where $m=1$ is assigned to ${\bf 3_1}$ for addition while ${\bf 3_1}$ plays the role of $2$ is for multiplication. The aim of this section is to find out a table of the relation between knots and numbers by using only the operations $\star$ and $\times$ on knots and by using the following data as the initial step for induction: 

{\bf Initial data for induction:} The prime knot ${\bf 3_1}$ is assigned with the number $1$ and it also plays the role of $2$. This means that the number $2$ is not assigned to other knots and is left for the prime knot ${\bf 3_1}$. 
$\diamond$

{\bf Remark}. We shall say that the prime knot ${\bf 3_1}$ is assigned with the number $1$ and is related to the prime number $2$. 
$\diamond$

We shall give an induction 
on the number $n$ of $2^n$ for establishing the table. For each induction step on $n$ because of the special role of the trefoil knot ${\bf 3_1}$ we let the composite knot ${\bf 3_1}^n$ obtained by repeatedly taking $\star$ operation $n-1$ times on the trefoil knot ${\bf 3_1}$ be assigned with the number $2^n$ in this induction. 

Let us first give the following table relating knots and numbers up to $2^5$ as a guide for the induction  
for establishing the whole classification table of knots:

\begin{displaymath}
\begin{array}{|c|c|c|c|} \hline
\mbox{Type of Knot}& \mbox{Assigned number} \,\, |m|
 &\mbox{Type of Knot}& \mbox{Assigned number} \,\, |m|
\\ \hline
{\bf 3_1} & 1 & {\bf 6_3} &  17\\ \hline

 &  2
&  {\bf 3_1\times 4_1} &  18 \\ \hline

{\bf 4_1} &  3 & {\bf 7_1} &  19 \\ \hline

{\bf 3_1\star 3_1} &  4 &
{\bf 4_1\star 5_1} &  20
\\ \hline

{\bf 5_1} & 5 & {\bf 4_1\star(3_1\star 4_1) } &  21
\\ \hline

{\bf 3_1\star 4_1} & 6 & {\bf 4_1\star 5_2} & 22 \\ \hline

{\bf 5_2} &  7 & {\bf 7_2} & 23 \\ \hline

{\bf 3_1\star 3_1\star 3_1} &  8 &
{\bf 3_1\star (3_1\times 3_1)}& 24 \\ \hline

{\bf 3_1\times 3_1} & 9 &
{\bf 3_1\star (3_1\star  5_1)}& 25 \\ \hline

{\bf 3_1\star 5_1} &  10 & {\bf 3_1\star 6_1}
& 26 \\ \hline

{\bf 6_1} &  11 & {\bf 3_1\star (3_1\star  5_2)}& 27
\\ \hline

{\bf 3_1\star 5_2} &  12 & {\bf 3_1\star 6_2} & 28
\\ \hline

{\bf 6_2} &  13 & {\bf 7_3} & 29 \\ \hline

{\bf 4_1\star 4_1} &  14 & { \bf (3_1\star 3_1)\star (3_1\star
4_1)} & 30 \\ \hline

{\bf 4_1\star (3_1\star 3_1)} & 15 & {\bf 7_4} & 31
\\ \hline

{ \bf (3_1\star 3_1)\star (3_1\star 3_1)} & 16 & {\bf (3_1\star
3_1)\star (3_1\star 3_1\star 3_1)} &  32
\\ \hline
\end{array}
\end{displaymath}

From this table we see that the $\star$ operation is similar to the usual multiplication $\cdot$ on numbers. Without the $\times$ operation this $\star$ operation would be exactly the usual multipilcation on numbers if this $\star$ operation is regarded as an operation on numbers.
From this table we see that comparable composite knots (in a sense
from the table and we shall discuss this point later) are grouped
in each of the intervals between two prime numbers. It is
interesting that in each interval composite numbers are one-to-one
assigned to the comparable composite knots while prime numbers are one-to-one
assigned to prime knots. Here a main point is to introduce the $\times$ operation while keeping composite knots correspond to composite numbers and prime knots correspond to prime numbers. To this end we need to have rooms at the positions of composite numbers for the introduction of composite knots obtained by the $\times$ operation. We shall show that these rooms can be obtained by using the special property of the trefoil knot which is assigned with the number $1$ (for the addition property of the $\star$ and $\times$ operations) while this trefoil knot is similar to the number $2$ for the multiplication property of the $\star$ operation.

Let us then carry out the induction steps for obtaining the whole table. To this end let us investigate in more detail the above comparable properties of knots. We have the following definitions and theorems.

{\bf Definition}.
We write $K_1<K_2$ if
$K_1$ is before $K_2$ in the ordering of knots; i.e. the number assigned to $K_1$ is less than the number assigned to $K_2$.

{\bf Definition (Preordering)}. Let two knots be
written in the form $K_1\star K_2$ and $K_1\star K_3$ where we
have determined the ordering of $K_2$ and $K_3$. Then we say that
$K_1\star K_2$ and $K_1\star K_3$ are in a preordering in the sense that we put the ordering of these two knots to follow the ordering of $K_2$ and $K_3$. If this preordering is not changed by conditions from other preorderings on these two knots (which are from other factorization forms of these two knots) then this preordering becomes the ordering of these two knots. 
We shall see that this preordering gives the comparable property in the above table. $ \diamond$

{\bf Remark}.
a) This definition is consistent since if $K_1$ is the unknot then we have 
$K_1\star K_2$=$K_2$ and $K_1\star K_3$=$K_3$ and thus the ordering of $K_1\star K_2$ and $K_1\star K_3$ follows the ordering of $K_2$ and $K_3$.

b) We can also define similarly the preordering of two knots $K_1\times K_2$ and 
$K_1\times K_3$ with the $\times$ operation. $\diamond$

We have the following theorem.
\begin{theorem}

Consider two knots of the form $K_1\star K_2$ and $K_1\star K_3$ where $K_1$, $K_2$ and $K_3$ are prime knots such that $K_2<K_3$. 
Then we have $K_1\star K_2 < K_1\star K_3$.
\label{pre}
\end{theorem}
{\bf Proof}.
Since $K_1$, $K_2$ and $K_3$ are prime knots there are no other factorization forms of the two knots $K_1\star K_2$ and $K_1\star K_3$. Thus these two forms of the two knots are the only way to give preordering to the two knots and thus there are no other conditions to change the preordering given by this factorization form of the two knots. Thus we have that $K_2<K_3$ implies $K_1\star K_2 < K_1\star K_3$.
$\diamond$

\begin{theorem}

Suppose two knots are written in the form $K_1\star K_2$ and $K_1\star K_3$ for determining their ordering and that the other forms of these two knots are not for determining their ordering. Suppose that $K_2<K_3$. Then we have $K_1\star K_2 < K_1\star K_3$.
\end{theorem}
{\bf Proof}.
The proof of this theorem is similar to the proof of the above theorem. Since the other factorization forms are not for the determination of the ordering of the two knots in the factorization form $K_1\star K_2$ and $K_1\star K_3$ we have that the preordering of these two knots in this factorization form becomes the ordering of these two knots. Thus we have $K_1\star K_2 < K_1\star K_3$. $\diamond$

As a generalization of theorem \ref{pre} we have the following theorem.

\begin{theorem}

Let two knots be of the form $K_1\star K_2$ and $K_1\star K_3$ where $K_2$ and $K_3$ are prime knots. Suppose that $K_2<K_3$. Then we have $K_1\star K_2 < K_1\star K_3$.
\label{pre2}
\end{theorem}
{\bf Proof}.
We have the preordering that $K_1\star K_2$ is before $K_1\star K_3$. Then since $K_2$ and $K_3$ are prime knots we have that the other preordering of $K_1\star K_2$ and $K_1\star K_3$ can only from the factorization of $K_1$. Without loss of generality let us suppose that $K_1$ is of the form $K_1=K_4\star K_5$ where 
$K_4< K_5$ and $K_4$ and $K_5$ are prime knots. Then we have the factorization 
$K_1\star K_2=K_4\star (K_5\star K_2)$ and $K_1\star K_3=K_5\star (K_4\star K_3)$. 
This factorization is the only factorization that might change the preordering that $K_1\star K_2$ is before $K_1\star K_3$. Then if $K_2\neq K_4$ or $K_3\neq K_5$
with this  factorization the two knots $K_1\star K_2$ and $K_1\star K_3$ are noncomparable in the sense that this factorization gives no preordering property and that the ordering of these two knots is determined by other conditions. Thus this factorization of $K_1\star K_2$ and $K_1\star K_3$ is not for the determination of the ordering of $K_1\star K_2$ and $K_1\star K_3$. Thus the preordering that $K_1\star K_2$ is before $K_1\star K_3$ is the ordering of $K_1\star K_2$ and $K_1\star K_3$. On the other hand if $K_2=K_4$ and $K_3=K_5$ then this factorization gives the same preordering that $K_1\star K_2$ is before $K_1\star K_3$. Thus for this case the preordering that $K_1\star K_2$ is before $K_1\star K_3$ is also the ordering of $K_1\star K_2$ and $K_1\star K_3$. Thus we have $K_1\star K_2 <K_1\star K_3$. $\diamond$

In addition to the above theorems we have the following theorems.
\begin{theorem}

Consider two knots of the form $K_1\times K_2$ and $K_1\times K_3$ where $K_1$, $K_2$ and $K_3$ are prime knots such that $K_2<K_3$. 
Then we have $K_1\times K_2 < K_1\times K_3$.
\label{1}
\end{theorem}
{\bf Proof}. By using a preordering property for knots with $\times$ operation as similar to that for knots with $\star$ operation we have that the proof of this theorem is similar to the proof of the above theorems. $\diamond$

\begin{theorem}

Let two knots be of the form $K_1\times K_2$ and $K_1\times K_3$ where $K_2$ and $K_3$ are prime knots. Suppose that $K_2<K_3$. Then we have $K_1\times K_2 < K_1\times K_3$.
\end{theorem}
{\bf Proof}. 
The proof of this theorem is also similar to the proof of the theorem \ref{pre2}. $\diamond$

These two theorems will be used for introducing and ordering knots involved with a $\times$ operation which will have the effect of pushing out composite knots with the property of jumping over (to be defined) such that knots are assigned with a prime number if and only if the knot is a prime knot.

Let us investigate more on the property of preordering. We consider the following

{\bf Definition (Preordering sequences)}.
 At the $n$th induction step let the prime knot ${\bf 3_1}$ take a $\star$ operation with the previous $(n-1)$th step. We call this obtained sequence of composite knots as a
preordering sequence. Thus from the ordering of the $(n-1)$th step we have a sequence of composite knots which will be for the construction of the $n$th step. 

Then we let the prime knot ${\bf 4_1}$ (or the knot assigned with a prime number which is $3$ in the $2$nd step as can be seen from the above table) take a $\star$ operation with the previous $(n-2)$th step. From this we get a sequence of composite knots for constructing the $n$th step. 
Then we let the prime knots ${\bf 5_1}$ and ${\bf 5_2}$ (which are prime knots in the same step assigned with a prime number which is $5$ or $7$ in the $3$rd step as can be seen from the above table) take a $\star$ operation with the previous $(n-3)$th step respectively. From this we get two sequences for constructing the $n$th step.

Continuing in this way until the sequences are obtained by a prime knot in the $(n-1)$th step taking a $\star$ operation with the step $n=1$ where the prime knot is assigned with a prime number in the $(n-1)$th step by induction (By induction each prime number greater than $2$ will be assigned to a prime knot). 

We call these obtained sequences of composite knots as 
the preordering sequences of composite knots for constructing the $n$th step. Also we call
the sequences truncated from these preordering sequences as preordering subsequences of composite knots for constructing the $n$th step. $\diamond$

We first have the follwing lemma on preordering sequence.
\begin{lemma}

Let $K$ be a knot in a preordering sequence of the $n$th step. 
Then there exists a room for this $K$ in the $n$th step in the sense that this $K$ corresponds to a number in the $n$th step or in the $(n-1)$th step.
\end{lemma}
{\bf Proof}. 
Let $K$ be of the form $K= {\bf3_1}\star K_1$ where $K_1$ is a knot in the previous $(n-1)$th step. By induction we have that $K_1$ is assigned with a number $a$ which is the position of $K_1$ in the previous $(n-1)$th step. Then since ${\bf3_1}$ corresponds to the number $2$ we have that $K$ corresponds to the number $2\cdot a$ in the $n$th step (We remark that 
$K$ may not be assigned with the number $2\cdot a$). Thus there exists a room for this $K$ in the $n$th step.

Then let $K$ be of the form $K= {\bf4_1}\star K_2$ where $K_2$ is a knot in the previous $(n-2)$th step. By induction we have that $K_2$ is assigned with a number $b$ which is the position of $K_2$ in the previous $(n-2)$th step. Since ${\bf4_1}$ is by induction assigned with the prime number $3$ we have $3\cdot b >3\cdot2^{n-3}>2\cdot2^{n-3}=2^{n-2}$. Also we have 
$3\cdot b <3\cdot2^{n-2}<2^2\cdot2^{n-2}=2^{n}$. Thus there exists a room for this $K$ in the $(n-1)$th step or the $n$th step.

Continuing in this way we have that this lemma holds. $\diamond$

{\bf Remark}. By using  this lemma we shall construct each $n$th step of the classification table by first filling the $n$th step with the preordering subsequences of the $n$th step. $\diamond$

{\bf Remark}.
When the number corresponding to the knot $K$ in the above proof is not in the $n$th step we have that the knot $K$ in the preordering sequences of the $n$th step has the function of pushing a knot $K^{\prime}$ out of the $n$th step where this knot $K^{\prime}$ is related to a number in the $n$th step in order for the knot $K$ to be filled into the $n$th step. 

As an example in the above table the knot $K= {\bf4_1}\star{\bf5_1} $ (related to the number $3\cdot5$) in a preordering sequence of the $5$th step pushes the knot $K^{\prime}= {\bf5_1}\star {\bf5_1}$ related to the number $5\cdot5$ in the $5$th step out of the $5$th step. This relation of pushing out is by the chain $3\cdot5 \to 2\cdot2\cdot5 \to 5\cdot5$.

As another example in the above table the knot $K= {\bf3_1}\star({\bf3_1}\times{\bf3_1})$ (correspoded to the number $2\cdot9$) in a preordering sequence of the $5$th step pushes the knot $K^{\prime}= {\bf3_1}\star({\bf4_1}\star {\bf5_1})$ related to the number $2\cdot3\cdot5$ in the $5$th step out of the $5$th step. This relation of pushing out is by the chain $2\cdot9 \to 2\cdot2\cdot2\cdot3 \to 2\cdot3\cdot5$.
$\diamond$

\begin{lemma}
 
For $n\geq 2$ the preordering subsequences for the $n$th step can cover the whole $n$th step.  
\end{lemma}
{\bf Proof}. 
For $n=2$ we have one preordering sequence with number of knots $=2^0$ which is obtained by the prime knot ${\bf 3_1}$ taking $\star$ operation with the step $n=2-1=1$.  
In addition we have the knot ${\bf 3_1}\star{\bf 3_1}$ which is assigned at the position of $2^n,n=2$ by the induction procedure.  Then since  the total rooms of this step $n=2$ is $2^1$ we have that
these two knots cover this step $n=2$.

For $n=3$ we have one preordering sequence with number of knots $=2^1$ which is obtained by the prime knot ${\bf 3_1}$ taking $\star$ operation with the step $3-1=2$. This sequence cover half of this step $n=3$ which is with $2^{3-1}=2^2$ rooms. 
Then we have one more preordering sequence  which is obtained by the knot ${\bf 4_1}$ taking $\star$ operation with step $n=1$ giving the number $2^0=1$ of knots. This covers half of the remaining rooms of the step $n=3$ which is with $2^{2-1}=2^1$ rooms. Then
in addition we have the knot ${\bf 3_1}\star{\bf 3_1}$ which is assigned at the position of $2^n,n=2$ by the induction procedure. The total of these four knots thus cover the step $n=3$.

For the $n$th step we have one preordering sequence with the number of knots $=2^{n-2}$ which is obtained by the prime knot ${\bf 3_1}$ taking $\star$ operation with the $n-1$th step. This sequence cover half of this $n$th step which is with $2^{n-1}$ rooms. Then we have a preordering sequence  which is obtained by the knot ${\bf 4_1}$ taking $\star$ operation with the $(n-2)$th step giving the number $2^{n-3}$ of knots.  This covers half of the remaining rooms of the $n$th step which is with the remaining $2^{n-2}$ rooms. Then we have one preordering sequence obtained by picking a prime knot (e.g.${\bf 5_1}$) which by induction is assigned with a prime number (e.g. the number 5) taking $\star$ operation with the $(n-3)$th step.
Continue in this way until the knot ${\bf 3_1}^n$ 
is by induction assigned at the position of $2^n$. The total number of these knots is $2^{n-1}$ and thus cover this $n$th step. This proves the lemma. $\diamond$

{\bf Remark}. 
Since there will have more than one prime number in the $k$th steps ($k>2$) 
in the covering 
of the $n$th step there will have knots from  
the preordering sequences in repeat and in overlapping. These knots in repeat and in overlapping  may be deleted when  the ordering of the subsequences of the preordering sequences has been determinated for the covering of the $n$th step. 

Also in the preordering sequences some knots which are in repeat and are not used for the covering of the $n$th step will be omitted when the ordering of the subsequences of the preordering sequences has been determinated for the covering of the $n$th step.
$\diamond$

Let us then introduce another definition for constructing the classification table of knots.

{\bf Definition (Jumping over of the first kind)}. At an induction $n$th step consider a knot
$K^{\prime}$ and the knot $K={\bf 3_1}^n$ which is a $\star$ product of $n$ knots
$3_1$. $K^{\prime}$ is said to jump over $K$, denoted by $K \prec
K^{\prime}$,
 if exist $K_2$ and $K_3$  such that
$K^{\prime}=K_2\star K_3$ and for any $K_0$, $K_1$
such that  $K = K_0\star K_1$ where $K_0$, $K_1$, $K_2$ and $K_3$
are not
equal to ${\bf 3_1}$ we have
\begin{equation}
2^{n_0}<p_1\cdot\cdot\cdot p_{n_2}, \quad
2^{n_1}>q_1\cdot\cdot\cdot q_{n_3}
\label{class1}
\end{equation}
or vice versa
\begin{equation}
2^{n_0}>p_1\cdot\cdot\cdot p_{n_2}, \quad
2^{n_1}<q_1\cdot\cdot\cdot q_{n_3}
\label{class2}
\end{equation}
where $2^{n_0}$, $2^{n_1}$ are the numbers assigned to $K_0$ and $K_1$ respectively ($n_0 +n_1=n$) and
\begin{equation}
K_2= K_{p_1}\star \cdot\cdot\cdot \star K_{p_{n_2}} \quad K_3=
K_{q_1}\star \cdot\cdot\cdot \star K_{q_{n_3}}
\label{class3}
\end{equation}
where $K_{p_i}$, $K_{q_j}$ are prime knots which have been
assigned with prime integers $p_i$, $q_j$ respectively; and the following inequality holds:
\begin{equation}
2^n=2^{n_0+n_1}>p_1\cdot\cdot\cdot p_{n_2}\cdot q_1\cdot\cdot\cdot
q_{n_3}
\label{class33}
\end{equation}
Let us call this definition as the property of jumping over of the first kind. $\diamond$

We remark that the definition of jumping over of the first kind is a generalization of the above
ordering of ${\bf 4_1}\star{\bf 5_1}$ and ${\bf 3_1\star 3_1\star
3_1\star 3_1}$ in the above table in the step $n=4$ of $2^4$. 
Let us consider some examples of this
definition. Consider the knots $K^{\prime}=K_2\star K_3={\bf
4_1\star 5_1}$ and $K={\bf 3_1\star 3_1\star 3_1\star 3_1}$. For
any $K_0$, $K_1$ which are not equal to ${\bf 3_1}$ such that
$K=K_0\star K_1$ we have $2^{n_0}< 5$ and $2^{n_1}> 3$ (or vice
versa) where $3$, $5$ are the numbers of ${\bf 4_1}$ and
${\bf 5_1}$ respectively. Thus we have that ${\bf (3_1\star
3_1)\star(3_1\star 3_1)} \prec {\bf 4_1\star 5_1}$.

As another example we have that
$ {\bf 3_1\star(3_1\star 3_1)\star(3_1\star 3_1)} \prec
{\bf 5_1\star 5_1}$, ${\bf 4_1\star 4_1\star 4_1 }$,  and
${\bf 3_1\star (4_1\star 5_1)}$.

{\bf A Remark on Notation}. At the $n$th step let a composite knot of the form $K_1\star K_2\star\cdot\cdot\cdot\star K_q$ where each $K_i$ is a prime knot such that $K_i$ is assigned with a prime number $p_i$ in the previous $n-1$ steps. Then in general $K_1\star K_2\star\cdot\cdot\cdot\star K_q$ is not assigned with the number $p_1\cdot\cdot\cdot p_q$. However with a little confusion and for notation convenience we shall sometimes use the notation $p_1\cdot\cdot\cdot p_n$ to denote the knot $K_1\star K_2\star\cdot\cdot\cdot\star K_q$ and we say that this knot is related to the number $p_1\cdot\cdot\cdot p_n$ (as similar to the knot ${\bf3_1}$ which is related to the number $2$ but is assigned with the number $1$) and we keep in mind that the knot $K_1\star K_2\star\cdot\cdot\cdot\star K_q$ may not be assigned with the number $p_1\cdot\cdot\cdot p_n$. With this notation then we may say that the composite number $3\cdot5$ jumps over the number $2^4$ which means that the composite knot ${\bf4_1}\star{\bf5_1}$ jumps over the knot ${\bf3_1}\star{\bf3_1}\star{\bf3_1}\star{\bf3_1}$. $\diamond$

{\bf Definition (Jumping over of the general kind)}. At the $n$th step let a composite knot $K^{\prime}$ be related with a number $p_1\cdot p_2\cdot\cdot\cdot p_m$ where the number $p_1\cdot p_2\cdot\cdot\cdot p_m$ is in the $n$th step. Then we say that the knot $K^{\prime}$ (or the number $p_1\cdot p_2\cdot\cdot\cdot p_m$) is of jumping over of the general kind (with respect to the knot $K$ in the definition of the jumping over of the first kind and we also write $K\prec K^{\prime}$) if 
$K$ satisfies one of the following conditions:

1) $K^{\prime}$ (or the number related to $K^{\prime}$) is of jumping over of the first kind; or

2) There exists a $p_i$ (for simplicity let it be $p_1$) and a prime number $q$ such that $p_1$ and $q$ are in the same step $k$ for some $k$ and $q$ is the largest prime number in this step such that the numbers $p_1\cdot p_2\cdot\cdot\cdot p_m$ and $q\cdot p_2\cdot\cdot\cdot p_m$ are also in the same step and that the knot $K_q^{\prime}$ related with $q\cdot p_2\cdot\cdot\cdot p_m$ is of jumping over of the first kind. $\diamond$

{\bf Remark}. The condition 2) is a natural generalization of 1) that if $K^{\prime}$ and the knot $K_q^{\prime}$ are as in 2) then they are both in the preordering sequences of an induction $n$th step or both not. 
Then since $K_q^{\prime}$  is of jumping over into an $(n+1)$th induction step and thus is not in the preordering sequences of the induction $n$th step we have that $K^{\prime}$ is also of jumping over into this $(n+1)$th induction step (even if $K^{\prime}$ is not of jumping over of the first kind). This means that $K^{\prime}$ is of jumping over of the general kind. $\diamond$

{\bf Example of jumping over of the general kind}. At an induction step let $K^{\prime}$ be represented by $11\cdot5\cdot5$ (where we let $p_1=11$) and let $K_q^{\prime}$ be represented by $13\cdot5\cdot5$ (where we let $q=13$). Then $K_q^{\prime}$ is of jumping over of the first kind. Thus we have that $K^{\prime}$ is of jumping over (of the general kind). $\diamond$

We shall show that if $K={\bf 3_1}^n \prec K^{\prime}$ then we can set
$K={\bf 3_1}^n<K^{\prime}$. Thus we have, in the above first example, ${\bf
(3_1\star 3_1)\star(3_1\star 3_1)}< {\bf 4_1\star 5_1}$ while $2^4
> 3\cdot 5$.  From this
property we  shall have rooms for the introduction of the $\times$
operation such that composite numbers are assigned to composite knots and prime numbers are assigned to prime knots. We have the following theorem.

\begin{theorem}

If $K ={\bf 3_1}^n\prec K^{\prime}$ then it is consistent with the preordering property that $K={\bf 3_1}^n <K^{\prime}$ for setting up the table. 
\end{theorem}

For proving this theorem let us first prove the following lemma.
\begin{lemma}

The preordering sequences for the construction of the $n$th step do not have knots of jumping over of the general kind.
\end{lemma}
{\bf Proof of the lemma}.
It is clear that the preordering sequence obtained by the ${\bf 3_1}$ taking a $\star$ operation with the previous $(n-1)$th step has no knots with the jumping over of the first kind property since ${\bf 3_1}$ is corresponded with the number $2$ and the  previous $(n-1)$th step has no knots with the jumping over of the first kind property for this $(n-1)$th step. Then preordering sequence obtained by the ${\bf 4_1}$ taking a $\star$ operation with the previous $(n-2)$th step has no knots with the jump over of the first kind property since ${\bf 4_1}$ is assigned with the number $3$ and $3<2^2$ and the  previous $(n-2)$th step has no knots with the jumping over  of the first kind property for this $(n-2)$th step. Continuing in this way we have that all the knots in these preordering sequences do not satisfy the property of jumping over of the first kind. Then let us show that these preordering sequences have no knots with the property of jumping over of the general kind. Suppose this is not true. Then there exists a knot with the property of jumping over of the general kind and let this knot be represented by a number of the form
$p_1\cdot p_2\cdot\cdot\cdot p_m$ as in the definition of jumping over of the general kind such that there exists a prime number $q$ and that $p_1$ and $q$ are in the same step $k$ for some $k$ and $q$ is the largest prime number in this step such that the numbers $p_1\cdot p_2\cdot\cdot\cdot p_m$ and $q\cdot p_2\cdot\cdot\cdot p_m$ are also in the same step and the knot $K_q$ represented by $q\cdot p_2\cdot\cdot\cdot p_m$ is of jumping over of the first kind. Then since $p_1$ and $q$ are in the same step $k$ we have that the two knots represented by $p_1\cdot p_2\cdot\cdot\cdot p_m$ and $q\cdot p_2\cdot\cdot\cdot p_m$ are elements of two preordering sequences for the construction of the same $n$th step. Now since we have shown that the preordering sequences for the construction of the $n$th step do not have knots of jumping over of the first kind we have that this is a contradiction. This proves the lemma.
$\diamond$

{\bf Proof of the theorem}.
By the above lemma if $K={\bf 3_1}^n\prec K^{\prime}$ then $ K^{\prime}$ is not in the preordering sequences for the $n$th step and thus is pushed out from the $n$th step by the preordering sequences for the $n$th step and thus we have $K={\bf 3_1}^n< K^{\prime}$, as was to be proved.
$\diamond$

{\bf Remark}. We remark that there may exist knots (or numbers related to the knots) which are not in the preordering sequences and are not of jumping over. An example of such special knot is the knot ${\bf4_1}\star{\bf5_1}\star{\bf5_1}$ related with $3\cdot 5\cdot 5$ (but is not assigned with this number). 
$\diamond$

{\bf Definition}.
When there exists a knot which is not in the preordering sequences of the $n$th step and is not of jumping over we put this knot back into the $n$th step to join the  preordering sequences for the filling and covering of the $n$th step. Let us call the preordering sequences together with the knots which are not in the preordering sequences of the $n$th step and are not of jumping over as the generalized preordering sequences (for the filling and covering of the $n$th step).
$\diamond$

{\bf Remark}. By using the generalized preordering sequences for the covering of the $n$th step we have that the knots (or the number related to the knots) in the $n$th step pushed out of the $n$th step by the generalized preordering sequences  are just the knots of jumping over (of the general kind). $\diamond$

Then we also have the following theorem.
\begin{theorem}

At each $n$th step ($n>3$) in the covering of the $n$th step ($n>3$) with the generalized preordering sequences there are rooms for introducing new knots with the $\times$ operations.
\label{times}
\end{theorem}
{\bf Proof}.
We want to show that at each $n$th step ($n>3$) there are rooms for
introducing new knots with the $\times$ operations.
At $n=4$ we have shown that there is the room at the position $9$ for introducing the knot ${\bf 3_1}\times{\bf 3_1}$ with the $\times$ operation.
Let us suppose that this property holds at an
induction step $n-1$. Let us then consider the induction step $n$.
For each $n$ because of the relation between $1$ and $2$ for ${\bf
3_1}$  as a part of the induction step $n$ the number $2^n$ is assigned to the knot
${\bf 3_1}^n$ which is a $\star$
product of $n$ ${\bf 3_1}$. Then we want to show that for this
induction step $n$ by using the $\prec$ property we have rooms for
introducing the $\times$ operation. Let $K^{\prime}$ be a knot
such that ${\bf 3_1}^{n-1}\prec K^{\prime}$ 
and $K^{\prime}=K_2\star K_3$
is as in the definition of $\prec$ of jumping over of the first kind 
such that $p_1\cdot\cdot\cdot p_{n_2}\cdot q_1\cdot\cdot\cdot q_{n_3}<2^{n-1}$ 
(e.g. for $n-1=4$
we have $K^4={\bf 3_1}\star{\bf 3_1}\star{\bf 3_1}\star{\bf 3_1}$
and $K^{\prime}=K_2\star K_3={\bf 4_1}\star{\bf 5_1}$). Then let
us consider $K^{\prime\prime}=({\bf 3_1}\star K_2)\star K_3$.
Clearly we have ${\bf 3_1}^n\prec K^{\prime\prime}$. Thus for each
$K^{\prime}$  we have a $K^{\prime\prime}$ such that
${\bf 3_1}^n\prec K^{\prime\prime}$. Clearly all these $ K^{\prime\prime}$ are different.

Then from $K^{\prime}$ let us
construct more $K^{\prime\prime}$, as follows. Let $K^{\prime}$ be a knot of jumping over
of the first kind. Let
$p_1\cdot\cdot\cdot p_{n_2}$ and $q_1\cdot\cdot\cdot q_{n_3}$ be
as in the definition of jumping over
of the first kind. Then as in the definition of jumping over
of the first kind (w.l.o.g)
we let 
\begin{equation}
2^{n_0}<p_1\cdot\cdot\cdot p_{n_2} \quad \mbox {and}\quad 
2^{n_1}>q_1\cdot\cdot\cdot q_{n_3}
\label{forjumpingover}
\end{equation} 

Then we have
\begin{equation}
2^{n_0+1}<(2\cdot p_1\cdot\cdot\cdot p_{n_2})-1\quad \mbox {and}\quad 
2^{n_1}>q_1\cdot\cdot\cdot q_{n_3}
\label{jumpover12}
\end{equation}
Also it is trivial that we have
$2^{n_0}<(2\cdot p_1\cdot\cdot\cdot p_{n_2})-1$ and
$2^{n_1+1}>q_1\cdot\cdot\cdot q_{n_3}$. This shows that ${\bf 3_1}^n\prec
K^{\prime\prime}:= K_{2a}\star K_{3}$ where $K_{2a}$ denotes the
knot with the number $(2\cdot p_1\cdot\cdot\cdot p_{n_2})-1$ as in
the definition of jumping over of the first kind (We remark that this $K^{\prime\prime}$
corresponds to the knot ${\bf 4_1}\star({\bf 4_1}\star{\bf 4_1})$
in the above induction step where $K_{2a}={\bf 4_1}\star{\bf 4_1}$
is with the number $2\cdot 5-1=3\cdot 3$). 

It is clear that all these more $K^{\prime\prime}$ are different from the above $K^{\prime\prime}$ constructed by the above method of taking a $\star$ operation with ${\bf 3_1}$. 
Thus there are more $ K^{\prime\prime}$ than $K^{\prime}$. Thus at this $n$th step there are rooms for introducing new knots with the $\times$ operations. 
This proves the theorem. $\diamond$

{\bf Remark}.
In the proof of the above theorem we have a way to construct the knots $ K^{\prime\prime}$ by replacing a number $a$ with the number $2a-1$. There is another way of constructing 
the knots $ K^{\prime\prime}$ by replacing a number $b$ with the number $2b+1$. For this way we need to check that the number related to $ K^{\prime\prime}$ is in the $(n-1)$th step for $ K^{\prime\prime}$ of jumping over into the $n$th step.

 As an example let us consider the knot $ K^{\prime}={\bf 4_1}\star{\bf 4_1}\star{\bf 4_1}$ of jumping over into the $6$th step with the following data: 
\begin{equation}
2^{3}<3\cdot3 \quad \mbox {and}\quad 
2^{2}>3
\label{jumpingover9}
\end{equation}

From this data we have:
\begin{equation}
2^{3+1}<2\cdot3\cdot3-1 =17\quad \mbox {and}\quad 
2^{2}>3
\label{jumpingover10a}
\end{equation}
This data gives a knot $ K^{\prime\prime}$ with the related number $3\cdot17 $.

On the other hand from the data (\ref{jumpingover9}) we have:
\begin{equation}
2^{3}<3\cdot3 \quad \mbox {and}\quad 
2^{2+1}>2\cdot3+1
\label{jumpingover10}
\end{equation}
Since $(3\cdot3)(2\cdot3+1)=(2\cdot5-1)(2\cdot3+1)=2\cdot5\cdot2\cdot3+2\cdot2-1<2\cdot2\cdot2^4-1<2^6$ we have that the knot $ K^{\prime\prime}={\bf 4_1}\star{\bf 4_1}\star{\bf 5_2}$ related with the number $3\cdot3\cdot7$ is of jumping over into the $7$th step (We shall show that ${\bf 5_2}$
is assigned with the number $7$). $\diamond$

{\bf Remark}. The above theorem shows that at each $n$th step there are rooms for introducing new knots with the $\times$ operations and thus we may establish a one-to-one correspondence of knots and numbers such that prime knots are bijectively assigned with prime numbers. Further to this theorem we have the following main theorem:

\begin{theorem}

A classification table of knots can be formed (as partly described by the above table up to $2^n$ with $n=5$) by induction on the number $2^n$
such that knots are one-to-one assigned with an integer and prime knots are bijectively assigned with prime numbers such that the prime number $2$ corresponds to the trefoil knot. This assignment is onto the set of positive integers except $2$ where the trefoil knot is assigned with 1 and is related to $2$ and at each $n$th induction step of the number $2^n$ there are rooms for introducing new knots with the $\times$ operations only.

Further this assignment of knots to numbers for the $n$th induction step of the number $2^n$ effectively includes the determination of the distribution of prime numbers in the $n$th induction step and is by induction determined by this assignment for the previous $n-1$ induction steps such that the assignment for the previous $n-1$ induction steps is inherited in this assignment for the $n$th induction step as the preordering sequences in the determination of this assignment for the $n$th induction step.
\label{maintheorem}
\end{theorem}

{\bf Remark}. Let us also call this assignment of knots to numbers as the structure of numbers obtained by assigning numbers to knots. This structure of numbers is the original number system together with the one-to-one assignment of numbers to knots.

{\bf Proof}.
By the above lemmas and theorems we have that
the generalized preordering sequences 
have the function of pushing out those composite knots 
of jumping over from the $n$th step. 
It follows that for step $n>3$ there must exist
 chains of transitions whose initial states are composite knots in repeat 
 (to be replaced by the new composite knots with $\times$ operations only); or the knots of
 jumping over into this $n$th step from the previous $(n-1)$th step; or the knots in the preordering sequences with the $\times$ operations; such that the composite knots 
 of jumping over are pushed out from the $n$th step by these chains. These chains are obtained by ordering the subsequences of preordering sequences 
 such that the preordering property holds in the  
 $n$th step. Further 
the intermediate states of the chains must be positions of composite numbers. This is because that if a chain is transited to an intermediate state which is a position of prime number then there are no composite knots related by  
this prime number and thus this chain can not be transited to the next state and  is stayed at the intermediate state forever and thus the chain can not push out the composite knot of jumping over. Then when a composite knot is  at the position of an intermediate state (which is a position of composite number as has just been proved) then this knot is definitely assigned with this composite number.
Then when a composite knot which is in repeat is  at the position of an intermediate state then this knot is also definitely assigned with this composite number. 
It follows that when the chains are completed we have that the ordering of the subsequences of preordering sequences is determined.

Then the remaining knots (which are not at the transition states of the chains) which are not in repeat are definitely assigned with the number of the positions of these knots in the $n$th step. For these knots the numbers of positions assigned to them are just the number related to them respectively.

Then the remaining knots (which are not at the transition states of the chains) which are in repeat must be replaced by new prime knots because of the repeat and that no other knots related with numbers in this $n$th step in the generalized preordering sequences can be used to replace the remaining knots.
 This means that the numbers of the positions of these remaining knots in repeat are prime numbers in this $n$th step. This is because that if the numbers of the positions assigned to the new prime knot is a composite number then the composite knot related  
with this composite number is either in a transition state or is not in transition.
If the composite knot is not in transition then the composite number related to 
 this composite knot is just the number assigning to this composite knot and since this number is  also assigned to the new prime knot that this is a contradiction. Then if this composite knot is in transition state then this means that the remaining knot is also in transition state and this is a contradiction since by definition the remaining knot is not at the transition states of the chains.

Thus prime numbers in the $n$th step are assigned and  are only assigned to prime knots which replace the remaining knots in repeat 
 in the $n$th step. 
 Thus from the preordering sequences we have determined the positions (i.e. the distribution) of prime numbers in the $n$th step. Now since the preordering sequences are constructed by  the previous steps  we have shown that the basic structure (in the sense of above proof) of this assignment of knots with numbers for the $n$th step (including the determination of the distribution of prime numbers in the $n$th step) is determined by this assignment of knots with numbers for the previous $n-1$ steps. In other words we have that the basic structure of the $n$th induction step 
 is determined by the structure of the previous $n-1$ steps.
 
 To complete the proof of this theorem let us  show that at each $n$th induction step ($n>3$) there are rooms for introducing new composite knots with the $\times$ operations only and we can determine the ordering of these composite knots with the $\times$ operations only in each $n$th induction step. 
 
In the above proof we have shown that the basic structure of the $n$th induction step is determined by the structure of the previous steps such that the positions of the composite knots with the $\times$ operations only in the $n$th induction step are correctedly determined by the structures of the previous steps. These positions are fitted for the corrected composite knots with the $\times$ operations only constructed (by the $\times$ operations) by knots in the previous steps. Thus for this $n$th induction step the introducing and the ordering of composite knots with the $\times$ operations only is also determined by the structures of the previous $n-1$ steps.
 
Further since
the structures of the previous steps are inherited in the structure of the $n$th induction step as the preordering sequences in the determination of the structure of the $n$th induction step we have that all the properties of the structures of the previous steps are inherited in the structure of the $n$th induction step in the determination of the structure of the $n$th induction step. Thus the new composite knots with the $\times$ operations only in the $n$th induction step inherit the ordering properties (such as the preordering property) of composite knots with the $\times$ operations only  
in the previous steps.
(These ordering properties of the composite knots with the $\times$ operations only can be used to find out the corrected composite knots with the $\times$ operations only to be assigned at the corrected positions in the $n$th step).

 With this fact let us then show that at each $n$th induction step ($n>3$) there are rooms for introducing new composite knots with the $\times$ operations only.
 As in the proof of the theorem \ref{times} we first construct more $ K^{\prime\prime}$ by the method following (\ref{forjumpingover}). Let us start at the step $n=4$. For this step we have the knot $K^{\prime}={\bf 4_1}\star{\bf 5_1}$ jumps over into the step $n=5$. For this $K^{\prime}$ we have the following data as in (\ref{forjumpingover}):
\begin{equation}
2^{2}<5 \quad \mbox {and}\quad 
2^{2}>3
\label{jumpingover3}
\end{equation}
From (\ref{jumpingover3}) we construct a $ K^{\prime\prime}$ for the step $n=5$ by the following data:
\begin{equation}
2^{2+1}<2\cdot5-1= 3\cdot3 \quad \mbox {and}\quad 
2^{2}>3
\label{jumpingover4}
\end{equation}
This data gives one more $ K^{\prime\prime}={\bf 4_1}\star{\bf 4_1}\star{\bf 4_1}$. Then from (\ref{jumpingover3}) we construct one more $ K^{\prime\prime}$ for the step $n=5$ by the following data:
\begin{equation}
2^{3}>5 \quad \mbox {and}\quad 
2^{1+1}<2\cdot3-1=5
\label{jumpingover5}
\end{equation}
This data gives one more $K^{\prime\prime}={\bf 5_1}\star{\bf 5_1}$. Thus in this step $n=5$ there are two rooms  for the two knots $ K^{\prime}={\bf 4_1}\star{\bf 5_1}$ and ${\bf 3_1}\star({\bf 3_1}\times{\bf 3_1})$ coming from the preordering sequences  and there exists exactly one room for introducing a  new composite knot with the $\times$ operations only (Recall that we also have a $K^{\prime\prime}={\bf 3_1}\star{\bf 4_1}\star{\bf 5_1}$). 
From the ordering of knots in the previous steps we determine that ${\bf 3_1}\times{\bf 4_1}$ is the composite knot with the $\times$ operations only for this step.

Thus at the $4$th and $5$th steps we can and only can introduce exactly one composite knot with the $\times$ operations only and they are the knots ${\bf 3_1}\times {\bf 3_1}$ and ${\bf 3_1}\times {\bf 4_1}$ respectively.
This shows that at the $4$th and the $5$th steps we can determine the number of prime knots with the minimal number of crossings $=3$ and $=4$ respectively (These two prime knots are denoted by ${\bf 3_1}$ and ${\bf 4_1}$ respectively and we do not distinguish knots with their mirror images for this determination of the ordering of knots with the $\times$ operations only. This also shows that there are rooms for introducing new composite knots with the $\times$ operations only in the $4$th and $5$th steps). 

Then since this property is inherited  
in the $6$th step  we can thus determine that the $6$th step is a step for introducing new composite knots with the $\times$ operations only of the form ${\bf 3_1}\times {\bf 5_{(\cdot)}}$ where ${\bf 5_{(\cdot)}}$ denotes a prime knot with the minimal number of crossings $=5$ (and thus there are rooms for introducing new composite knots with the $\times$ operations only in this $6$th step). Also since the properties in the $4$th and $5$th steps are inherited  
in the $6$th step we can determine
the number of prime knots with the minimal number of crossings $=5$ by the knots of the form ${\bf 3_1}\times {\bf 5_{(\cdot)}}$ as 
this is a property of knots with the $\times$ operations only in the $4$th and $5$th steps
(In the classification table  in the next section we show that there are exactly two composite knots of the form ${\bf 3_1}\times {\bf 5_1}$ and ${\bf 3_1}\times {\bf 5_2}$ in the $6$th step whose ordering are determined by the preordering property of knots and the structure of the $6$th step. This thus shows that there are exactly two prime knots with
the minimal number of crossings $=5$ and they are denoted by ${\bf 5_1}$ and ${\bf 5_2}$ respectively).

Then since the properties of the $4$th, $5$th and $6$th steps are inherited  
in the  $7$th step we can determine that  the $7$th step is a step for introducing new composite knots with the $\times$ operations only of the form ${\bf 3_1}\times {\bf 6_{(\cdot)}}$ where ${\bf 6_{(\cdot)}}$ denotes a prime knot with the minimal number of crossings $=6$ (and thus there are rooms for introducing new composite knots with the $\times$ operations only in this $7$th step). Also since the properties in the $4$th, $5$th and $6$th steps are inherited 
in the $7$th step we can determine
the number of prime knots with the minimal number of crossings $=6$ by the knots of the form ${\bf 3_1}\times {\bf 6_{(\cdot)}}$ as 
this is a property of knots with the $\times$ operations only in the $4$th, $5$th and $6$th steps (In the classification table  in the next section we show that there are exactly three composite knots of the form ${\bf 3_1}\times {\bf 6_1}$, ${\bf 3_1}\times {\bf 6_2}$ and ${\bf 3_1}\times {\bf 6_3}$ in the $7$th step whose ordering are determined by the preordering property of knots and the structure of the $7$th step. This thus shows that there are exactly three prime knots with
the minimal number of crossings $=6$ and they are denoted by ${\bf 6_1}$, ${\bf 6_2}$ and ${\bf 6_3}$ respectively).

Continuing in this way we thus show that at each $n$th induction step $(n>3)$ we can determine the number of prime knots with the minimal number of crossings $=n-1$ and there are rooms for introducing new composite knots with the $\times$ operations only.
This proves the theorem.
$\diamond$

{\bf Example}. Let us consider the above table up to $2^5$ (with $n$ up to $5$) as an example.

For the induction step at $n=2$ (or at $2^2$) we have one preordering sequence obtained by letting ${\bf 3_1}$ to take a $\star$ operation with the step $n=1$ (For the step $n=1$ the number $2^1$ is related to the trefoil knot ${\bf 3_1}$): ${\bf 3_1}\star {\bf 3_1}$. Then we fill the step $n=2$ with this preordering sequence and we have the following ordering of knots for this step $n=2$:
\begin{equation}
{\bf 3_1}\star {\bf 3_1},
{\bf 3_1}\star  {\bf 3_1}
\label{step2}
\end{equation}
where the first ${\bf 3_1}\star  {\bf 3_1}$ placed at the position $3$ is the preordering sequence while the second ${\bf 3_1}\star  {\bf 3_1}$ placed at the position $2^2$ is required by the induction procedure. For this step there is no numbers of jumping over.  Then we have that the first ${\bf 3_1}\star  {\bf 3_1}$ is a repeat of the second ${\bf 3_1}\star  {\bf 3_1}$.  
Thus this repeat one must be replaced by a new prime knot. Let us choose the prime knot ${\bf 4_1}$ to be this new prime knot since ${\bf 4_1}$ is the smallest of prime knots other than the trefoil knot. Then this new prime knot must be at the position of a prime number, as we have proved in the above theorem. Thus we have determined that $3$ is a prime number in this step $n=2$ by using the structure of numbers of step $n=1$ which is only with the prime number $2$.

Then for the induction step at $n=3$ (or at $2^3$) we have two preordering sequence obtained by letting ${\bf 4_1}$ to take a $\star$ operation with the step $n=1$ and by letting ${\bf 3_1}$ to take a $\star$ operation with the step $n=2$: 
\begin{equation}
{\bf 4_1}\star{\bf 3_1}; {\bf 3_1}\star{\bf 4_1},
{\bf 3_1}\star ({\bf 3_1}\star {\bf 3_1})
\label{step3a}
\end{equation}
where the first knot is the preordering sequence obtained by letting ${\bf 4_1}$ to take a $\star$ operation with the step $n=1$ and the second and third knots is the preordering sequence obtained by letting ${\bf 3_1}$ to take a $\star$ operation with the step $n=2$. 

For this step there is no numbers of jumping over and thus there are no chains of transition. Thus the ordering of the above three knots in this step follow the usual ordering of numbers.  
Thus the number assigned to the knot ${\bf 4_1}\star {\bf 3_1}={\bf 3_1}\star{\bf 4_1}$ must be assigned with a number less than that of ${\bf 3_1}\star {\bf 3_1}\star {\bf 3_1}$ by the ordering of ${\bf 3_1}\star{\bf 4_1}$ and ${\bf 3_1}\star {\bf 3_1}\star {\bf 3_1}$ in the second preordering sequence. 
By this ordering of the two preordering sequences  we have that the step $n=3$ is of the following form:
\begin{equation}
{\bf 4_1}\star{\bf 3_1};{\bf 3_1}\star{\bf 4_1},
{\bf 3_1}\star ({\bf 3_1}\star {\bf 3_1});
{\bf 3_1}\star {\bf 3_1}\star {\bf 3_1}
\label{cover}
\end{equation}
where the fourth knot ${\bf 3_1}\star {\bf 3_1}\star {\bf 3_1}$ is put at the position of $2^3$ and is assigned with the number $2^3$ as required by the induction procedure.
 Thus the third knot ${\bf 3_1}\star ({\bf 3_1}\star {\bf 3_1})$  is a repeated one and thus must be replaced by a prime knot and the position of this prime knot is determined to be a prime number.  Thus we have determined that the number $7$ is a prime number.   
 Then since there are no chains of transition we have that the composite knot ${\bf 3_1}\star{\bf 4_1}$ must be assigned with the number related to this knot and this number is $2\cdot 3=6$. Thus the composite knot ${\bf 3_1}\star{\bf 4_1}$ is at the position of $6$ and that 
the first knot ${\bf 4_1}\star{\bf 3_1}$ is a repeat of the second knot and thus must be replaced by a prime knot. Then since this prime knot is at the position of $5$ we have that $5$ is determined to be a prime number. Now the two prime knots at $5$ and $7$ must be the prime knots ${\bf 5_1}$ and ${\bf 5_2}$ respectively since these two knots are  the smallest prime knots other than ${\bf 3_1}$ and ${\bf 4_1}$ (We may just put in two prime knots and then later determine what these two knots will be. If we put in other prime knots then this will not change the distribution of prime numbers determined by the structure of numbers of the previous steps 
and it is only that the prime knots are assigned with incorrect prime numbers. Further as shown in the above proof by using knots of the form ${\bf 3_1}\times{\bf 5_{(\cdot)}}$  we can determine that there are exactly two prime knots with minimal number of crossings $=5$ and they are denoted by ${\bf 5_1}$ and ${\bf 5_2}$ respectively. From this we can then determine that these two prime knots are ${\bf 5_1}$ and ${\bf 5_2}$).
 Thus we have
the following ordering for $n=3$:
\begin{equation} 
{\bf 5_1}<{\bf 3_1}\star{\bf 4_1}<
{\bf 5_2}<{\bf 3_1}\star {\bf 3_1}\star {\bf 3_1}
\label{order1}
\end{equation} 
where ${\bf 5_1}$ is assigned with the prime number $5$ and ${\bf 5_2}$ is assigned with the prime number $7$. This gives the induction step $n=3$. For this step there is no knot with $\times$ operation since there is no knots of jumping over.

Let us then consider the step $n=4$ (or $2^4$). For this step we have the following three preordering sequences obtained from the steps $n=1,2,3$: 
\begin{equation}
\begin{array}{rl}
&{\bf 5_1}\star {\bf3_1};\\
&{\bf 4_1}\star {\bf 4_1},{\bf 4_1}\star {\bf 3_1}\star {\bf 3_1};\\
&{\bf 3_1}\star {\bf5_1},{\bf 3_1}\star {\bf 3_1}\star {\bf 4_1},{\bf 3_1}\star {\bf 5_2},{\bf 3_1}\star {\bf 3_1}\star {\bf 3_1}\star{\bf 3_1};\\
 \end{array}
\label{order2}
\end{equation}
where the third sequence is obtained by taking $\star$ operation of the knot ${\bf 3_1}$ with step $n=3$ while the third sequence is obtained by taking $\star$ operation of the knot ${\bf 4_1}$ with the step $n=2$ and the first sequence is obtained by taking $\star$ operation of the knot ${\bf 5_1}$ with step $n=1$. Then as required by the induction procedure the knot ${\bf 3_1}\star {\bf 3_1}\star {\bf 3_1}\star {\bf 3_1}$ is assigned at the position of $2^4$. The total number of knots in (\ref{order2}) plus this knot is exactly $2^3$ which is the total number of this step $n=4$. 

{\bf Remark}. We  have one more preordering sequence  obtained by taking $\star$ operation of the knot ${\bf 5_2}$ with step $n=1$. This preordering sequence gives the knot
${\bf 5_1}\star {\bf3_1}$. However since the knots in (\ref{order2}) and the knot ${\bf 3_1}\star {\bf 3_1}\star {\bf 3_1}\star {\bf 3_1}$ assigned at the position of $2^4$ are enough for covering this step $n=4$ and that the knot ${\bf 5_1}\star {\bf3_1}$ of this preordering sequence is a repeat of the knot ${\bf 5_1}\star {\bf3_1}$ in (\ref{order2}) that this preordering sequence obtained by taking $\star$ operation of the knot ${\bf 5_2}$ with step $n=1$ can be omitted. $\diamond$

Then to find the chains of transition for this step let us order the three preordering sequences with the following ordering where we rewrite the preordering sequences in column form and the knot ${\bf 3_1}\star {\bf 3_1}\star {\bf 3_1}\star {\bf 3_1}$ assigned at the position of $2^4$ is put to follow the three sequences:
\begin{equation}
\begin{array}{rl}
&{\bf 5_1}\star {\bf3_1};\\
&{\bf 3_1}\star {\bf5_1},\\
&{\bf 3_1}\star {\bf 3_1}\star {\bf 4_1},\\
&{\bf 3_1}\star {\bf 5_2},\\
&{\bf 3_1}\star {\bf 3_1}\star {\bf 3_1}\star{\bf 3_1};\\
&{\bf 4_1}\star {\bf 4_1},\\
&{\bf 4_1}\star {\bf 3_1}\star {\bf 3_1};\\
&{\bf 3_1}\star {\bf 3_1}\star {\bf 3_1}\star{\bf 3_1}
 \end{array}
\label{step4}
\end{equation}
We notice that this column exactly fills the step $n=4$.

For this step we have that the number $3\cdot5$ (or the knot ${\bf 4_1}\star {\bf 5_1}$ related with $3\cdot5$ ) is of jumping over.
From (\ref{step4}) we have the following chain of transition for pushing out ${\bf 4_1}\star {\bf 5_1}$ at $3\cdot5$ by a knot with the $\times$ operation replacing the repeated knot ${\bf 5_1}\star {\bf 3_1}$ at the position of $9=3\cdot3$:
\begin{equation}
{\bf 3_1}\times {\bf 3_1} (\mbox {at} 3\cdot3)\to {\bf 4_1}\star {\bf 4_1}(\mbox {at} 2\cdot7) \to{\bf 3_1}\star {\bf 5_2} (\mbox {at} 2\cdot2\cdot3)\to
{\bf 3_1}\star {\bf 3_1}\star {\bf 4_1} (\mbox {at} 3\cdot5)\to {\bf 4_1}\star {\bf 5_1} (\mbox {pushed out})
\label{Chain2}
\end{equation}
where we choose the knot ${\bf 3_1}\times{\bf 3_1}$ as the knot with the $\times$ operation since ${\bf 3_1}\times{\bf 3_1}$ is the smallest one of such knots.
For this chain the intermediate states are at positions of composite numbers $2\cdot7$, $2\cdot2\cdot3$ and $3\cdot5$. Thus the knots in this chain at the positions of these composite numbers are assigned with these composite numbers respectively.

Then once this chain of pushing out ${\bf 4_1}\star {\bf 5_1}$ at $3\cdot5$ is set up we have that the other knots in repeat must by replaced by prime knots and that their positions must be prime numbers.
These positions are at $11$ and $13$ and thus $11$ and $13$ are determined to be prime numbers (The knot ${\bf 3_1}\star{\bf 3_1}\star{\bf 3_1}\star{\bf 3_1}$ at the end of this step must be assigned with $2^4=16$ by the induction procedure and thus the knot at $13$ is a repeat).
Then  the new prime knots ${\bf 6_1}$ and ${\bf 6_2}$ are suitable knots corresponding to the prime numbers $11$ and $13$ respectively since they are the smallest  prime knots other than ${\bf 3_1}$. ${\bf 4_1}$, ${\bf 5_1}$ and ${\bf 5_2}$ (As the above induction step we may just put in two prime knots and then later determine what these two prime knots will be. As shown in the above proof by using knots of the form ${\bf 3_1}\times{\bf 6_{(\cdot)}}$  we can determine that there are exactly three prime knots with minimal number of crossings $=6$ and they are denoted by ${\bf 6_1}$, ${\bf 6_2}$ and ${\bf 6_3}$ respectively. From this we can then determine that these two prime knots are ${\bf 6_1}$ and ${\bf 6_2}$).

 This completes the step $n=4$. 
Thus the structure of numbers 
of this step (including distribution of prime numbers in this step) is  determined by the structure of numbers of the previous induction steps.
 
Let us then consider the step $n=5$. For this step we have the following four preordering sequences from the previous steps $n=1,2,3,4$:
\begin{equation}
{\bf 6_1}\star{\bf 3_1}
\label{step5a}
\end{equation}
and
\begin{equation}
\begin{array}{rl}
 & {\bf 5_2}\star{\bf 4_1},\\
 & {\bf 5_2}\star ({\bf 3_1}\star {\bf 3_1})
\end{array}
\label{step5b}
\end{equation}
and
\begin{equation}
\begin{array}{rl}
 & {\bf 4_1}\star{\bf 5_1},\\
 & {\bf 4_1}\star ({\bf 3_1}\star{\bf 4_1}),\\
 & {\bf 4_1}\star {\bf 5_2},\\
 & {\bf 4_1}\star ({\bf 3_1}\star {\bf 3_1}\star {\bf 3_1})
 \end{array}
\label{CC2}
\end{equation} 
and
\begin{equation}
\begin{array}{rl}
 & {\bf 3_1}\star ({\bf 3_1}\times {\bf 3_1}),\\
 & {\bf 3_1}\star({\bf 3_1}\star {\bf 5_1}),\\
 & {\bf 3_1}\star {\bf 6_1},\\
 & {\bf 3_1}\star ({\bf 3_1}\star{\bf 5_2}),\\
 & {\bf 3_1}\star{\bf 6_2}, \\
 & {\bf 3_1}\star ({\bf 4_1}\star {\bf 4_1}),\\ 
 & {\bf 3_1}\star ({\bf 3_1}\star {\bf 3_1}\star {\bf 4_1}),\\
 & {\bf 3_1}\star ({\bf 3_1}\star {\bf 3_1}\star {\bf 3_1}\star {\bf 3_1})
\end{array}
\label{CC4}
\end{equation}

The total number of knots (including repeat) in the above sequences plus the knot ${\bf 3_1}\star {\bf 3_1}\star {\bf 3_1}\star {\bf 3_1}\star {\bf 3_1}$ to be assigned at the position of $2^5$ exactly cover this $n=5$ step. 

{\bf Remark}. As similar to the step $n=4$ two preordering sequences ${\bf 5_1}\star {\bf 4_1}, {\bf 5_1}\star {\bf 3_1}\star {\bf 3_1}$ and ${\bf 6_2}\star{\bf 3_1}$ are omitted since these sequences are with knots which are repeats of the knots in the above preordering sequences. $\diamond$

Then to find the chains of transition for this step let us order these four preordering sequences with the following ordering where the knot ${\bf 3_1}\star {\bf 3_1}\star {\bf 3_1}\star {\bf 3_1}\star {\bf 3_1}$ assigned at the position of $2^5$ is put to follow the four sequences:
\begin{equation}
\begin{array}{rl}
& {\bf 6_1}\star {\bf 3_1}; \\
& {\bf 5_2}\star{\bf 4_1},\\
& {\bf 5_2}\star {\bf 3_1}\star {\bf 3_1};\\
& {\bf 4_1}\star{\bf 5_1}, \\
& {\bf 4_1}\star ({\bf 3_1}\star{\bf 4_1}),\\
& {\bf 4_1}\star {\bf 5_2}, \\
& {\bf 4_1}\star ({\bf 3_1}\star {\bf 3_1}\star {\bf 3_1});\\
& {\bf 3_1}\star ({\bf 3_1}\times {\bf 3_1}), \\
& {\bf 3_1}\star({\bf 3_1}\star {\bf 5_1}),\\
& {\bf 3_1}\star {\bf 6_1}, \\
& {\bf 3_1}\star ({\bf 3_1}\star{\bf 5_2}),\\
& {\bf 3_1}\star{\bf 6_2}, \\
& {\bf 3_1}\star ({\bf 4_1}\star {\bf 4_1}), \\
& {\bf 3_1}\star ({\bf 3_1}\star {\bf 3_1}\star {\bf 4_1}),\\
& {\bf 3_1}\star ({\bf 3_1}\star {\bf 3_1}\star {\bf 3_1}\star {\bf 3_1});\\
& ({\bf 3_1}\star {\bf 3_1})\star {\bf 3_1}\star {\bf 3_1}\star {\bf 3_1}
\end{array}
\label{CC5}
\end{equation}

For this step we have three composite knots ${\bf 3_1}\star ({\bf 4_1}\star{\bf 5_1})$, ${\bf 5_1}\star {\bf 5_1}$ and ${\bf 4_1}\star ({\bf 4_1}\star{\bf 4_1})$ (related with $2\cdot3\cdot5$,$5\cdot5$ and $3\cdot3\cdot3$ respectively) of jumping over and there are two new knots ${\bf 4_1}\star {\bf 5_1}$ and ${\bf 3_1}\star ({\bf 3_1}\times{\bf 3_1})$ coming from the previous step. Thus there is a room for the introduction of new knot obtained only by the $\times$ operation. Then this new knot must be the composite knot ${\bf 3_1}\times {\bf 4_1}$ since besides the composite knot ${\bf 3_1}\times {\bf 3_1}$ it is the smallest of composite knots of this kind.

From (\ref{CC5})
there is a chain of transition given by $18\to 21\to 22\to 26\to 28 \to 27$ and the composite knot ${\bf 4_1}\star ({\bf 4_1}\star{\bf 4_1})$ related with $27=3\cdot3\cdot3$ is pushed out into the next step by the composite knot ${\bf 5_2}\star{\bf 4_1}$ at the starting position $18$. Then this repeated knot must be replaced by a new composite knot obtained by the $\times$ operation only and this new composite knot must be the knot ${\bf 3_1}\times {\bf 4_1}$. 

Then the composite knots at the intermediate states are assigned with the numbers of these states respectively.

In addition to the above chain there are two more chains: $24\to 30$ and $20\to 25$.
The
chain $24\to 30$ starts from ${\bf 3_1}\star ({\bf 3_1}\times{\bf 3_1})$ at $24$ and the composite knot ${\bf 3_1}\star ({\bf 4_1}\star{\bf 5_1})$ at $30$ is pushed out by the composite knot ${\bf 3_1}\star ({\bf 3_1}\star{\bf 3_1}\star{\bf 4_1})$. 
Then the chain $20\to 25$ starts from ${\bf 4_1}\star {\bf 5_1}$ at $20$ and  the composite knot 
${\bf 5_1}\star {\bf 5_1}$ at $25$ is pushed out by the composite knot ${\bf 3_1}\star ({\bf 3_1}\star{\bf 5_1})$.
 
Then the knots ${\bf 3_1}\star ({\bf 3_1}\star{\bf 3_1}\star{\bf 4_1})$ and ${\bf 3_1}\star ({\bf 3_1}\star{\bf 5_1})$ at the intemediate states of these two chains are assigned with the numbers $30=2\cdot3\cdot5$ and $25=5\cdot5$ respectively.

Now the remaining repeated composite knots at the positions $17,19, 23,29,31$ must be replaced by new prime knots  and thus $17,19, 23,29,31$ are determined to be prime numbers and they are determined by the prime numbers in the previous induction steps. Then we may follow the usual table of knots to determine that the new prime knots for the prime numbers $17,19, 23,29,31$ are ${\bf 6_3}$, ${\bf 7_1}$, ${\bf 7_2}$, ${\bf 7_3}$ and ${\bf 7_4}$ respectively (As the above induction steps we may just put in five prime knots and then later determine what these five prime knots will be. As shown in the above proof by using knots of the form ${\bf 3_1}\times{\bf 7_{(\cdot)}}$  we can determine the number of prime knots with minimal number of crossings $=7$. From this we can then determine these five prime knots). 

In summary we have the following form of the step $n=5$:
\begin{equation}
\begin{array}{rl}
& {\bf 6_3} \\
& {\bf 3_1}\times{\bf 4_1}\\
& {\bf 7_1}\\
& {\bf 4_1}\star{\bf 5_1} \\
& {\bf 4_1}\star ({\bf 3_1}\star{\bf 4_1})\\
& {\bf 4_1}\star {\bf 5_2} \\
& {\bf 7_2}\\
& {\bf 3_1}\star ({\bf 3_1}\times {\bf 3_1}) \\
& {\bf 3_1}\star({\bf 3_1}\star {\bf 5_1})\\
& {\bf 3_1}\star {\bf 6_1} \\
& {\bf 3_1}\star ({\bf 3_1}\star{\bf 5_2})\\
& {\bf 3_1}\star{\bf 6_2} \\
& {\bf 7_3} \\
& {\bf 3_1}\star ({\bf 3_1}\star {\bf 3_1}\star {\bf 4_1})\\
& {\bf 7_4}\\
& ({\bf 3_1}\star {\bf 3_1})\star {\bf 3_1}\star {\bf 3_1}\star {\bf 3_1}
\end{array}
\label{CC6}
\end{equation}

This completes the induction step at $n=5$. We have that the structure of numbers of this step (including distribution of prime numbers in this step) is determined by the structure of numbers of the previous induction steps.
$\diamond$

\section{A Classification Table of Knots II}\label{sec2a}

Following the above classification table up to $2^5$ let us in this section give the table up to $2^7$.
Again we shall see from the table that the preordering property is clear. At the $7$th step  there is a special composite knot ${\bf 4_1}\star{\bf 5_1}\star{\bf 5_1}$ which is not of jumping over and is not in the preordering sequences (On the other hand the knot ${\bf 5_1}\star{\bf 5_1}\star{\bf 5_1}$ is of jumping over).

We remark again that it is interesting that (by the ordering of composite knots with the $\times$ operation only) at the $6$th step we require exactly two prime knots with minimal number of crossings $=5$ to form the two composite knots obtained by the $\times$ operation only. From this we can determine the number of prime knots with minimal number of crossings $=5$ without using the actual contruction of these prime knots. We then denote these two prime knots by ${\bf 5_1}$ and ${\bf 5_2}$ respectively and the two composite knots obtained by the $\times$ operation only by ${\bf 3_1}\times{\bf 5_1}$ and ${\bf 3_1}\times{\bf 5_2}$ respectively. Similarly at the $7$th step we can determine that there are exactly three prime knots with minimal number of crossings $=6$ and we denote these three prime knots by ${\bf 6_1}$ and ${\bf 6_2}$ and ${\bf 6_3}$ respectively. These three prime knots give the composite knots ${\bf 3_1}\times{\bf 6_1}$, ${\bf 3_1}\times{\bf 6_2}$ and ${\bf 3_1}\times{\bf 6_3}$ respectively. We can then expect that at the next $8$th step we may determine that the number of prime knots with minimal number of crossings $=7$ is $7$ and then at the next $9$th step the number of prime knots with minimal number of crossings $=8$ is $21$, and so on; as we know from the well known table of prime knots \cite{Rol}.
Here the point is that we can determine the number of prime knots with the same minimal number of crossings without using the actual construction of these prime knots (and by using only the classification table of knots).
\begin{displaymath}
\begin{array}{|c|c|c|} \hline
\mbox{Type of Knot}& \mbox{Assigned number} \,\, |m|
 &\mbox{Repeated Knots being replaced}
\\ \hline

{\bf 3_1\star 6_3} & 33 & {\bf } \\ \hline

{\bf 3_1\star(3_1\times 4_1)} & 34 & {\bf } \\ \hline

{\bf 3_1\star 7_1} & 35 & {\bf } \\ \hline

{\bf 3_1\times 5_1} & 36 & {\bf 3_1\star(4_1\star 5_1)} \\ \hline

{\bf 7_5} & 37 & {\bf 3_1\star(4_1\star 3_1\star 4_1)} \\ \hline

{\bf 3_1\times 5_2} & 38 & {\bf 3_1\star(4_1\star 5_2)} \\ \hline

{\bf 3_1\star 7_2} & 39 & {\bf } \\ \hline

{\bf 3_1\star(3_1\star 3_1\times 3_1)} &40  & {\bf } \\ \hline

{\bf 7_6} & 41 & {\bf 3_1\star(3_1\star 3_1\star 3_1\star 5_1)} \\ \hline

{\bf 5_1\star 5_1} & 42 & {\bf } \\ \hline

{\bf 7_7} & 43 & {\bf 5_1\star (3_1\star4_1)} \\ \hline

{\bf 5_1\star 5_2} & 44 & {\bf } \\ \hline

{\bf 4_1\times 4_1} & 45 & {\bf 5_1\star(3_1\star3_1\star3_1), 5_2\star(3_1\star4_1)} \\ \hline

{\bf 5_2\star5_2} & 46 & {\bf } \\ \hline

{\bf 8_1} & 47 & {\bf 5_2\star(3_1\star3_1\star3_1)},{\bf 4_1\star(3_1\times3_1)}\\ \hline

{\bf 4_1\star(3_1\star5_1)} & 48 & {\bf } \\ \hline

{\bf 4_1\star6_1} & 49 & {\bf } \\ \hline

{\bf 4_1\star(3_1\star5_2)} & 50 & {\bf } \\ \hline

{\bf 4_1\star6_2} & 51 & {\bf } \\ \hline

{\bf 4_1\star(4_1\star4_1)} & 52 & {\bf } \\ \hline

{\bf 8_2} & 53 & {\bf 4_1\star(4_1\star3_1\star3_1)} \\ \hline

{\bf 3_1\star(3_1\star3_1\star5_1)} & 54 & {\bf } \\ \hline

{\bf 3_1\star(3_1\star6_1)} & 55 & {\bf } \\ \hline

{\bf 3_1\star(3_1\star3_1\star5_2)} & 56 & {\bf } \\ \hline

{\bf 3_1\star(3_1\star6_2)} & 57 & {\bf } \\ \hline

{\bf 3_1\star7_3} & 58 & {\bf } \\ \hline

{\bf 8_3} & 59 & {\bf 3_1\star(3_1\star3_1\star3_1\star4_1)} \\ \hline

{\bf 3_1\star7_4} &60  & {\bf } \\ \hline

{\bf 8_4} &61  & {\bf 3_1\star(3_1\star3_1\star3_1\star 3_1\star3_1) } \\ \hline

{\bf 4_1\star(4_1\star3_1\star3_1)} & 62 & {\bf } \\ \hline

{\bf 4_1\star(3_1\star3_1\star3_1\star3_1)} & 63 & {\bf } \\ \hline

{\bf 3_1\star(3_1\star3_1\star3_1\star3_1\star3_1)} &64  & {\bf } \\ \hline
\end{array}
\end{displaymath}

\begin{displaymath}
\begin{array}{|c|c|c|} \hline
\mbox{Type of Knot}& \mbox{Assigned number} \,\, |m|
 &\mbox{Repeated Knots being replaced}
\\ \hline
{\bf 3_1\star(3_1\star6_3)} & 65 & {\bf } \\ \hline

{\bf 3_1\times(3_1\times3_1)} & 66 & {\bf 3_1\star(3_1\star3_1\times4_1)} \\ \hline

{\bf 8_5} & 67 & {\bf 3_1\star(3_1\star7_1)} \\ \hline

{\bf 4_1\times5_1} &68  & {\bf 3_1\star(3_1\star4_1\star5_1)} \\ \hline

{\bf 4_1\times(3_1\star4_1)} &69  & {\bf 3_1\star(3_1\star4_1\star3_1\star4_1)} \\ \hline

{\bf 4_1\times5_2} & 70 & {\bf 3_1\star(3_1\star4_1\star5_2)} \\ \hline

{\bf 8_6} & 71 & {\bf 3_1\star(3_1\star7_2)} \\ \hline

{\bf 4_1\star6_3} & 72 & {\bf } \\ \hline

{\bf 8_7} &73  & {\bf 4_1\star(3_1\times4_1)} \\ \hline

{\bf 5_1\star(3_1\times3_1)} &74  & {\bf } \\ \hline

{\bf 5_1\star(3_1\star5_1)} & 75 & {\bf } \\ \hline

{\bf 5_1\star6_1} & 76 & {\bf } \\ \hline

{\bf 5_1\star(3_1\star5_2)} & 77 & {\bf } \\ \hline

{\bf 5_1\star6_2} &78  & {\bf } \\ \hline

{\bf 8_8} & 79 & {\bf 5_1\star(4_1\star4_1), 5_2\star(3_1\star5_1)} \\ \hline

{\bf 5_2\star6_1} &80  & {\bf } \\ \hline

{\bf 5_2\star(3_1\star5_2)} & 81 & {\bf } \\ \hline

{\bf 5_2\star6_2} &82  & {\bf } \\ \hline

{\bf 8_9} &83  & {\bf 5_2\star(4_1\star4_1)} \\ \hline

{\bf 4_1\star7_1} &84  & {\bf } \\ \hline

{\bf 4_1\star(4_1\star5_1)} & 85 & {\bf } \\ \hline

{\bf 4_1\star(4_1\star4_1\star3_1)} &86  & {\bf } \\ \hline

{\bf 4_1\star(4_1\star5_2)} & 87 & {\bf } \\ \hline

{\bf 4_1\star7_2} &88  & {\bf } \\ \hline

{\bf 8_{10}} & 89 & {\bf 4_1\star(3_1\star3_1\times3_1)} \\ \hline

{\bf 4_1\star(3_1\star3_1\star5_1)} &90  & {\bf } \\ \hline

{\bf 4_1\star(3_1\star6_1)} & 91 & {\bf } \\ \hline

{\bf 4_1\star(3_1\star3_1\star5_2)} & 92 & {\bf } \\ \hline

{\bf 4_1\star(3_1\star6_2)} & 93 & {\bf } \\ \hline

{\bf 4_1\star7_3} & 94 & {\bf } \\ \hline

{\bf 4_1\star(4_1\star3_1\star3_1\star3_1)} &95  & {\bf } \\ \hline

{\bf 4_1\star7_4} & 96 & {\bf } \\ \hline
\end{array}
\end{displaymath}

\begin{displaymath}
\begin{array}{|c|c|c|} \hline
\mbox{Type of Knot}& \mbox{Assigned number} \,\, |m|
 &\mbox{Repeated Knots being replaced}
\\ \hline

{\bf 8_{11}} & 97 & {\bf 4_1\star(3_1\star3_1\star3_1\star3_1\star3_1)} \\ \hline

{\bf 4_1\star(5_1\star5_1)} &98  & {\bf 3_1\star(3_1\star6_3)} \\ \hline

{\bf 3_1\star(3_1\star3_1\times4_1)} &99  & {\bf } \\ \hline

{\bf 3_1\star(3_1\star7_1)} &100  & {\bf } \\ \hline

{\bf 8_{12}} & 101 & {\bf3_1\star(3_1\times5_1) } \\ \hline

{\bf 3_1\star7_5} & 102 & {\bf } \\ \hline

{\bf 8_{13}} &103  & {\bf 3_1\star(3_1\times5_2)} \\ \hline

{\bf 3_1\star(3_1\star7_2)} & 104 & {\bf } \\ \hline

{\bf 3_1\star(3_1\star3_1\star3_1\times3_1)} & 105 & {\bf } \\ \hline

{\bf 3_1\star7_6} & 106 & {\bf } \\ \hline

{\bf 8_{14}} &107  & {\bf 3_1\star(5_1\star5_1)} \\ \hline

{\bf 3_1\star7_7} & 108 & {\bf } \\ \hline

{\bf 8_{15}} &109  & {\bf 3_1\star(5_1\star5_2)} \\ \hline

{\bf 3_1\times6_1} &110  & {\bf 3_1\star(5_1\star5_2)} \\ \hline

{\bf 3_1\times(3_1\star5_2)} & 111 & {\bf 3_1\star(4_1\times4_1)} \\ \hline

{\bf 3_1\times6_2} &112 & {\bf 3_1\star(5_2\star5_2)} \\ \hline

{\bf 8_{16}} & 113 & {\bf 3_1\star(5_2\star5_2)} \\ \hline

{\bf 3_1\star8_1} & 114 & {\bf } \\ \hline

{\bf 3_1\times6_3} & 115 & {\bf 3_1\star(4_1\star3_1\star5_1),3_1\star(4_1\star4_1\star4_1)} \\ \hline

{\bf 3_1\star8_2} &116  & {\bf } \\ \hline

{\bf 3_1\star(3_1\star3_1\star3_1\star5_1)} &117  & {\bf } \\ \hline

{\bf 3_1\star(3_1\star3_1\star6_1)} & 118 & {\bf } \\ \hline

{\bf 3_1\star(3_1\star3_1\star3_1\star5_2)} & 119 & {\bf } \\ \hline

{\bf 3_1\star(3_1\star3_1\star6_2)} &120  & {\bf } \\ \hline

{\bf 3_1\star(3_1\star3_1\star7_3)} &121  & {\bf } \\ \hline

{\bf 3_1\star8_3} & 122 & {\bf } \\ \hline

{\bf 3_1\star(3_1\star7_4)} &123  & {\bf } \\ \hline

{\bf 3_1\star8_4} &124  & {\bf } \\ \hline

{\bf 3_1\times(3_1\times4_1)} &125  & {\bf 3_1\star(4_1\star4_1\star3_1\star3_1)} \\ \hline

{\bf 3_1\star(4_1\star3_1\star3_1\star3_1\star3_1)} & 126 & {\bf } \\ \hline

{\bf 8_{17}} &127  & {\bf 3_1\star(3_1\star3_1\star3_1\star3_1\star3_1\star3_1)} \\ \hline

{\bf 3_1\star(3_1\star3_1\star3_1\star3_1\star3_1\star3_1)} & 128 & {\bf } \\ \hline

\end{array}
\end{displaymath}

\section{Examples of Quantum Links and Link Invariant}\label{sec12}

Let us extend the above quantum knots and knot invariant to quantum links and link
invariant. Let us first consider some examples to see how the quantum link and link invariant is defined.
Let us consider the link in Fig.6a.  
We may let the two
knots of this link be with $z_1$ and $z_4$ as the initial and final end point respectively.
We let the ordering
of these two knots be such that when the $z$ parameter goes
 one loop on one knot then the $z$ parameter for another knot also goes one loop.
 The trace invariant (\ref{m14}) for this link is given by:
\begin{equation}
\begin{array}{rl}
&Tr W(z_3,w_1)W(w_1,z_2)W(z_1,w_1)W(w_1,z_4)\cdot\\
& \qquad W(z_4,w_2)W(w_2,z_1)W(z_2,w_2)W(w_2,z_3)
\end{array}
\label{l1a}
\end{equation}
We let the ordering of the Wilson lines in (\ref{l1a})
be such that $W(z_1,z_2)$ and $W(z_4,z_3)$
start first. Then next $W(z_2,z_1)$ and $W(z_3,z_4)$
follows. Form this ordering we have that (\ref{l1a}) is
equal to:
\begin{equation}
\begin{array}{rl}
&Tr RW(z_1,w_1)W(w_1,z_2)W(z_3,w_1)W(w_1,z_4)\cdot\\
& \qquad W(z_4,w_2)W(w_2,z_3)W(z_2,w_2)W(w_2,z_1)R^{-1}
 \\
=&Tr W(z_1,z_2)W(z_3,z_4)
W(z_4,z_3)W(z_2,z_1)
 \\
=&Tr W(z_2,z_2)W(z_3,z_3)
\end{array}
\label{l1}
\end{equation}
where we have used (\ref{m7a}) and (\ref{m8a}).
Since by definition (\ref{m14}) we have that $Tr W(z_2,z_2)W(z_3,z_3)$ is the knot
invariant for two unlinking trivial knots, equation (\ref{l1}) shows that the link in 
Fig.6a 
is topologically equivalent to two unlinking trivial knots.
Similarly we can show that the link in Fig.6b 
is topologically equivalent to two unlinking trivial knots.

\begin{figure}[hbt]
\centering
\includegraphics[scale=0.4]{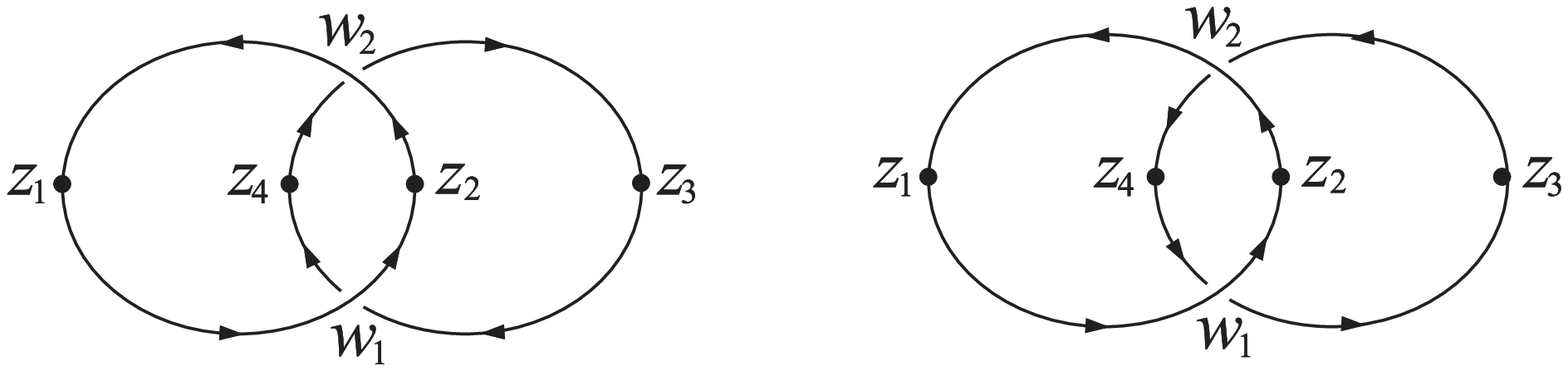}

Fig.6a   
\hspace*{3.5cm} Fig.6b 
\end{figure}

\begin{figure}[hbt]
\centering
\includegraphics[scale=0.4]{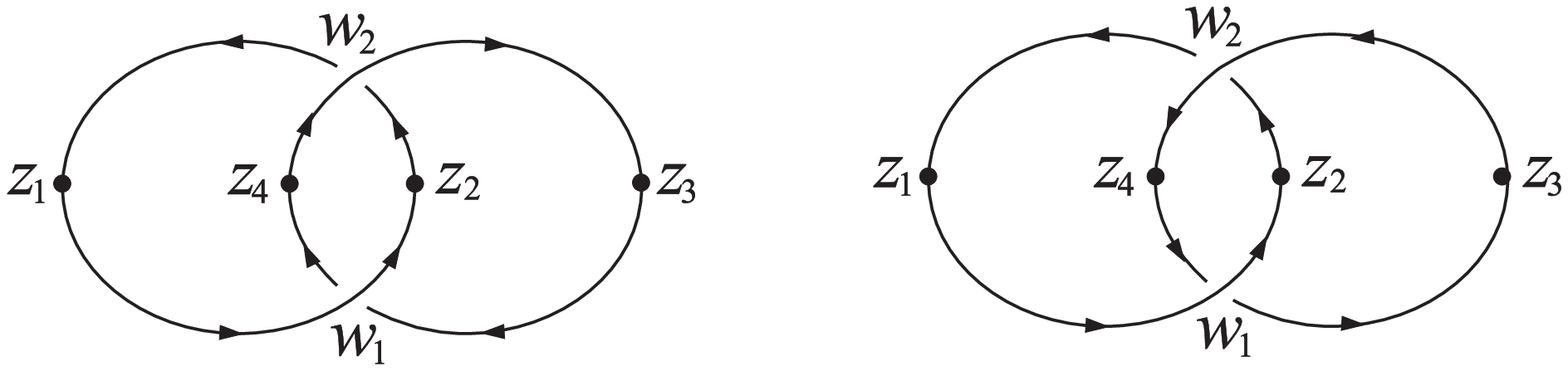}

Fig.7a   
\hspace*{3.5cm} Fig.7b 
\end{figure}

Let us then consider the Hopf link in 
Fig.7a. 
The
trace invariant (\ref{m14}) for this link is given by:
\begin{equation}
\begin{array}{rl}
& Tr W(z_3,w_1)W(w_1,z_2)W(z_1,w_1)W(w_1,z_4)\cdot\\
& \qquad W(z_2,w_2)W(w_2,z_3)W(z_4,w_2)W(w_2,z_1)
\end{array}
\label{l3}
\end{equation}
The ordering of the Wilson lines in (\ref{l3}) is such that
$W(z_1,z_2)$ starts first and $W(z_3,z_4)$ follows it. Then next
we let $W(z_2,z_1)$ starts first and $W(z_4,z_3)$ follows it. Let
us call this ordering as the simultaneous ordering which will be
used to define the braiding formulas for a crossing between two
knot components of a link. This ordering has a property that when
the $z$ parameter has traced one loop in one knot of the link we
have that the $z$ parameter has also traced one loop on the other
knot. From this ordering we have that (\ref{l3}) is equal to:
\begin{equation}
\begin{array}{rl}
&Tr RW(z_1,w_1)W(w_1,z_2)W(z_3,w_1)W(w_1,z_4)\cdot\\
& \qquad
 W(z_2,w_2)W(w_2,z_1)W(z_4,w_2)W(w_2,z_3)R^{-1} \\
=&Tr W(z_1,z_2)W(z_3,z_4)
W(z_2,z_1)W(z_4,z_3)
\end{array}
\label{l4}
\end{equation}
Then let us consider the following trace:
\begin{equation}
\begin{array}{rl}
& Tr
R^{-2}W(z_3,w_1)W(w_1,z_2)W(z_1,w_1)W(w_1,z_4)\cdot\\
& \qquad W(z_2,w_2)W(w_2,z_3)W(z_4,w_2)W(w_2,z_1)
\end{array}
\label{l5}
\end{equation}
We let the ordering of the Wilson lines in (\ref{l5})
be such that $W(z_1,z_2)$ starts first and $W(z_4,z_3)$
follows it. Then next $W(z_2,z_1)$ starts first and
$W(z_3,z_4)$ follows it. From this ordering we have that
(\ref{l5}) is equal to:
\begin{equation}
\begin{array}{rl}
&Tr
R^{-2}RW(z_1,w_1)W(w_1,z_2)W(z_3,w_1)W(w_1,z_4)\cdot\\
& \qquad W(z_2,w_2)W(w_2,z_1)W(z_4,w_2)W(w_2,z_3)R \\
=&Tr W(z_1,z_2)W(z_3,z_4)
W(z_2,z_1)W(z_4,z_3)
\end{array}
\label{l6}
\end{equation}
Then since in (\ref{l6}) the crossings between the two knot components have been eliminated we have that the two knot components are independent 
and thus the starting points for the two knot components are independent
and thus  (\ref{l6}) is equal to (\ref{l4}).

On the other hand from the ordering of (\ref{l5}) we
have that (\ref{l5}) is equal to:
\begin{equation}
\begin{array}{rl}
&Tr R^{-2}W(z_3,w_1)RW(z_1,w_1)W(w_1,z_2)R^{-1}\\
& \qquad W(w_1,z_4)W(z_2,w_2)
RW(z_4,w_2)W(w_2,z_3)R^{-1}W(w_2,z_1) \\
=&Tr R^{-2}W(z_3,w_1)RW(z_1,z_2)
R^{-1}W(w_1,z_4)W(z_2,w_2)R
W(z_4,z_3)R^{-1} \\
  & \qquad W(w_2,z_1) \\
=&Tr R^{-2}W(z_3,w_1)RW(z_1,z_2)
W(z_2,w_2)W(w_1,z_4)W(z_4,z_3)R^{-1}W(w_2,z_1) \\
=&Tr R^{-2}W(z_3,w_1)
RW(z_1,w_2)W(w_1,z_3)R^{-1}W(w_2,z_1) \\
=&Tr R^{-2}W(z_3,w_1)
W(w_1,z_3)W(z_1,w_2)W(w_2,z_1)\\
=&Tr R^{-2}W(z_3,z_3)W(z_1,z_1)
\end{array}
\label{l7}
\end{equation}
where we have repeatedly used (\ref{m9}). From (\ref{l4}),
(\ref{l6}) and (\ref{l7}) we have that the knot invariant for the
Hopf link in 
Fig.7a  
is given by:
\begin{equation}
Tr R^{-2}W(z_3,z_3)W(z_1,z_1)
\label{l8}
\end{equation}
We remark that in (\ref{l8}) since $R$ is a $R$-matrix between two
knot components of the Hopf link we have that  $R$ acts on
$W(C_1):=W(z_3,z_3)$ or on $W(C_2):=W(z_1,z_1)$.
In this case we say that the domain of $R$ is $\{W(C_1), W(C_2)\}$.

From this property of $R$ we have that the $R$ and the monodromies
$R_i, i=1,2$ for $W(C_1)$ and $W(C_2)$ in (\ref{l8}) are
independent.

 Then let us consider the Hopf link in 
 Fig.7b. 
 The
correlation for this link is given by
\begin{equation}
\begin{array}{rl}
& Tr W(z_4,w_1)W(w_1,z_2)W(z_1,w_1)W(w_1,z_3)\cdot\\
& \qquad W(z_2,w_2)W(w_2,z_4)W(z_3,w_2)W(w_2,z_1)
\end{array}
\label{l9}
\end{equation}
By a derivation which is dual to the above derivation
for the Hopf link in Fig.7a we have that (\ref{l9})
is equal to
\begin{equation}
Tr R^{2}W(z_4,z_4)W(z_1,z_1)
\label{l10}
\end{equation}
where the $R$ and the monodromies for $W(z_4,z_4)$ and
$W(z_1,z_1)$ in (\ref{l10}) are independent. We see that the
invariants for the above two Hopf links are different. This agrees
with the fact that these two links are not topologically
equivalent.

As more examples let us consider the linking of two trivial
knots with linking number $2$ as in 
Fig.8a (The reader may skip the following of this section for the first reading). 
Similar to the above
computations we have that this link which is analogous to
the Hopf link in Fig.7a 
is with an invariant equals
to $Tr R^{-4}W(z_4,z_4)W(z_1,z_1)$. Also
for the link in 
Fig.8b 
which is analogous to the Hopf link
in Fig.7b 
is with an invariant equals to
$Tr R^{4}W(z_4,z_4)W(z_1,z_1)$.

\begin{figure}[hbt]
\centering
\includegraphics[scale=0.5]{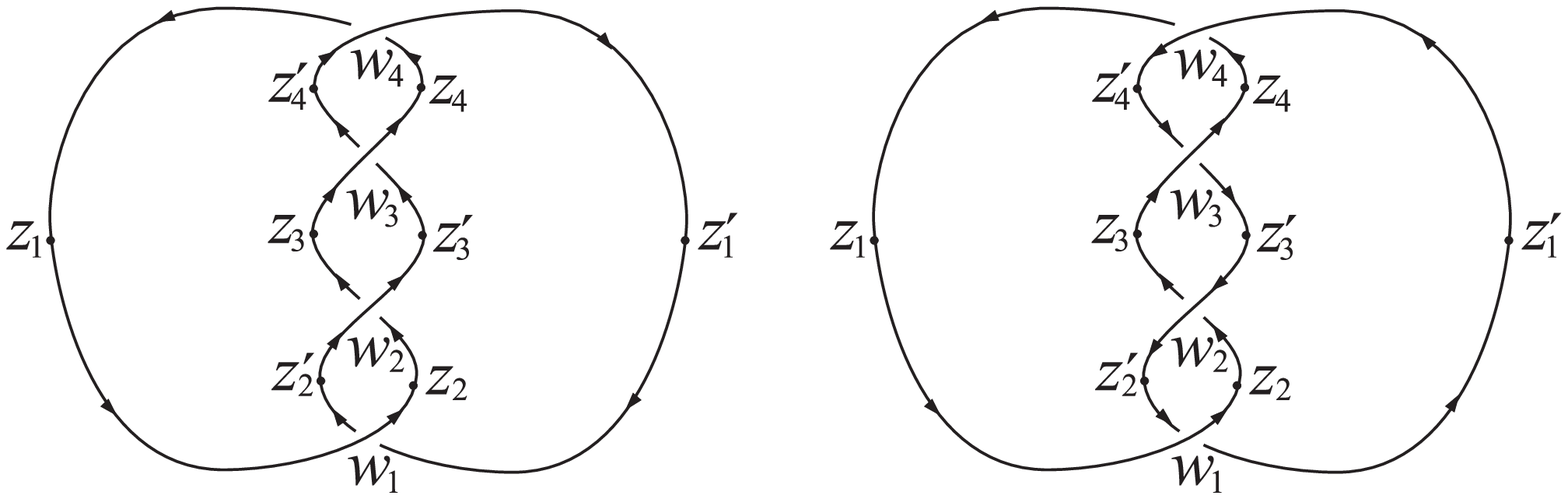}

Fig.8a   
\hspace*{4.5cm} Fig.8b 
\end{figure}

Let us give a computation of the link invariant of the link diagram in
Fig.8a, as follows.
By definition the link invariant of this diagram is the
 trace of the generalized Wilson loop  of this link diagram
which is given by:
\begin{equation}
\begin{array}{rl}
 & Tr W(z_1^{'},w_1)W(w_1,z_2)W(z_1,w_1)W(w_1,z_2^{'})\cdot \\
 & W(z_2,w_2)W(w_2,z_3^{'})W(z_2^{'},w_2)W(w_2,z_3)\cdot \\
  & W(z_3^{'},w_3)W(w_3,z_4)W(z_3,w_3)W(w_3,z_4^{'})\cdot \\
 & W(z_4,w_4)W(w_4,z_1^{'})W(z_4^{'},w_4)W(w_4,z_1)
\end{array}
\label{link}
\end{equation}
where the ordering is such that $W(z_1,z_2)$ stars first and $W(z_1^{'},z_2^{'})$
follows it. Then next we let $W(z_2,z_3)$ stars first and
$W(z_2^{'},z_3)^{'}$ follows it. Continuing in this way we have an
ordering such that when the $z$ parameter has traced one loop
we have that the $z^{'}$ parameter has also traced one loop.
From the ordering  and the braiding formulas (\ref{m7a}), (\ref{m8a})
we have that (\ref{link}) is equal to:
\begin{equation}
\begin{array}{rl}
 & Tr RW(z_1,w_1)W(w_1,z_2)W(z_1^{'},w_1)W(w_1,z_2^{'})\cdot \\
 & W(z_2,w_2)W(w_2,z_3)W(z_2^{'},w_2)W(w_2,z_3^{'})R^{-1}\cdot \\
  & RW(z_3,w_3)W(w_3,z_4)W(z_3^{'},w_3)W(w_3,z_4^{'})\cdot \\
 & W(z_4,w_4)W(w_4,z_1)W(z_4^{'},w_4)W(w_4,z_1^{'}) \\
= & Tr W(z_1,z_2)W(z_1^{'},z_2^{'})\cdot
 W(z_2,z_3)W(z_2^{'},z_3^{'})\cdot \\
  & W(z_3,z_4)W(z_3^{'},z_4^{'})\cdot
  W(z_4,z_1)W(z_4^{'},z_1^{'})
\end{array}
\label{link2}
\end{equation}
On the other hand as similar to the Hopf link let us consider the following
trace:
\begin{equation}
\begin{array}{rl}
 & Tr R^{-2}W(z_1^{'},w_1)W(w_1,z_2)W(z_1,w_1)W(w_1,z_2^{'})\cdot \\
 & W(z_2,w_2)W(w_2,z_3^{'})W(z_2^{'},w_2)W(w_2,z_3)\cdot \\
  & W(z_3^{'},w_3)W(w_3,z_4)W(z_3,w_3)W(w_3,z_4^{'})\cdot \\
 & W(z_4,w_4)W(w_4,z_1^{'})W(z_4^{'},w_4)W(w_4,z_1)
\end{array}
\label{link3}
\end{equation}
where we let the ordering be such that $W(z_1,z_2)$ starts first and
$W(z_4^{'},z_1^{'})$ follows it. Continuing in this way the ordering
in (\ref{link3}) is then determined. From this ordering we have that
(\ref{link3}) is equal to:
\begin{equation}
\begin{array}{rl}
 & Tr R^{-2}RW(z_1,w_1)W(w_1,z_2)W(z_1^{'},w_1)W(w_1,z_2^{'})\cdot \\
 & W(z_2,w_2)W(w_2,z_3)W(z_2^{'},w_2)W(w_2,z_3^{'})R^{-1}\cdot \\
  & RW(z_3,w_3)W(w_3,z_4)W(z_3^{'},w_3)W(w_3,z_4^{'})\cdot \\
 & W(z_4,w_4)W(w_4,z_1)W(z_4^{'},w_4)W(w_4,z_1^{'})R \\
=  & Tr W(z_1,z_2)W(z_1^{'},z_2^{'})\cdot
 W(z_2,z_3)W(z_2^{'},z_3^{'})\cdot \\
  & W(z_3,z_4)W(z_3^{'},z_4^{'})\cdot
 W(z_4,z_1)W(z_4^{'},z_1^{'})
\end{array}
\label{link4}
\end{equation}
Then since the two knot components are independent we have that
the starting points for the two knot components are independent
and we thus have that (\ref{link4}) is equal to (\ref{link2}).

On the other hand from the ordering of (\ref{link3}) we have that
(\ref{link3}) is equal to:
\begin{equation}
\begin{array}{rl}
 & Tr R^{-2}W(z_1^{'},w_1)RW(z_1,w_1)W(w_1,z_2)R^{-1}W(w_1,z_2^{'})\cdot \\
 & W(z_2,w_2)RW(z_2^{'},w_2)W(w_2,z_3^{'})R^{-1}W(w_2,z_3)\cdot \\
  & W(z_3^{'},w_3)RW(z_3,w_3)W(w_3,z_4)R^{-1}W(w_3,z_4^{'})\cdot \\
 & W(z_4,w_4)RW(z_4^{'},w_4)W(w_4,z_1^{'})R^{-1}W(w_4,z_1)\\

=& Tr R^{-2}W(z_1^{'},w_1)RW(z_1,z_2)R^{-1}W(w_1,z_2^{'})\cdot \\
 & W(z_2,w_2)RW(z_2^{'},z_3^{'})R^{-1}W(w_2,z_3)\cdot \\
  & W(z_3^{'},w_3)RW(z_3,z_4)R^{-1}W(w_3,z_4^{'})\cdot \\
 & W(z_4,w_4)RW(z_4^{'},z_1^{'})R^{-1}W(w_4,z_1)\\

= & Tr R^{-2}W(z_1^{'},w_1)RW(z_1,z_2)W(z_2,w_2)W(w_1,z_2^{'})\cdot \\
 & W(z_2^{'},z_3^{'})R^{-1}W(w_2,z_3)\cdot \\
  & W(z_3^{'},w_3)RW(z_3,z_4)W(z_4,w_4)W(w_3,z_4^{'})\cdot \\
 & W(z_4^{'},z_1^{'})R^{-1}W(w_4,z_1)\\

= & Tr R^{-2}W(z_1^{'},w_1)RW(z_1,w_2)W(w_1,z_3^{'})R^{-1}W(w_2,z_3) \\
 & W(z_3^{'},w_3)RW(z_3,w_4)W(w_3,z_1^{'})R^{-1}W(w_4,z_1)\\

= & Tr R^{-2}W(z_1^{'},w_1)W(w_1,z_3^{'})W(z_1,w_2)W(w_2,z_3) \\
 & W(z_3^{'},w_3)W(w_3,z_1^{'})W(z_3,w_4)W(w_4,z_1)\\

= & Tr R^{-2}W(z_1^{'},z_3^{'})W(z_1,z_3) \\
 & W(z_3^{'},z_1^{'})W(z_3,z_1)\\

= & Tr R^{-2}R^{-2}W(z_1^{'},z_1^{'})W(z_1,z_1)
\end{array}
\label{link5}
\end{equation}
where the final step is from the above derivation of the invariant
of the Hopf link. This shows that the invariant of the knot
diagram (a) in Fig.8  
is equal to $Tr
R^{-4}W(z_1^{'},z_1^{'})W(z_1,z_1)$ where $R$ is independent of
the monodromies of $W(z_1^{'},z_1^{'})$ and $W(z_1,z_1)$.
Similarly we can show that the invariant of the knot diagram (b)
in Fig.8  
is equal to $Tr R^{4}W(z_1^{'},z_1^{'})W(z_1,z_1)$.

Let us generalize the Hopf link to the case with linking number $n$.
Then  by induction on the above results we have that the two generalized Hopf links
with linking number $n$
are respectively with invariants
\begin{equation}
Tr R^{-2n}W(z_1^{'},z_1^{'})W(z_1,z_1), \qquad
Tr R^{2n}W(z_1^{'},z_1^{'})W(z_1,z_1)
\label{b}
\end{equation}

\section{Classification of Links}

 Similar to the case of knot for each link $L$
let us construct the generalized Wilson loop $W(L)$. For the case
of link in  constructing the generalized Wilson loop we need to
consider the crossings between two knot components of a link. As
shown in the Hopf link example for a crossing between two knot
components of a link we give it a simultaneous ordering such that
the braiding formulas for such crossing are defined.
When the  braiding formulas are defined we have then completely
represented this crossing by its Wilson product. Once a crossing
between two knot components of a link $L$ is completely
represented by its Wilson product we can then follow the
orientations of each knot component of this link $L$ to write out
the sequence of Wilson products for the sequence of crossings on
each knot component of the link. By writing out all these
sequences of Wilson products of each knot component one by one
until all crossings have been counted we have that the generalized
Wilson loop $W(L)$ of $L$ is then formed. In this process of
counting the crossings we have that the crossings which have been
counted once will not been counted again when they reappear. From
these reappearances we have the property of circling and
sub-circling of the link.

Let us consider some examples to illustrate the construction of $W(L)$.
As a simple example let us consider again the Hopf links in
Fig.7. 
We let an ordering
be such that $W(z_1,z_2)$ starts first and $W(z_3,z_4)$
follows it simultaneously. This is by definition a simultaneous ordering
of $W(z_1,z_2)$ and $W(z_3,z_4)$.
Then next  we let $W(z_2,z_1)$ starts first
and $W(z_4,z_3)$ follows it simultaneously. This is by definition a simultaneous ordering
of $W(z_2,z_1)$ and $W(z_4,z_3)$.

For the Hopf link if we let $1$ denote the
crossing of $W(z_1,z_2)$ with $W(z_3,z_4)$ and let $2$
denote the
crossing of $W(z_2,z_1)$ with $W(z_4,z_3)$.
Then we have $W(L)=12$. 

Let us also denote the corresponding crossings of the Hopf link by $1$ and $2$ respectively.
Then for the Hopf link  we have the following circling
property:
\begin{equation}
12 = 21 = 12 = \cdot\cdot\cdot
\label{o3}
\end{equation}

Then by using exactly the same method for proving the circling property of $W(K)$ of a knot $K$ we can show that the generalized Wilson loop $W(L)=12$ of the Hopf link $L$ also has the circling property (\ref{o3}).

In the following let us consider more examples of $W(L)$ and the circling property of links.

{\bf Examples of $W(L)$ and the circling property of links}.

 As an example
let us consider the link $L$ in Fig.9. 

\begin{figure}[hbt]
\centering
\includegraphics[scale=0.8]{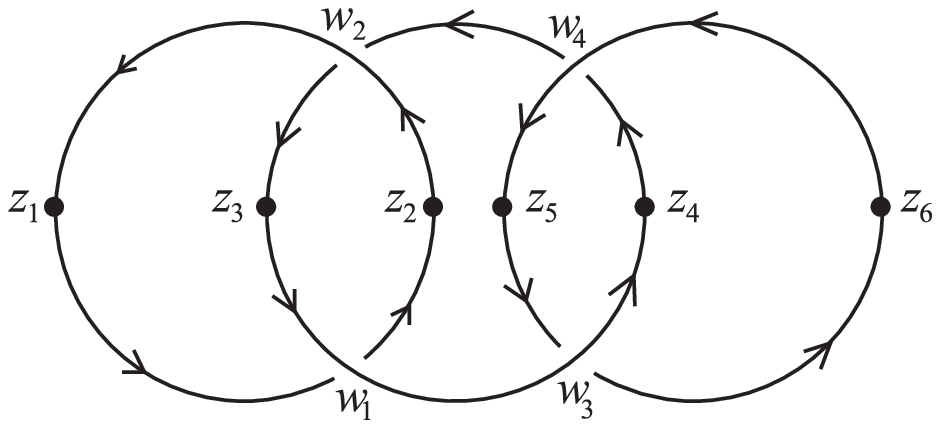}

Fig.9
\end{figure}

Let $i$ denote the
crossing at $w_i, i=1,2,3,4$. Then the generalized
Wilson loop $W(L)$ is given by $W(L)=1234$. Then by using exactly the same method for proving the circling property of $W(K)$ of a knot $K$ we can show that the generalized
Wilson loop $W(L)=1234$ satisfies the following three circlings of
$L$:
\begin{equation}
\begin{array}{rl}
& [12]34 = [21]34 = [12]34 = \cdot\cdot\cdot
\\
& 1234 = 4123 = 3412 = 2341 = 1234 = \cdot\cdot\cdot
\\
& 12[34] = 12[43] = 12[34] = \cdot\cdot\cdot
\end{array}
\label{o4}
\end{equation}
where the sequences in the bracket $[\,\, ]$ form a circling.

As one more example let us consider the Borromeo ring $L$ in
Fig.10.

\begin{figure}[hbt]
\centering
\includegraphics[scale=0.8]{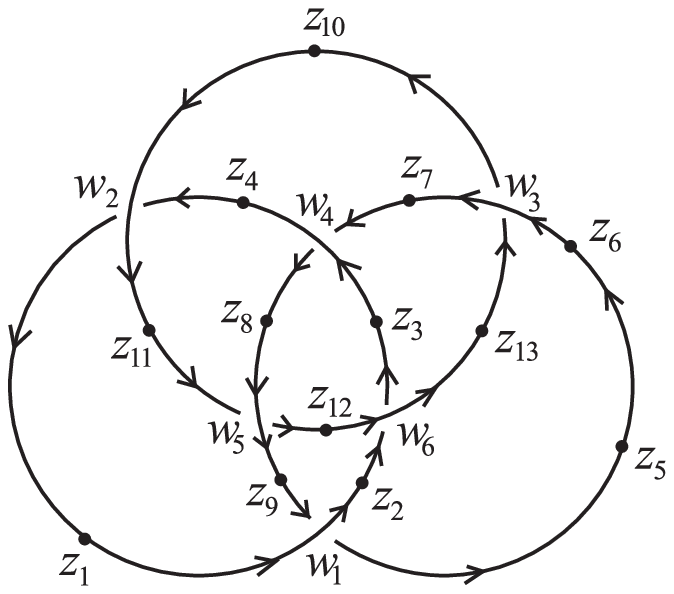}

Fig.10
\end{figure}

Let $i$ denote the crossing at $w_i, i=1,...,6$. Then we
have many ways to write out $L$
which are all equivalent by the circling property of $L$. As an
example we have the following three circlings of $L$ where each
sequence of crossings representing $L$:
\begin{equation}
\begin{array}{rl}
& [1642]5(6)3 = [2164]5(6)3 = [4216]5(6)3 =
[6421]5(6)3=[1642]5(6)3 =\cdot\cdot\cdot
\\
& [2563]4(5)1 = [3256]4(5)1 = [6325]4(5)1 = [5632]4(5)1
=[2563]4(5)1= \cdot\cdot\cdot
\\
& [3451]6(4)2 = [1345]6(4)2 = [5134]6(4)2 = [4513]6(4)2
=[3451]6(4)2 = \cdot\cdot\cdot
\end{array}
\label{o4a}
\end{equation}
In the above sequences the notation $(x)$ means that the number
$x$ is circled to the position of $x$ without $(\,\,)$ as indicated in
the sequences (Also this notation $(x)$ means that the $x$ in $(\,\,)$
reappears and is not counted).

Also the above three circlings of sequences can be
circled to each other. For example we have the following circling:
\begin{equation}
\begin{array}{rl}
&[1642]5(6)3=[1(6)42]563
 = 1(6)4[2563] = 14[2563] =4[2563]1 \\
 =&[2563]14=[2563]1(3)4(5)
 =[2563]4(5)1(3) =[2563]4(5)1=\cdot\cdot\cdot
\end{array}
 \label{o5}
\end{equation}
This shows that the  first and the second circlings of the above
three circlings can be circled to each other. 

Then each of the above sequences can be as the generalized Wilson loop $W(L)$ and by using exactly the same method for the case of knot for proving the cirling property  we can show that the
generalized Wilson loop $W(L)$ also has the above circling properties (\ref{o4}) and (\ref{o5}).
$\diamond$

In general by using exactly the same method for the case of knot for proving the cirling property of the
generalized Wilson loop of a knot we can show that the
generalized Wilson loop of
a general link also has the above circling property of this general link.
 With this circling property of the  generalized
Wilson loop as similar to the case of knot we have the following
theorem for links:
\begin{theorem}
Each link $L$ can be faithfully represented by its
generalized Wilson loop $W(L)$
in the sense that if two link diagrams have the
same generalized Wilson loop then these two link diagrams must be
equivalent.
\end{theorem}

From this theorem on faithful representation of links
we  have the following classification theorem for links.
\begin{theorem}
Let $W(L)$ denote the generalized Wilson loop of a
link $L$ with component knots $K_j, j=1,... n$.
Then  $W(L)$ is a link invariant which classifies links.
We have that $W(L)$ can  be written in the following
form:
\begin{equation}
W(L)= P_L\prod_{i=1}^n W(K_i)
:= R_1^{-m_1}\cdot\cdot\cdot  R_q^{-m_q}
W(K_1)\cdot\cdot\cdot W(K_n)
\label{AA2a}
\end{equation}
where $R_i, i=1,...,q$ are  
monodromies of the KZ equation which come from the linkings of
$K_j, j=1,...,n$ where  the domains of $R_i$ are subsets of $\{W(K_1), ... , W(K_n)\}$.
Also these monodromies $R_i$ and monodromies for $W(K_j), j=1,...,n$ are independent of each other  in the sense that the acting domains of these monodromies are different.

 Then the trace $Tr W(L)$
is also a link invariant which  classifies
links.
\end{theorem}

{\bf Proof.}
The proof of this theorem is similar to the case of
knot. As similar to the case of knot
let us  first find the following expression
for  $ W(L)$:
\begin{equation}
 W(L)=
R_1^{-m_1}\cdot\cdot\cdot  R_q^{-m_q}
 W(K_1)\cdot\cdot\cdot
 W(K_n)
\label{AA3a}
\end{equation}
where $R_i$ are independent monodromies which are also independent
of the monodromies for $W(K_j)$.
  From this expression and the
above theorem on faithful representation of link we
then have that the link invariant
$Tr W(L)$ classifies links.

Let $L$ be a trivial link
with two unlinking component  $K_1$ and $K_2$.
We may suppose that $K_1$ and $K_2$ have no crossings
to each other. Then  $W(L)$ is in the following form:
\begin{equation}
 W(L)= W(K_1) W(K_2)
\label{c5}
\end{equation}
where we have
\begin{equation}
W(K_j) =R_j^{k_j}A_j, \qquad j=1,2
\label{c6}
\end{equation}
for some $k_j, j=1,2$.
We have that the $R_j, j=1,2$ matrices operating on two independent operators
$A_1$ and $A_2$ respectively.

Let $K_1$ and $K_2$ be linked to form a link $L$.
Then from the theorem in the
section on solving the KZ equations we have that $W(L)$
is in a tensor product form. Since $K_1$ and $K_2$ are two closed curves
we have that this tensor product reduces its degree to a product with a
tensor product of the form
 $A_1\otimes A_2$ where $A_1$ and $A_2$ are two independent operators for
 $K_1$ and $K_2$ respectively. Then since the matrices $\Phi_{ij}$
 and $\Psi_{ij}$ act on either $A_1$ or $A_2$ we have that they commute
with $A_1\otimes A_2$ and thus we can write $W(L)$ in the following
form:
\begin{equation}
W(L)=R_{a_1}^{-p_{a_1}}\cdot\cdot\cdot R_{a_b}^{-p_{a_b}}A_1\otimes A_2
\label{CL}
\end{equation}
where the monodromies $R_{a_i}$ acts either on  $A_1$ or $A_2$.
Since $W(K_j)=R_j^{-k_j}A_j, j=1,2$ for some $R_j$ from (\ref{CL}) we can
write $W(L)$ in the following form:
\begin{equation}
W(L)=R_{1}^{-m_{1}}\cdot\cdot\cdot R_{q}^{-m_{q}}W(K_1)W(K_2)
\label{CL1}
\end{equation}
where $W(K_1)W(K_2)=W(K_1)\otimes W(K_2)$.
We have that the monodromies $R_{i}$ in (\ref{CL1}) must be
independent of the monodromies $R_j, j=1,2$ since if $R_{i}=R_j$
then it will be absorbed by $W(K_j)$ to form a generalized
Wilson loop $W(K')$ for some knot $K'$ which is not equivalent to $W(K_j)$.
This is impossible since $L$ is not formed with this knot $K'$.
On the other hand the monodromies $R_{i}$ in (\ref{CL1}) can be set to
be independent of each other (in the sense that the acting domains of these monodromies are different)
since if two $R_{i}$ are the same then
they can be merged into one $R_{i}$.
This form (\ref{CL1}) of $W(L)$ is just the required form (\ref{AA3a}).
For a general $L$ the proof of this form of $W(L)$ is similar.
This proves the theorem. $\diamond$

{\bf Remark}. Let us give
more details on the domain of a monodromy, as follows.
For simplicity let us consider the above link $L$ with two components
$K_1$ and $K_2$. We have that
a monodromy  $R_{i}$ acts on $A_1$ or $A_2$ (or acts on $W(K_1)$ or
$W(K_2)$. Thus the domain of $R_{i}$ is actually a subset of
$\{W(K_1), W(K_2)\}$.   
Let us consider the Hopf link as a simple example. For the Hopf link $L$ we have that
$W(L)= R^{\pm 2}W(K_1)W(K_2)$. In this $W(L)$ the monodromy $R$
is with domain  $\{W(K_1), W(K_2)\}$ since $R$ is obtained by braiding
between the Wilson lines of $K_1$ and  $K_2$. On the other hand the
$R_j, j=1,2$ are for the forming of $W(K_j)$. In this case we say that $R_j$ are with domains  $\{W(K_j)\}, j=1,2$
respectively.

\section{Quantum Invariants of 3-manifolds and Classification}\label{sec13a}

In this section we derive quantum invariants of closed 3-manifolds
from the above quantum invariants of links.
We have the Lickorish-Wallace theorem  which states that any closed (oriented and
connected) 3-manifold $M$ can be obtained from a Dehn surgery on a
framed link $L$ \cite{Lic1}\cite{Rol}\cite{Wal}\cite{Pra}.

Let us first consider 3-manifolds obtained from surgery on 
framed
knots $K^{\frac pq}$ where $p$ and $q$ are co-prime integers.  We
have the following expression of the generalized Wilson loop
$W(K^{\frac pq})$ of $K^{\frac pq}$:
\begin{equation}
\begin{array}{rl}
W(K^{\frac pq})&=R^{-2p}
R_{3}^{m_{3}}W(K)
R_{3}^{-m_{3}}W(K_c)\\
 &=R^{-2p}
 R_{3}^{m_{3}}R_1^{-m_1}W(C_1)
 R_{3}^{-m_{3}}R_2^{-m_2}W(C_2)
\end{array} 
\label{q11a}
\end{equation}
where the $R$-matrix $R$ denotes the linking matrix which acts on
$(W(K),W(K_c))$ where $K_c$ denotes the partner (or company) of $K$
for the framed knot $K^{\frac pq}$. 
In (\ref{q11a}) we write
$W(K)=R_1^{-m_1}W(C_1)$  and $W(K_c)=R_2^{-m_2}W(C_2)$. 
The integers
 $m_1$ and $m_2$ are indexes for the knot $K$ and its partner $K_c$
respectively.

On the other hand  in (\ref{q11a}) the $R$-matrix $R_{3}$ acting on $W(K)$ and $W(K_c)$ is from the linking of $K$ and $K_c$ (and from the number $q$) and is for the effect of giving $0$ linking number. 

Then since $K_c$ is as the partner of $K$ in the construction of $M$ 
we have that the $R$-matrices $R_i$  acting
on $W(C_i)$ $(i=1,2)$ respectively are such that $R_1=R_2$. 
From this  we have that $R_3$ and $R_1=R_2$ are as the same function on $W(C_1)$ and $W(C_2)$ and thus $R_1=R_2=R_3$. Thus from (\ref{q11a}) we have the following representation of $M$:
\begin{equation}
W(K^{\frac pq})=R^{-2p}W(K)W(K_c)
 =R^{-2p}R_1^{-m_1}W(C_1)R_2^{-m_2}W(C_2)
\label{q11}
\end{equation}
where we have absorbed  $R_{3}$ into the matrices $R_i$ $(i=1,2)$ and for simplicity the resulting indexes are still denoted by $m_1$ and $m_2$ respectively (For simplicity we still use $W(K)$ to denote $R_1^{-m_1}W(C_1)$ for the resulting index $m_1$ which may be different from the original $m_1$ for $W(K)$. Similarly we still use $W(K_c)$ to denote $R_2^{-m_2}W(C_2)$ for the resulting index $m_2$ which may be different from the original $m_2$ for $W(K_c)$).

Thus in the case of the linking of $K$ and $K_c$ giving the effect of $0$ linking number there may have many surgeries on different $K$ (with different original $m_1$ and $m_2$ but giving the same resulting indexes $m_1$ and $m_2$) giving the same $M$ by this degeneration \cite{Osi}. We notice that all these surgeries are with the same representation (\ref{q11}) of $M$.

Then from Kirby calculus \cite{Kir} we have that (\ref{q11}) may still be a many-to-one
representation of 3-manifold $M$ obtained from surgery on
$K^{\frac pq}$. Let us from (\ref{q11}) find a one-to-one
representation (or invariant) of 3-manifold $M$. To this end let
us investigate homeomorphisms (or symmetries) on $M$ which are the origins of many-to-one of (\ref{q11}). Since (\ref{q11}) is a
surjective representation of 3-manifold $M$ (in the sense that if
$M_1$ is not homeomorphic to $M_2$ then their representations
(\ref{q11}) are not equal) we can investigate these symmetries
from the form of (\ref{q11}).

From (\ref{q11}) we see that there are three independent degrees
of freedom from the three indexes of the three $R$-matrices: $R$,
$R_1$ and $R_2$.  
Then we notice that there is a degenerate degree
of freedom between $K$ and $K_c$ that  
$R_1=R_2$ while $R_1$ and $R_2$ act on $W(C_1)$ and
$W(C_2)$ respectively. 
Further this degenerate degree of freedom
is the only degenerate degree of freedom of (\ref{q11}).   
This 
degenerate degree of freedom of (\ref{q11}) is the source of all
symmetries on $M$ when $M$ is represented by (\ref{q11}) since (\ref{q11}) is a surjective representation of $M$ that it contains all the topological properties of $M$ and that any nontrivial symmetry of $M$ when $M$ is represented by (\ref{q11}) needs a degenerate degree of freedom in (\ref{q11}) for its existence.

Let us consider a symmetry from this degenerate degree of freedom.
If we reverse the orientation of $K_c$ then we have a symmetry
that the obtained manifold is the same manifold as the original
one. Under this symmetry we have that
 $p$ is changed to $-p$. Let us call this symmetry as the
$\pm p$-symmetry.

Then there is a symmetry between $K$ and its mirror
image $\overline{K}$ that surgery on $K^{\frac pq}$ and on
$\overline{K}^{\frac{-p}{q}}$  give the same 3-manifold. This
symmetry is also from the changing of $p$ to $-p$ and can be
regarded as a part of the whole $\pm p$-symmetry. Then when $K$ is an amphichieral knot that $K=\overline{K}$ there is a further degeneration that 
$K^{\frac pq}=K^{\frac{-p}{q}}$.

Let us consider the $\pm p$-symmetry. By the
$\pm p$-symmetry we have that surgery on $K^{\frac pq}$ and on
 $K^{\frac {-p}{-q}}$ give the same manifold. Let us consider the
 details of this symmetry step by step. Let us first consider the
 case that $q=1$. Let $K$ be a nontrivial knot. Then by the
 $\pm p$-symmetry we have
that $W(K^{p})=R^{-2p}W(K)W(K_c)=R^{-2p}W(K)W(K)$ where $K=K_c$
represents the same 3-manifold as that of
$W(K^{\frac{-p}{-1}})=R^{2p}W(K)W(K_r)=R^{2p}W(K)W(K_{cr})$ where
$K_r$ denotes the knot which is obtained from $K$ by reversing the
orientation of $K$ (We have $K_{cr}=K_r$ since $K=K_c$ in this
case). This reversing of orientation is from the $\pm p$-symmetry
that $q$ is changed to $-q$. Let us then consider the structure of
$W(K_r)$. Since $W(K)=R_1^{-m_1}W(C_1)$ and that a general form of
$R$-matrix is of the form $R_1^a$ for some integer $a$ and for
some $R$-matrix $R_1$ we have that $W(K_r)$ is of the following
form:
\begin{equation}
W(K_r)=(R_2^a)^{-m_1}W(C_2) \label{q13}
\end{equation}
where we let $R_2=R_1$ and the $R_1$ of $W(K)=R_1^{-m_1}W(C_1)$ is
replaced by $R_1^a=R_2^a$ for some integer $a$ which is as a new
variable to be determined and we let $C_2$ be a copy of $C_1$ with
reversing orientation and $R_i$ acts on $W(C_i)$ for $i=1,2$.

We notice that  we now have the vector $(m_1, am_1)$ which is as the
index vector for the $R_1$ and $R_2$ matrices where the integer
$a$ is as a new degree of freedom when $m_1\neq 0$ ($m_1\neq 0$ corresponds to a nontrivial knot). This implies
that the integer $am_1$ is a new degree of freedom when $m_1\neq 0$.
Thus we have that
from the $\pm p$-symmetry a new degree of freedom is introduced
and this completely eliminates the property of degenerate degree
of freedom of (\ref{q11}).

Now since there are no more degenerate degree of freedom left in
the form (\ref{q13}) we have that the $\pm p$-symmetry is the only
nontrivial symmetry which can be derived from (\ref{q11}) when
$m\neq 0$ where by the term nontrivial symmetry we mean a symmetry
which can transform a form of (\ref{q11})
to another distinct form of (\ref{q11}).

Let us now determine the integer $a$ for $W(K_r)$ for a given $p$.
To this end let us consider some consequences of the $\pm
p$-symmetry, as follows.

From the $\pm p$-symmetry we have the degenerate property that
$R^{-2p}W(K)W(K)$ and $R^{2p}W(K)W(K_r)$ represent the same
manifold. We have $W(K)W(K)= R_1^{-m_1}W(C_1)R_2^{-m_1}W(C_2)$ and
$W(K)W(K_r)= R_1^{-m_1}W(C_1)(R_2^a)^{-m_1}W(C_2)$. For simplicity let
us sometimes omit the factor $W(C_1)W(C_2)$. Then we have that the
two distinct products $R^{-2p}R_1^{-m_1}R_2^{-m_1}$ and
$R^{2p}R_1^{-m_1}(R_2^a)^{-m_1}$ ($a\neq 1$) represent the same
manifold. This is as the degenerate property of the $\pm
p$-symmetry.

We have that the factor $R_1^{-m_1}=R_2^{-m_1}$ of
$R^{-2p}R_1^{-m_1}R_2^{-m_1}$ has already been in
the representation $R^{2p}R_1^{-m_1}(R_2^a)^{-m_1}$. 
Thus the representation
$R^{2p}R_1^{-m_1}(R_2^a)^{-m_1}$ contains the information of the
representation $R^{-2p}R_1^{-m_1}R_2^{-m_1}$ and the information of all
degenerate properties of the $\pm p$-symmetry. Thus we may use the
representation $R^{2p}R_1^{-m_1}(R_2^a)^{-m_1}$ only to represent the
manifold and that this representation $R^{2p}R_1^{-m_1}(R_2^a)^{-m_1}$
contains the information of all degenerate property of the $\pm
p$-symmetry about the manifold.

Now suppose we have another representation
$R^{2p}R_1^{-k}(R_2^b)^{-k}$ where $k\neq m_1$ and $b\neq 1$. Then
since $k\neq m_1$ there are at least three distinct integers from
the set ${m_1,am_1,k,bk}$ and thus this is over the maximal degenerate
property of the $\pm p$-symmetry represented by
$R^{2p}R_1^{-m_1}(R_2^a)^{-m_1}$ which has the maximal degenerate
property of the $\pm p$-symmetry of allowing at most two distinct
integers $m_1,am_1$. It then follows that the two representations
$R^{2p}R_1^{-m_1}(R_2^a)^{-m_1}$ and $R^{2p}R_1^{-k}(R_2^b)^{-k}$
represent two nonequivalent 3-manifolds for $k\neq m_1$.

Further since the $\pm p$-symmetry is the only symmetry of
(\ref{q11}) and the matrix $R^{2p}$ as a function of $p$ is a
one-to-one mapping we have that the two representations
$R^{2p}R_1^{-m_1}(R_2^a)^{-m_1}$ and
$R^{2p^{\prime}}R_1^{-k}(R_2^b)^{-k}$ represent two nonequivalent
3-manifolds where $a,b\neq 1$ when $p\neq p^{\prime}$ where $p$
and $p^{\prime}$ are of the same sign. In summary we have the
following theorem:
\begin{theorem}
The representation $R^{2p}R_1^{-m_1}(R_2^a)^{-m_1}$
indexed by the integers $m_1, am_1$ where $a\neq 1$ represents
3-manifolds in a one-to-one way in the sense that if $k\neq m_1$
then $R^{2p}R_1^{-m_1}(R_2^a)^{-m_1}$ and $R^{2p}R_1^{-k}(R_2^b)^{-k}$
represent two nonequivalent 3-manifolds where $a,b\neq 1$.

Further we have that the two representations
$R^{2p}R_1^{-m_1}(R_2^a)^{-m_1}$ and
$R^{2p^{\prime}}R_1^{-k}(R_2^b)^{-k}$ represent two nonequivalent
3-manifolds where $a,b\neq 1$ when $p\neq p^{\prime}$ where $p$
and $p^{\prime}$ are of the same sign.
\end{theorem}

Now let us determine the property of the number $a$. 
We have that $a$ always exists since it is for the representation $W(K_r)$ of $K_r$. 
Let us consider the representation
$R^{2p}R_1^{-m_1}(R_2^a)^{-m_1}$ for $W(K^{\frac{-p}{-1}})$
where 
we consider the case that $a\neq 1$.
We want to find out the
property of $a$ and the relation between $a$ and $m_1$. 
 
For this $a$ as similar to the role of $m_1$ let us also
construct a product $R^{2p}R_1^{-a}(R_2^d)^{-a}$. Then when $a\neq
m_1$ and $m_1\neq 1$ and $d\neq 1$ by the above theorem we have that
the two products $R^{2p}R_1^{-m_1}(R_2^a)^{-m_1}$ and
$R^{2p}R_1^{-a}(R_2^d)^{-a}$ cannot represent the same 3-manifold.
Thus for $R^{2p}R_1^{-m_1}(R_2^a)^{-m_1}$ and
$R^{2p}R_1^{-a}(R_2^d)^{-a}$ represent the same 3-manifold we must
have either $a=m_1$ when $m_1\neq 1$ (This implies that $d=a=m_1$) or
$a=-m_1=-1$ when $m_1=1$ (This implies $d=a=-m_1$) or $d=1$.

For the case $d=1$ we have that the product
$R^{2p}R_1^{-a}(R_2^d)^{-a}$ equals $R^{2p}R_1^{-a}R_2^{-a}$ which
represents the framed knot $H^p$ for an amphichieral knot $H$ with
the property that $H_r=H$ (and $\overline{H}=H$) and that
$W(H^p)=R^{-2p}W(H)W(H)=R^{2p}R_1^{-a}R_2^{-a}$ and
$W(H^{\frac{-p}{-1}})=W(H^{\frac{-p}{1}})=R^{2p}R_1^{-a}R_2^{-a}$
represent the same manifold. For this amphichieral knot $H$ we
have that the representation $R^{2p}R_1^{-a}R_2^{-a}$ contains
only one integer $a$ and thus its information is contained in
$R^{2p}R_1^{-m_1}(R_2^a)^{-m_1}$ which contains two integers $m$ and
$a$ and thus $R^{2p}R_1^{-a}R_2^{-a}$ and
$R^{2p}R_1^{-m_1}(R_2^a)^{-m_1}$ represent the same manifold where the
representation $R^{2p}R_1^{-m_1}(R_2^a)^{-m_1}$ is with the maximal
non-degenerate property  in the sense that it has the index vector
$(m_1,am_1)$ where $m_1\neq am_1$ and $m_1, am_1\neq 0$ such that no more
degenerate degree of freedom left in this representation.

For the case $a=m_1$ we have the representation
$R^{2p}R_1^{-m_1}(R_2^{m_1})^{-m_1}$ for $K^{\frac{-p}{-1}}$ which
represents the same manifold as that of $R^{-2p}R_1^{-m_1}R_2^{-m_1}$
for $K^p$.

For the case $a=-m_1=-1$ when $m_1=1$ and  $d=a$  
the representation
$R^{2p}R_1^{-m_1}(R_2^{-m_1})^{-m_1}=R^{2p}R_1^{-1}R_2^{1}$ for
$K^{\frac{-p}{-1}}$ which represents the same manifold as that of
$R^{-2p}R_1^{-m_1}R_2^{-m_1}=R^{-2p}R_1^{-1}R_2^{-1}$ for $K^p$.

These three cases  then determine the property of $a$ and its
relation with $m_1$. We have that the case $d=1$ corresponds to an
amphichieral knot $H$. 

On the other hand each amphichieral knot $H$ gives the representation $W(H^{\frac{-p}{-1}})=W(H^{\frac{-p}{1}})=R^{2p}R_1^{-a}R_2^{-a}$
of a manifold $M$ which is a degenerate result of the $\pm p$-symmetry. This degeneratation is as a part of a whole $\pm p$-symmetry. Thus each amphichieral knot $H$ gives a different nontrivial homepmorphism.
Thus from the above analysis there must exist a framed knot $K^{\frac{p}{q}}$ 
which gives the same manifold $M$ 
where more generally we let $q\geq 1$ for some $q$.
When $q=1$ from the above analysis we have that $K$ is non-amphichieral (and thus $K\neq H$)  and $M$ is
with the representation
$R^{2p}R_1^{-m_1}(R_2^a)^{-m_1}$ where $m_1$ is the integer indexing $K$. Then when $q> 1$ we have that $K^{\frac{p}{q}}$ is represented by $R^{-2p}R_1^{-m_1}R_2^{-b}$ where $b$ is the integer indexing $K_c$.
If $K_c=K$ (or $b=m_1$) then we have the same result as the case $q=1$ that $K$ is non-amphichieral and $M$ is with the representation $R^{2p}R_1^{-m_1}(R_2^a)^{-m_1}$. If $K_c\neq K$ (or $b\neq m_1$) then from the above analysis we have that $M$ is with the representation $R^{2p}R_1^{-m_1}(R_2^{m_1})^{-b}$. 
Then from the above analysis we must have $a= m_1$ or $a=b$ (i.e. $H=K$ or $H=K_c$). If $a= m_1$ then $K=H$ and 
$M$ is with the representation
$R^{2p}R_1^{-a}(R_2^a)^{-b}=R^{2p}R_1^{-a}(R_2^b)^{-a}$. Comparing to the representation $R^{2p}R_1^{-a}R_2^{-a}$  we have that $b=1$. Thus if $H=K$ the framed knot $K^{\frac{p}{q}}$ gives no new information to eliminate the degeneration. Thus we must have $H=K_c$ (or $a=b$). Thus $M$ is with the representation $R^{2p}R_1^{-m_1}(R_2^{m_1})^{-a}$. 
Thus we have the following theorem:
\begin{theorem}
Let $M$ be a closed (oriented and connected) 3-manifold which is constructed by a Dehn surgery on a framed knot $K^p$ (or on $K^{\frac{-p}{-1}}$) where $K$ is a nontrivial knot. Then $M$ can be uniquely represented by a representation of the form $R^{2p}R_1^{-m_1}(R_2^a)^{-m_1}$ where $m_1$ and $a$ are integers and $a\neq 1$.

Further when $M$ is not obtained from surgery on $H^p$ where $H$ is an amphichieral knot we have that $a=m_1$. On the other hand when $M$ is obtained from an amphichieral knot $H$ we have that the integer $a$ is from the representation $R_1^a$ of the amphichieral knot $H$ and the integer $m_1$ is from the representation $R_1^{-m_1}$ of another knot $K$ where $M$ is also obtained from surgery on $K^{\frac{p}{q}}$ (or on $K^{\frac{-p}{-q}}$) for some $q\geq 1$ by the $\pm p$-symmetry (when $q > 1$ we have $H=K_c$ where $K_c$ denotes the partner of $K$ for $K^{\frac{p}{q}}$).

Furthermore we have that the two representations $R^{2p}R_1^{-m_1}(R_2^a)^{-m_1}$ and $R^{2p^{\prime}}R_1^{-k}(R_2^b)^{-k}$ represent two nonequivalent 3-manifolds where $a,b\neq 1$ when $p\neq p^{\prime}$ where $p$ and $p^{\prime}$ are of the same sign.
\end{theorem}

We remark that by the $\pm p$-symmetry we may fix the sign of $p$
(For example we may fix the sign of $p$ such that $p>0$) to obtain
the manifold $M$ if $M$ is obtained from $K^p$ or $K^{-p}$.

We shall also write the above invariant
$R^{2p}R_1^{-m_1}(R_2^a)^{-m_1}$ in the following  complete form:
\begin{equation}
\overline{W}(K^{\frac p1}):=R^{2p}R_1^{-m_1}(R_2^a)^{-m_1}W(C_1)W(C_2)
\label{in1}
\end{equation}

Let us then consider the case that $M$ is obtained from
surgery on $K^{\frac{p}{q}}$ where $K$ is a nontrivial knot and
$q> 1$ is an integer which is co-prime with respect to $p$.
From the $\pm p$-symmetry we have that $K^{\frac{p}{q}}$ and
$K^{\frac{-p}{-q}}$ give the same manifold. Then we have the
representations $W(K^{\frac{p}{q}})=R^{-2p}W(K)W(K_c)$ and
$W(K^{\frac{-p}{-q}})=R^{2p}W(K)W(K_{cr})$ for the same manifold.
We write
$W(K_c)=R_2^{-m_2}W(C_2)$ where $m_2$ is the integer indexing
$K_c$ and is related to the integers $q$ and $m_1$ where $m_1$ is the
integer indexing $K$. 

Then if $m_2=0$ we have that $M$ is a lens space and let us consider this subcase later.

For the subcase $m_2\neq 0$ suppose that there exists an amphichieral knot $H$ related to $K$ for $q> 1$ as above. Then by the above theorem $M$ is with the representation
$R^{2p}R_1^{-m_1}(R_2^{m_1})^{-a}$ where $a=m_2$ and $a$ is the
integer indexing $H$. On the other hand suppose that there does not exist an amphichieral knot $H$ related to $K$ for $q> 1$ as above. Then by following the above analysis for $q= 1$ we have that $M$ is simply with the representation $R^{2p}R_1^{-m}(R_2^{m_1})^{-a}$ where $a=m_2$. 

Then as similar to the above case $q= 1$ 
since $\pm p$-symmetry is the only symmetry (when $m_1\neq 0$) and the matrix $R^{2p}$ as a function of $p$ is a one-to-one mapping we have that two representations $R^{2p}R_1^{-m_1}(R_2^{m_1})^{-a}$ and $R^{2p^{\prime}}R_1^{-k}(R_2^k)^{-b}$ represent two nonequivalent 3-manifolds  when $p\neq p^{\prime}$ where $p$ and $p^{\prime}$ are of the same sign.

Thus as similar to the above theorem  we have the following theorem:
\begin{theorem}
Let $M$ be a closed (oriented and connected) 3-manifold which is constructed by a Dehn surgery on a framed knot $K^{\frac{p}{q}}$ (or on $K^{\frac{-p}{-q}}$) where $K$ is a nontrivial knot and
$q\neq 1$ and $M$ is not a lens space. 
Then $M$ can be uniquely represented by a representation of the form $R^{2p}R_1^{-m_1}R_2^{-am_1}$ where $m_1\neq 0$ is an integer for the representation $W(K)=R_1^{-m_1}W(C_1)$ of $K$ and $a=m_2\neq 0$ 
is an integer for the representation $W(K_c)=R_2^{-m_2}W(C_2)$ of $K_c$ .

Further we have that the two representations $R^{2p}R_1^{-m_1}(R_2^{m_1})^{-a}$ and $R^{2p^{\prime}}R_1^{-k}(R_2^k)^{-b}$  represent two nonequivalent 3-manifolds  when $p\neq p^{\prime}$ where $p$ and $p^{\prime}$ are of the same sign.
\end{theorem}

From the above  
theorems we have the following theorem of
one-to-one representation of 3-manifolds obtained from framed
knots $K^{\frac{p}{q}}$:
\begin{theorem}
Let $M$ be a closed (oriented and connected) 3-manifold which is
constructed by a Dehn surgery on a framed knot $K^{\frac{p}{q}}$
 where $K$ is a nontrivial knot and $M$ is not a lens space.
 Then we have the following one-to
 one representation (or invariant) of $M$ (We call this invariant
 as the invariant with the maximal non-degenerate property):
\begin{equation}
R^{2p}R_1^{-m_1}R_2^{-am_1} \label{qq1}
\end{equation}
where $m_1\neq 0$ is the integer indexing $K$  
and
$a=m_2\neq 0$ such that  $a\neq 1$
is an integer related to $m_1$ and $q$ such that $a$ is either the integer indexing an amphichieral knot $H$ giving the same $M$ by Dehn surgery on $H^p$ or is the integer indexing the knot $K_c$ of $K$; 
and we
choose a convention that $p>0$.
\end{theorem}

{\bf Proof.} From the above two 
theorems we have that $M$
constructed by a Dehn surgery on a framed knot $K^{\frac{p}{q}}$
can be represented by (\ref{qq1}) where the expression
$R_1^{-m_1}R_2^{-am_1}$ absorbs all the degenerate amphichieral cases which are the only degenerate cases when $m_1\neq 0$, $m_2\neq 0$ and $p>0$. Thus (\ref{qq1}) is a one-to-one representation of $M$. This proves the theorem. $\diamond$

We remark that we shall also write the above invariant
$R^{2p}R_1^{-m_1}R_2^{-am_1}$
in the
following complete form:
\begin{equation}
\overline{W}(K^{\frac{p}{q}}):=R^{2p}R_1^{-m_1}R_2^{-am_1}W(C_1)W(C_2)
\label{in2}
\end{equation}

{\bf Remark}. There exist nontrivial knots $K$ such that the manifold $M$ obtained from
$K^{\frac {p}{q}}$ with $m_2=0$ is a lens space \cite{Rol}\cite{Ber}. 
 Before the investigation of the case of lens space let us first consider a well-known example
of the above $\pm p$ symmetry from \cite{Rol} and \cite{Pra}.
 $\diamond$

{\bf Example.} Let $K_{RT}^{-1}$ denote the right trefoil knot
$K_{RT}$ with framing $-1$ and let $H_E^{+1}$ denote the
figure-eight knot $H_E$ with framing $+1$. Then as shown in
\cite{Rol} we have that surgery on $K_{RT}^{-1}$ and on $H_E^{+1}$
give the same 3-manifold $M$. Then as a part of the $\pm p$
symmetry ($p=1$) we have that $K_{LT}^{+1}$ gives the same
manifold $M$ as that of $K_{RT}^{-1}$ where $K_{LT}$ denotes the
left trefoil knot. Then since $H_E$ is an amphichieral knot which
is equivalent to its mirror image by the same reason we have that
$H_E^{-1}$ gives the same manifold $M$ as that of $H_E^{+1}$. Let
us investigate this example to illustrate the above invariant (or
representation) for this manifold $M$. We have that the
generalized Wilson loop for $K_{RT}^{-1}$ is given by
\begin{equation}
W(K_{RT}^{-1})=R^{2}R_1^{-1}R_2^{-1}W(C_1)W(C_2) \label{ex1}
\end{equation}
Similarly the generalized Wilson loop for $K_{LT}^{+1}$ is given
by
\begin{equation}
W(K_{LT}^{+1})=R^{-2}R_1^{1}R_2^{1}W(C_1)W(C_2) \label{exa}
\end{equation}
 Then the generalized Wilson loop for $H_E^{+1}$ is given
by
\begin{equation}
W(H_E^{+1})=R^{-2}R_1^{-3}R_2^{-3}W(C_1)W(C_2) \label{ex2}
\end{equation}
where the index of $H_E$ is $3$.  Then as shown in the above
construction of invariant since $H_E^{+1}$ and $H_E^{-1}$ give the
same manifold we have that the index number $3$ for $H_E$ is
absorbed to the generalized Wilson loop
$\overline{W}(K_{LT}^{+1})$:
 \begin{equation}
\overline{W}(H_E^{\pm 1})=\overline{W}(K_{RT}^{-1})
=\overline{W}(K_{LT}^{+1}):=R^{2}W(K_{LT})W(K_{LTcr})
=R^{2}R_1^{1}R_2^{3}W(C_1)W(C_2) \label{ex3}
\end{equation}
This generalized Wilson loop is then by definition the unique
invariant (or representation) for the manifold $M$ constructed by
$K_{RT}^{-1}$,
 $H_E^{\pm 1}$ and $K_{LT}^{+1}$. We notice that for this example we have that
$a=3$ which is the index of $H_E^{\pm 1}$ and thus the indexes of
$R_1$ and $R_2$ are different. This is the  maximal non-degenerate
property of the invariant (\ref{ex3}) in the sense that it
contains the two indexes $1$ and $3$ for $K_{RT}$ and $H_E^{\pm
1}$. $\diamond$

Let us then consider the case $m_2=0$ for lens spaces.
We have that all lens spaces can be constructed by framed knots of the form  $C^{\frac{p}{q^{\prime}}}$ where $C$ denotes a trivial knot. 
Then we have that $m_2=0$ represents a trivial knot. 
Then $m_1 \neq 0$ represents a nontrivial knot. 
Thus the representation (\ref{qq1}) can represent a lens space with linking number $p$ when $m_2$ is related to $m_1 $ such that $m_2=0$. Thus the representation (\ref{qq1}) gives a 
 representation of 3-manifolds $M$ including all the lens spaces.

Let us then determine the number $m_1$ which can give $m_2=0$ (and thus the knot $K$ indexed by $m_1$ can give lens space). 

For $p=0$ we have that $m_1=0$ gives $m_2=0$ and thus gives the lens space $S^2\times S^1$.

For $p=1$ the case $m_1=0$ is excluded since the unknot $C$ with framing $p=1$ can be deleted (and thus is not minimal where we shall give details on the concept of minimal link). Then the $S^3$ is represented by the constant $1$ (and is not represented by $W(C^1)$).

Let us then consider $p>1$. From the property of lens spaces we have that $C^{\frac{p}{q^{\prime}}}$ and $C^{\frac {p}{-q^{\prime}}}$ give the same lens space. This symmetry (or homeomorphism) can be described by the following relation \cite{Mur}:
\begin{equation}
q_1q_2=\pm 1 +np
\label{lens}
\end{equation}
for some integer $n$; and $1\leq q_1 \leq p-1$ where $q^{\prime}=q_1+n_1p$ for some integer $n_1$ and we choose the region (mod $p$) of $q_2$ such that $q_1<q_2 $.
On the other hand since $C$ is a (trivial) amphichieral knot from the above analysis on amphichieral knot we have that there exists a nontrivial knot $K$ indexed by $m_1\neq 0$ (and an integer $q\geq 1$) such that $W(K^{\frac{p}{q}})=R^{-2p}W(K)W(K_c)=R^{-2p}R_1^{-m_1}W(C_1)W(C_2)$ and  $\overline{W}(K^{\frac{p}{-q}}):=R^{2p}R_1^{-am_1}W(C_1)W(C_2)$ represent the same lens space $M$ (which is constructed by $C^{\frac{p}{q^{\prime}}}$) where $m_1$ is replaced by $am_1$ for some integer $a$ such that $a$ is related to $C^{\prime}$ giving the reversing of $q^{\prime}$ and $-q^{\prime}$ as described by (\ref{lens}). Thus from (\ref{lens}) we have that $am_1=q_1q_2$. Thus $a=q_1$ and $m_1=q_2$ or $a=q_2$ and $m_1=q_1$.
Let us fix the choice that $m_1=q_1$ and $a=q_2>1$.

Then for fixed $p>1$ the numbers $q_1, q_2$ with (\ref{lens}) determine $K^{\frac{p}{q}}$ and that $m_1=q_1$ mod $p$. Then since $K^{\frac{p}{q}}$ is also determined by $m_1=q_1, q$ and that $q$ and $q_2$ are both the longitude variables for the construction of $K^{\frac{p}{q}}$ we have that $q=q_2$ mod $p$.

Thus we have the following theorem:
\begin{theorem}
The representation (\ref{qq1}) can also represent all the lens spaces when $m_2$ is related to $m_1 $ such that $m_2=0$ where we let the lens space $S^2\times S^1$ be represented by (\ref{qq1}) with $p=m_1=m_2=0$ and for all other 3-manifold $M$ we let $m_1 \neq 0$. 

Then for $p>1$ we have that the lens space $M$ constructed by $C^{\frac{p}{q^{\prime}}}$ is uniquely (in the sense of mod $p$) represented by the following invariant:
\begin{equation}
\overline{W}(K^{\frac{-p}{q}}):=R^{2p}R_1^{-q_1q_2}W(C_1)W(C_2) \label{qq2b}
\end{equation}
where $1\leq q_1\leq p-1$ and $q^{\prime}=q_2+n_2p$ for some integer $n_2$ and $q_2$ is restricted to a region mod $p$ such that $q_1<q_2$; and the nontrivial knot $K$ is indexed by the integer $m_1$ where $am_1=q_1q_2$ for some integer $a$ such that $a=q_2>1$ and 
$m_1=q_1$;
and $q=q_2$ mod $p$.
\end{theorem}

{\bf Remark}. We do not count $S^3$ as a lens space and that $S^3$ is simply represented by the constant $1$. $\diamond$

{\bf Remark}. In the above representation we choose $q_1$ such that $1\leq q_1\leq p-1$. It is clear that we may choose other regions (mod $p$) for $q_1$.
$\diamond$

{\bf Remark}. We may write (\ref{qq2b}) in the form
\begin{equation}
\overline{W}(K^{\frac{-p}{q}}):=R^{2p}W(C_1)R_2^{-am_1}W(C_2) 
\label{qq2c}
\end{equation}
or simply in the form $R^{2p}R_2^{-am_1}$. This is a degenerate form of the general form of (\ref{qq1}) with the degeneration that the variable $R_1^{-m_1}$ does not appear. $\diamond$

 Let us then consider a
3-manifold $M$ which is obtained from a framed link $L$ with the
minimal number $n$ of component knots where $n\geq 2$. From the
second Kirby moves we may suppose that $L$ is in the form that the
components $K_i^{\frac{p_i}{q_i}},i=1,...,n$ of $L$ do not wind
each other in the form described by the second Kirby moves (We remark that the usual second Kirby move is for framed links with integral framings. In the following lemma we may generalize it to framed links with rational framings). Let us
say that this minimal $L$ is in the form of maximal non-degenerate
state where the degenerate property is from the winding of one
component knot with the other component knot by the second Kirby
moves. Thus this $L$ has both the minimal and maximal properties as
described. Then we want to find a one-to-one representation (or
invariant) of $M$ from this $L$. Before this let us first prove the following lemma which extends the Kirby theorem:
\begin{lemma}
Let a move generalize the usual second Kirby
move to framed links with rational framings. Then we have that any homeomorphism on a 3-manifold can be written as a sequence of first Kirby move and this generalized second Kirby move.
\end{lemma}
{\bf Proof}. Let two framed knots  $K_i^{\frac{p_i}{q_i}},i=1,2$ be
with (coprime) rational framings (for generalzing the second Kirby move from integral framings to rational framings). 
Let the quantum invariant form of these two framed knots be given by:
\begin{equation}
\overline{W}(K_i^{\frac{p_i}{q_i}}):=R_i^{2p_i}R_{i1}^{-m_{i1}}R_{i2}^{-a_im_{i1}}W(C_{i1})W(C_{i2})
\label{kirby1}
\end{equation}
for $i=1,2$ where $R_i$, $R_{i1}$, $R_{i2}$ for $i=1,2$ are independent (The term $R_{i1}^{-m_{i1}}$ disappears when $K_i^{\frac{p_i}{q_i}}$ is for a lens space).
In terms of these two quantum invariant forms the (generalized) second Kirby move is described by changing $\overline{W}(K_1^{\frac{p_1}{q_1}})$ to the following form:
\begin{equation}
\overline{W}(K_1^{\prime\frac{p_1+p_2}{q_1+q_2-1}}):=R_1^{2(p_1 +p_2 )}
R_{11}^{-m^{\prime}}R_{12}^{-a^{\prime}m^{\prime}}W(C_{11})W(C_{12})
\label{kirby2}
\end{equation}
where the knot $K_1^{\prime}$ is obtained by winding $K_1$ to $K_2$ (by the connected sum operation) as described by the usual second Kirby move 
\cite{Kir}\cite{Pra}; $m^{\prime}$ denotes the assigned integer of $K_1^{\prime}$ and $a^{\prime}$ is the number
correspoding  to the  
number $q_1+q_2-1$;
and from the winding of $K_1$ and its partner 
to $K_2$ and its partner respectively
 we have the degeneration that  
in (\ref{kirby2}) $R_{11}=R_{21}$ and $R_{12}=R_{22}$; and the linking number between
$K_1^{\prime\frac{p_1+p_2}{q_1+q_2-1}}$ 
and $K_2^\frac{p_2}{q_2}$ is determined by the winding of $K_1$ to $K_2$.

Then from (\ref{kirby1}) for $i=1,2$ we can construct (\ref{kirby2}) for
$K_1^{\prime\frac{p_1+p_2}{q_1+q_2-1}}$.
 Conversely (by the degeneration 
 $R_{11}=R_{21}$ and $R_{12}=R_{22}$) from (\ref{kirby2}) for
 $K_1^{\prime\frac{p_1+p_2}{q_1+q_2-1}}$
 and (\ref{kirby1}) for $K_2^{\frac{p_2}{q_2}}$ we can reconstruct the data for $K_1^{\frac{p_1}{q_1}}$ and thus the quantum invariant (\ref{kirby1}) for $K_1^{\frac{p_1}{q_1}}$. Thus these two representations are equivalent.

Then the degeneration  
$R_{11}=R_{21}$ and $R_{12}=R_{22}$ of this winding of generalized second Kirby move gives a symmetry (or homeomorphism).
Conversely this winding  
is the only way to get a degenerate form which is equivalent to the two quantum invariant forms (\ref{kirby1}).
 Thus this winding of generalized second Kirby move is the only source for introducing symmetry relating two framed knots. Thus any homeomorphism on a 3-manifold can be written as a sequence of first Kirby move and this generalized second Kirby move. This proves the lemma. $\diamond$

For simplicity let us call this  generalized second Kirby move  
as the second Kirby move.
Then we want to find a one-to-one representation (or
invariant) of $M$ from the given $L$.
Let us write $W(L)$ in the form:
\begin{equation}
W(L)=P_L \prod_i W(K_i^{\frac{p_i}{q_i}}) \label{qq2a}
\end{equation}
where $P_L$ denotes a product of $R$-matrices acting on a subset
of $\{W(K_i),W(K_{ic}), i =1,...,n \}$ where $W(K_i^{\frac{p_i}{q_i}})$ are
independent (This is from the form of $L$ that the component knots
$K_i$ are independent in the sense that they do not wind each
other by the second Kirby moves). Then we consider the following representation (or
invariant) of $M$:
\begin{equation}
\overline{W}(L):=P_L \prod_i
\overline{W}(K_i^{\frac{p_i}{q_i}}) \label{qq2}
\end{equation}
where we define $\overline{W}(K_i^{\frac{p_i}{q_i}})$ by
(\ref{qq1}) and they are independent. We have the following
theorem:

\begin{theorem}
Let $M$ be a closed (oriented and connected) 3-manifold which is
constructed by a Dehn surgery on a framed $L$ with the minimal
number $n$ of component knots ($L$ has both the minimal and maximal properties). Then we have that (\ref{qq2}) is a one-to-one representation (or invariant) of $M$.
\end{theorem}

{\bf Proof}. We want to show
that (\ref{qq2}) is a one-to-one representation (or invariant) of
$M$. Let $L^{\prime}$ be another framed link for $M$ which is also
with the minimal number $n$ (and with the maximal property). Then
we want to show $\overline{W}(L)=\overline{W}(L^{\prime})$.
Suppose that each component $\overline{W}(K_i^{\frac{p_i}{q_i}})$ of $\overline{W}(L)$ does not represent a lens space. Then these components $\overline{W}(K_i^{\frac{p_i}{q_i}})$ are invariants of the components of $L$ respectively. 
Then since the components of $L$ do not wind each other as described by the second Kirby move we have that the components of $L$ are independent of each other. Thus there is no nontrivial homeomorphism changing these components $\overline{W}(K_i^{\frac{p_i}{q_i}})$  except those homeomorphisms involving 
the second Kirby moves for the winding of the components of $L$ with each other. 
Then under the second Kirby moves of these homeomorphisms we have that the components of $L$ wind each other and thus will reduce the independent degree of freedom to be less than $n$. Thus to restore the degree of freedom to $n$ these homeomorphisms must also contain the first Kirby moves of adding unknots with framing $\pm 1$. In this case these unknots can be deleted and thus $L$ 
 is not minimal and this is a contradiction. Thus there is no nontrivial homeomorphism changing the components $\overline{W}(K_i^{\frac{p_i}{q_i}})$ of $\overline{W}(L)$ 
 except those homeomorphisms consist of only the second Kirby moves for the winding of the components of $L$ with each other. 

 Now suppose that $\overline{W}(L)\neq \overline{W}(L^{\prime})$. Then there exists nontrivial homeomorphism of changing $L$ to $L^{\prime}$ for changing the components $\overline{W}(K_i^{\frac{p_i}{q_i}})$ of $\overline{W}(L)$ to the components of 
 $\overline{W}(L^{\prime})$. This is impossible since there are no  nontrivial homeomorphsm for changing these components $\overline{W}(K_i^{\frac{p_i}{q_i}})$ except those homeomorphisms consist of only the second Kirby moves for the winding of the components of $L$ with each other. Thus $\overline{W}(L)=\overline{W}(L^{\prime})$.
 
 Then let us suppose that there exists a component $\overline{W}(K_i^{\frac{p_i}{q_i}})$ of $\overline{W}(L)$ representing a lens space. Then this component $K_i^{\frac{p_i}{q_i}}$ must not be linked with the other components of $L$. Suppose not. Then by 
 the Rolfsen twist on this component $K_i^{\frac{p_i}{q_i}}$ such linking changes this component and the components linking to this component and thus the Rolfsen twist is a nontrivial homeomorphism on $L$ and thus by the above lemma it must contains a first Kirby move of adding a framed unknot $C^{\pm 1}$ with framing $\pm 1$. Then this framed unknot $C^{\pm 1}$ can be deleted by the first Kirby move and thus $L$ is not minimal. This is a contradiction. 
 
 Thus $L$ must be in the form that it is the sum of two parts where one part is only with components which do not represent  lens spaces and are of maximal nondegenerate form and the other part is formed by the components of $L$ representing lens spaces and each component of $L$ whenever representing a lens space must be unlinked with the other components of $L$.
Further since $L$ is minimal these unlinked components of $L$ can not be combined with each other to form another minimal representation. Thus the part of $\overline{W}(L)$ formed by the components of $\overline{W}(L)$ representing lens spaces is unique. On the other hand we have also shown that the part of $\overline{W}(L)$ with only components which do not represent  lens spaces and are of maximal nondegenerate form is unique. Thus we have that $\overline{W}(L)$ is unique and $\overline{W}(L)=\overline{W}(L^{\prime})$.
 
In conclusion we have that (\ref{qq2}) is a one-to-one representation (or invariant) of $M$, as was to be proved. 
$\diamond$

As a converse to the above theorem let us suppose that  the
representation (\ref{qq2}) uniquely represents $M_L$ in the sense
that there are no nontrivial homeomorphism transforming the $n$
independent components of $\overline{W}(L)$ to other $n$
independent components of $\overline{W}(L^{\prime})$ where the
link $L^{\prime}$ also gives the manifold $M_L$. Then from the
above proof we see that the link $L$ is a minimal (and maximal) link for
obtaining $M_L$.

{\bf Remark}. Let $L$ be a minimal (and maximal) framed link. Then from the above proof we have that the components of $L$ are independent of each other in the sense that if we transform a component framed knot of $L$ to an equivalent framed knot by a homeomorphism then the other components of $L$ are not affected by this transformation. $\diamond$

From the above theorems we then have the following
classification theorem:
\begin{theorem}
Let $M$ be a closed (oriented and connected) 3-manifold which is
not homeomorphic to $S^3$. Then  the representation 
consists of (\ref{qq2a}), (\ref{qq1}), (\ref{qq2b}) (or (\ref{qq2c})) 
is a one-to-one invariant of $M$. 
This quantum
invariant of $M$ has the following general expression (which is the representation (\ref{qq2})):
\begin{equation}
\overline{W}(L)=
P_L\prod_{i=1}^n\overline{W}(K_i^{\frac{p_i}{q_i}})
 \label{Pr1}
\end{equation}
where $L$ denotes a minimal surgery link for $M$ and $n\geq 1$
is the minimal number for $M$ (and $L$ is with the maximal property).

For $M=S^3$ we have that the invariant for $M$ is $1$.
\end{theorem}

\section{Proof of Poincar\'{e} Conjecture}\label{sec15}

Let us apply the above classification of closed 3-manifolds to prove
the Poincar\'{e} conjecture.
Let $M$ be a closed 3-manifold obtained from surgery on a framed
nontrivial knot $K^{\frac pq}$ which is the minimal link for $M$ with minimal
number $n=1$. From the above section we have the one-to-one generalized Wilson loop 
representation $\overline{W}(K_i^{\frac{p_i}{q_i}})$ as invariant of $M$ (From this we have that $K^{\frac pq}$ is a minimal link for $M$ with minimal
number $n=1$ and when $M$ is a lens space we also use this invariant of $M$ to represent $M$). Let us from this invariant to show that $M$ is non-simply connected.
 
To this end let us consider the fundamental group of $M$. We
recall that the fundamental group of $M$ can be obtained from the
knot group $G$ of $K$ by adding relations to the generators of $G$
where these additional relations are from the partner knot $K_c$
of $K$. By these additional relations we have that the fundamental
group of $M$ is formed as a subgroup of $G$. As an example we have
that the fundamental group of the Poincare sphere $M$ is given by
$\pi_1(M)=\{x,y,z| xy=yz=zx, [K_c]= x^{-2}yxz=1\}$ where $x,y,z$
are generators of the knot group $G$ of the right trefoil knot $K$
and $[K_c]= x^{-2}yxz=1$ is the additional relation. We have that
the generators of $G$ are distinguished by the crossings of the
knot $K$. Then we have that $K$ is represented by the generalized
Wilson loop $W(K)=R_1^{-m_1}W(C_1)$ which is in a form that the
crossings of $K$ have been equivalently eliminated such that $K$
is represented by a circle $C_1$ which winds with additional $m_1$
times by the factor $R_1^{-m_1}$ (Similarly we have that
$K^{\frac{p}{q}}$ represented by (\ref{qq1}) is in a form that the
crossings have been equivalently eliminated). Now let $x$ be a
generator of $G$ of $K$. Then we have that $x$ is represented as a
generator of a knot group of $C_1$ in this representation
$W(K)=R_1^{-m_1}W(C_1)$ of $K$. In this representation we have
that $x$ is represented as a circle encircling $C_1$. Let us write
$W(C_1)=R_1^{-n_1}A$ for some variable integer $n_1$ which is a
form that $C_1$ winds $n_1$ times. Thus the total winding is
$m_1+n_1$ times. Then as an equivalence we may regard $C_1$ winds
one time and $x$ is represented as a circle encircles $C_1$ with
$m_1+n_1$ times. Now while the generators $x$ of $G$ of $K$ are
distinguished by the crossings of $K$ we have that in the
representation $W(K)=R_1^{-m_1}W(C_1)$ of $K$ these $x$ are not
distinguished when they are generators for $C_1$ since $C_1$ has
no crossings.

Similarly for a generator $y$ of the knot group of $K_{cr}$
represented by  
$W(K_{cr})=R^{-am_1}W(C_2)=R^{-am_1}R^{-n_2}A$
we regard
$y$ as a circle encircles $C_2$ with $am_1+n_2$ times while $C_2$
winds one time.

Now in the generalized Wilson loop representation
$\overline{W}(K^{\frac pq})=R^{2p}R_1^{-m_1}R_2^{-am_1}W(C_1)W(C_2)$
of $K^{\frac pq}$ we have
that all the crossings of $K^{\frac pq}$ are eliminated (When $K^{\frac pq}$ is for a lens space we have the form (\ref{qq2c}) that the factor $R_1^{-m_1}$ disappears). This representation is similar to the representation $W(C_1)W(C_2)$ of the manifold $S^2\times S^1$ that all the crossings of $C_1$, $C_2$ and between $C_1$ and $C_2$ are eliminated. Thus as similar to the case of the manifold $S^2\times S^1$ we have that in this representation
the additional relations among the generators $x$ of $G$ for $K$ for defining the
fundamental group $\pi_1(M)$ of $M$ from $G$ are 
equivalently
eliminated and 
equivalently transformed to a relation
from which  the generators of $K$ are related to generators of
$K_{c}$ such that the generators of the fundamental group of $M$
are formed from the generators of $K$. Let us determine this
relation in this Wilson loop representation, as
follows. Since $\pi_1(M)$ is a subgroup of $G$ we have that a
generator $g$ of $\pi_1(M)$ in this generalized Wilson loop
representation is of the form that $g$ is a multiple product of
$x$ that $g$ is a circle encircles $C_1$ with $r_1(m_1+n_1)$ times
for some integer $r_1$. Then we have that the additional relation
gives a relation of $g$ to the generators $y$ of $K_{c}$. Now in
this generalized Wilson loop representation we have that the only
way that $g$ is related to the generators $y$ of $K_{c}$ is that
$g$ is also a multiple product of $y$ that $g$ is a circle
encircles $C_2$ with 
$r_2(am_1+n_2)$ 
times for some integer
$r_2$. Then since  
$am_1$ contains the factor $m_1$ we may choose
the integer $n_2=an_1$ such that $m_1+n_1$ is a factor of 
$am_1+n_2$
that $am_1+n_2=a(m_1+n_1)$.  
It follows
that we have the existence of $g$ of the form that $g$ is as a
circle encircles $C_1$ and $C_2$ with 
$a(m_1+n_1)$ times.

Now since $G$ is a nontrivial group with nontrivial generators $x$
we have that the fundamental group $\pi_1(M)$ of $M$ is with the
existence of nontrivial generators $g$ which in the generalized
Wilson loop representation $\overline{W}(K^{\frac pq})$ are just
some circles encircle $C_1$ and $C_2$ with 
$a(m_1+n_1)$ times. This shows that $M$ is non
simply-connected.
Thus we have proved
the following property P conjecture:
\begin{theorem}(Property P Conjecture)
Let $K$ be a nontrivial knot. Then the 3-manifold $M$ obtained
from Dehn surgery on $K^{\frac pq}$ is non simply-connected.
\end{theorem}

{\bf Remark}. When $K^{\frac pq}$ is for a lens space and $p=0$ we have $M=S^2\times S^1$. Then we have that $g$ is a
circle encircles $C_1$ and $C_2$ for $n_2$ times.  
Then any element of
$\pi_1(M)$ is of the form $g^k$ which is a circle encircles $C_1$
and $C_2$ for $kn_2$ times. Thus $\pi_1(M)$ is the group $ Z$ of
integers. $\diamond$

{\bf Remark}. When $K^{\frac pq}$ is for a lens space $M$ and $p>1$ then  the generator $g$ of the fundamental group $\pi_1(M)$ is a
circle encircles $C_1$ and $C_2$ for
$-q_1q_2+n_2=1+n^{'}p+n_2=a(m_1+n_1)$ times for $n_2=an_1$. Then  we have that
$g^p$ is a circle encircles $C_1$ and $C_2$ for
$p+n^{'}p^2+n_2 p$ times which is of the form
$0+n^{''}p+n_2 p=0$ mod $p$. Thus $g^p$ can be identified as
the identity $e$ of a quotient group which is a cyclic group with
$p$ elements and with $g$ as the generator. $\diamond$

{\bf Remark}.
For the Poincar\'{e} sphere $M_P$ obtained from surgery on $K_{RT}^1$ we have that
$am_1=-1$ 
where $m_1=1$ is the index for the right trefoil knot $K_{RT}$.
Thus the knot group of $K_{RT}$ and the fundamental group
$\pi_1(M_P)$ of $M_P$ are with the generators $g$ and $g^{-1}$ respectively and thus can be
with the same generator $g$. Then since $\pi_1(M_P)$ is a proper subgroup
of the knot group of $K_{RT}$ we have that the representation of the fundamental group
$\pi_1(M_P)$ must be a quotient group with finite elements of the representation
of the knot group of $K_{RT}$ which is an infinite cyclic
group generated by $g$. We have that this quotient group is with element of the
form $g^{k}$ with windings $k + [n_2]$ where we choose $n_2=120n_3$ for some
integer variable $n_3$. Then when $k=120$ we have $k=0$ mod $120$
and $g^{120}=e$ where $120$ is the number of elements of  $\pi_1(M_P)$.

Similarly for a 3-manifold $M$ obtained from $K_{RT}^{p}$ for $p>1$ (The Poincar\'{e} sphere $M_P$ is with $p=1$) we have that $a= -1$ and thus the fundamental group $\pi_1(M)$ is a finite nontrivial group.

On the other hand when a 3-manifold $M$ obtained from $K^{\frac
pq}$ (which is not a lens space) and is not homeomorphic to 3-manifolds  obtained from
$K_{RT}^{p}$ for $p\geq 1$ we have that $am_1\neq \pm m_1$. Thus
in this case we have that the generator $g$ of the representation
of $\pi_1(M)$ of $M$ is only a subgenerator of the representation
of the knot group of $K$. Thus we have that this representation of
$\pi_1(M)$ which is generated by $g$ is already a proper cyclic
subgroup of the representation of the knot group of $K$. In this
case if there are no further conditions on this representation of
$\pi_1(M)$ to be a quotient subgroup of the representation of the
knot group of $K$ then we have that this representation of
$\pi_1(M)$ is an infinite cyclic group. From this we then have
that the fundamental group $\pi_1(M)$ of $M$ is an infinite group.
$\diamond$

Now let $M$ be a closed 3-manifold which is classified by the
following  minimal invariant (\ref{Pr1}):
\begin{equation}
\overline{W}(L)=
P_L\prod_{i=1}^n\overline{W}(K_i^{\frac{p_i}{q_i}})
 \label{Pr2}
\end{equation}
where $L$ denotes a surgery link for $M$ and $n\geq 2$ is the
minimal number for $M$. From the above section we have that this
is a one-to-one invariant for $M$ that the framed knot components
$K_i^{\frac{p_i}{q_i}}, i=1,...,n$ cannot be eliminated.

Then from this invariant (or representation) we have that the
framed knot components $K_i^{\frac{p_i}{q_i}}, i=1,...,n$ are
independent of each other in the sense that their forms are not
changed by each other  though they are linked together to form the
linked $L$. Thus we have that the Wilson loop representation of
the fundamental group of $M$ contains the Wilson loop
representation of the fundamental groups of the manifolds
constructed from the framed knot components
$K_i^{\frac{p_i}{q_i}}, i=1,...,n$.

Then since the Wilson loop representation of the fundamental group
of $M$ contains the Wilson loop representation of the
fundamental groups of the framed knot components
$K_i^{\frac{p_i}{q_i}}, i=1,...,n$ and these fundamental groups of
the framed knot components $K_i^{\frac{p_i}{q_i}}$ are nontrivial
we have that the fundamental group of $M$ is nontrivial and $M$ is
non-simply connected.

Now let $M$ be a simply connected closed (orientable and
connected) 3-manifold. We want to show that it is homeomorphic to
$S^3$. Let us suppose that $M$ is not homeomorphic to $S^3$. Then
from the above classification theorem we have that $M$ is
classified by a quantum invariant of the form (\ref{Pr2}) for
$n\geq 1$. Thus we have that the fundamental group of $M$ is
nontrivial and thus $M$ is not simply connected. This is a
contradiction. Thus  $M$ must be homeomorphic to $S^3$, as was to
be proved.  This proves the following Poincar\'{e} Conjecture:

\begin{theorem}[Poincar\'{e} Conjecture]
Let $M$ be a closed (orientable and connected) and
simply connected 3-manifold. Then $M$ is homeomorphic
to the 3-sphere $S^3$.
\end{theorem}

\section{Conclusion}\label{sec16}

In this paper from a quantum gauge model we derive a conformal field theory
structure from which we derive
 a knot invariant
related to the HOMFLY polynomial.
The relation between these two
invariants is that both  the HOMFLY  polynomial
and this  knot  invariant can be
derived by using two Knizhnik-Zamolodchikov (KZ) equations
 which are dual to
each other and are derived from the quantum gauge model. In this
derivation an important concept called the Wilson lines and the
generalized Wilson loops is introduced such that each knot diagram
is represented by a generalized Wilson loop where the upper
crossing, zero crossing and undercrossing of two curves can be
represented by the orderings of two Wilson lines represent these
two curves. We show that this invariant can  classify knots by
showing that the generalized Wilson loop of a knot  faithfully
represents the knot in the sense that if two knot diagrams have
the same generalized Wilson loop then these two knot diagrams must
be equivalent. This invariant is in terms of the monodromy $R$
of the two Knizhnik-Zamolochikov equations. In the case of knots
this invariant can be written in the form $Tr R^{-m}$. From
this invariant we may classify knots with the integer $m$. A
classification table of knots can then be formed where prime knots
are classified with odd prime numbers $m$ and non-prime knots are
classified with non-prime numbers $m$.

Then from the quantum link invariant we can construct quantum
invariant of 3-manifolds . We first construct quantum
invariant of closed three-manifolds obtained by Dehn surgery on
framed knots.  We then introduce the concept of minimal link to
construct quantum invariant of closed three-manifolds obtained by
Dehn surgery on framed links. Then by using the Lickorish-Wallace theorem we show that this quantum invariant of 3-manifolds gives a one-to-one classification of closed 3-manifolds.  From this
classification of closed 3-manifolds we can then prove the
Poincar\'{e} Conjecture.

\end{document}